\documentclass[reqno,11pt]{amsart}
\usepackage{amsmath}
\usepackage{amsfonts,amsthm,amssymb}
\usepackage{bm}
\usepackage{amscd}
\usepackage{t1enc}
\usepackage[mathscr]{eucal}
\usepackage{indentfirst}
\usepackage{graphicx}
\usepackage{graphics}
\usepackage{pict2e}
\usepackage{epic}
\usepackage{comment}
\numberwithin{equation}{section}
\usepackage{geometry}
\usepackage{fullpage}
\usepackage{epstopdf} 
\usepackage[colorlinks,linkcolor=blue]{hyperref}
\usepackage[capitalise,noabbrev]{cleveref}

\crefformat{equation}{(#2#1#3)}
\crefmultiformat{equation}{(#2#1#3)}{ and~(#2#1#3)}{, (#2#1#3)}{ and~(#2#1#3)}
\crefrangeformat{equation}{(#3#1#4) to~(#5#2#6)}

\usepackage[normalem]{ulem}

\allowdisplaybreaks

\theoremstyle{plain}
\newtheorem{Thm}{Theorem}[section]
\newtheorem*{Thm*}{Theorem}
\newtheorem{Lem}[Thm]{Lemma}

\newtheorem{Prop}[Thm]{Proposition}

\theoremstyle{definition}

\newtheorem{Rem}[Thm]{Remark}
\newtheorem{?}[Thm]{Problem}

\newcommand{\p}{\partial}
\newcommand{\R}{\mathbb{R}}

\newcommand{\M}{\mathbf{M}}

\newcommand{\D}{\mathbf{D}}\newcommand{\Do}{\mathbf{D}_0}\newcommand{\Dn}{\mathbf{D}_{\neq}}
\newcommand{\G}{\mathbf{G}}

\newcommand{\Torus}{\mathbb{T}}

\newcommand{\dv}{\text{div}}

\newcommand{\mr}{\mathring}

\newcommand{\rhot}{\tilde{\rho}}

\newcommand{\rhob}{\bar{\rho}}

\newcommand{\ut}{\tilde{u}}

\newcommand{\uf}{{u}}

\newcommand{\thetat}{\tilde{\theta}}
\newcommand{\rhom}{\mathring{\rho}}
\newcommand{\um}{\mathring{u}}

\newcommand{\phim}{\mathring{\phi}}

\newcommand{\zetam}{\mathring{\zeta}}
\newcommand{\Pb}{\mathbf{P}}
\newcommand{\px}{\partial_1}
\newcommand{\pt}{\partial_t}

\newcommand{\T}{\theta}

\newcommand{\Tm}{\mathring{\theta}}
\newcommand{\Tt}{\tilde{\theta}}

\newcommand{\mm}{\mathring{m}}
\newcommand{\mt}{\tilde{m}}
\newcommand{\mtb}{\mathbf{\tilde{m}}}

\newcommand{\Em}{\mathring{E}}
\newcommand{\Et}{\tilde{E}}
\newcommand{\Eb}{\bar{E}}

\newcommand{\mbb}{\mathbf{\bar{m}}}

\newcommand{\abs}[1]{\left\lvert#1\right\rvert}

\newcommand{\norm}[1]{\left\lVert#1\right\rVert}

\usepackage{orcidlink}
\usepackage{cite}
\newcommand{\LL}{\mathbf{L}}

\begin{document}

\title{Nonlinear Stability of Planar Shock Waves for the 3-D Boltzmann Equation}
	
\author[D.-Q. Deng]{Dingqun Deng
			\orcidlink{0000-0001-9678-314X}
}
\address[D.-Q. Deng]{Department of Mathematics, Pohang University of Science and Technology, Pohang, Republic of Korea (South), ORCID: \href{https://orcid.org/0000-0001-9678-314X}{0000-0001-9678-314X}}
\email{dingqun.deng@postech.ac.kr}\email{dingqun.deng@gmail.com}

	\author[L.-D. Xu]{Lingda Xu \orcidlink{0000-0002-2982-2994}}
\address[L.-D. Xu]{Department of Mathematics, Yau Mathematical Sciences Center, Tsinghua University, Beijing 100084, China, ORCID: \href{https://orcid.org/0000-0002-2982-2994}{0000-0002-2982-2994}}
\email{xulingda@bimsa.cn}		

\date{\today}

\subjclass[2020]{Primary: 35Q20; Secondary 76L05, 76P05, 35L67}

\keywords{3-D Boltzmann equation, planar shock waves, large-time behavior, diffusion waves, micro-macro decomposition}

\begin{abstract}
This paper studies the stability and large-time behavior of the three-dimensional (3-D) Boltzmann equation near shock profiles. We prove the nonlinear stability of the composite wave consisting of two shock profiles under general perturbations without the assumption of integral zero of macroscopic quantities. To address the challenge caused by the compressibility of shock profiles, we apply the method of anti-derivative based on macro-micro decomposition. However, the system of anti-derivatives presents certain difficulties. Firstly, general perturbations may generate diffusion waves that evolve and interact with shock profiles, resulting in errors that are not controllable. We therefore introduce a set of coupled diffusion waves to cancel out these poor errors and perform careful estimates on wave interactions. Secondly, we perform diagonalized system estimates to fully exploit the compressibility of shock profiles and control terms that decay slowly. Thirdly, the presence of diffusion waves causes critical terms with decay $(1+t)^{-1}$, and we introduce a Poincar\'e type of inequality to address these terms. Finally, estimates on anti-derivatives can only control terms along the propagation direction, while for transversal directions, we use the entropy-entropy flux pair as well as the Poincar\'e inequality to control the lower order terms using diffusion terms. As a result, we obtain nonlinear stability through the energy method, which is the first stability result for the planar shock of the multi-dimensional Boltzmann equation to the best of our knowledge.

\end{abstract}

\maketitle
\tableofcontents

\section{Introduction}
In this work, we consider the 3-D Boltzmann equation, which reads
\begin{align}\label{be}
	f_t+v\cdot\nabla_xf=Q(f,f),\ \ \ \ \ f(0,x,v)=f_0(x,v), 
\end{align}
where $f(t, x, v)$ is the particle distribution function at time $t \geqq 0,$ position $x=(x_1,x_2,x_3) \in \mathbb{D}$ with velocity $v=\left(v_{1}, v_{2}, v_{3}\right) \in \mathbb{R}^{3} .$ Here, $\mathbb{D}:=\R\times\Torus^2$ is the infinitely long flat torus and $\Torus:=(\mathbb{R}/\mathbb{Z})$ is the one-dimensional (1-D) periodic domain. The Boltzmann collision operator $Q(\cdot,\cdot)$ is given by
\begin{align}\label{Q}
	Q(f, g)(v)& =\frac{1}{2} \int_{\mathbb{R}^{3} \times \mathbb{S}^{2}} B\left(\left|v-v_*\right|, \vartheta\right) \big\{f\left(v^{\prime}\right) g\left(v_*^{\prime}\right)+f\left(v_*^{\prime}\right) g\left(v^{\prime}\right)  \\
	&\qquad \qquad\qquad\qquad\qquad\qquad -f(v) g\left(v_*\right)-f\left(v_*\right) g(v)\big\} \ \mathrm{d} v_* \mathrm{d} \omega \notag \\
	& =: Q_{+}^{1}(f, g)(v) +Q_{+}^{2}(f, g)(v) +Q_{-}^{1}(f, g)(v) +Q_{-}^{2}(f, g)(v), \notag
\end{align}
where $\mathbb{S}^{2}$ denotes the unit sphere in $\mathbb{R}^{3}$, $f(v)=f(t, x, v)$, and $(v',v_*')$ are given by 
\begin{align}
	v^{\prime}=v-\left[\left(v-v_*\right) \cdot \omega\right] \omega, \quad v_*^{\prime}=v_*+\left[\left(v-v_*\right) \cdot \omega\right] \omega.
\end{align}
The conservation laws of momentum and energy give
\begin{align}
	v'+v_*'=v+v_*,\ \ \ \ \ \ \ \ \left|v'\right|^2+\left|v_*'\right|^2=\left|v\right|^2+\left|v_*\right|^2.
\end{align}
The Boltzmann collision kernel $B=B(v-v_*, \vartheta)$ in \cref{Q} depends only on $|v-v_*|$ and $\vartheta$ with $\cos \vartheta=(v-v_*) \cdot \omega /|v-v_*| .$ 
In this work, We consider the Grad's angular cut-off assumption: 
\begin{align}
	B(v-v_*, \vartheta)=|v-v_*|^{\gamma} b(\vartheta),
\end{align}
with
\[
-3<\gamma \leq 1, \quad 0 \leq b(\vartheta) \le C|\cos \vartheta|.
\]
In this paper, we consider the hard sphere model, that is $\gamma=1$.

\subsection{Literature} Here is a brief introduction to the literature on shock waves, anti-derivatives, and the zero-mass condition.

\smallskip 
It is well known that there are close relations between the Boltzmann equation \cref{be} and the systems of fluid mechanics, such as the compressible Euler system and the Navier-Stokes equations. In fact, through the famous Hilbert expansion and the Chapman-Enskog expansion, one knows that the compressible Euler equations are in the leading order of the Boltzmann equation with respect to the Knudsen number. Note that the compressible Euler system is a typical hyperbolic conservation law for which the important issues are the formation and evolution of shock waves. For the viscous conservation law, for example, the compressible Navier-Stokes equations, the shock wave is smoothed to the traveling wave solution due to the dissipation effect. 

\smallskip 
We first review some notable work in one-dimensional viscous conservation laws and the important case, Navier-Stokes equations. In the viscous case, \cite{Ilin1960} proved the stability of the shock wave for the scalar conservation law with a comparison principle. This result has been extended to cases of the system in \cite{G,MN} by introducing the method of anti-derivatives, that is, denoting
\begin{align}
    \Phi(x_1,t)=\int_{-\infty}^{x_1}\phi(y,t)dy,
\end{align}
For the perturbation $\phi.$ This method is very strong in the study of shock waves, but requires a zero-mass condition for initial perturbations, i.e.
\begin{align}
    \int_{-\infty}^{\infty}\phi(x,0)dx=0.
\end{align}
Note that this condition is not physical in some senses. Since then, many efforts have been made by excellent mathematicians and a lot of remarkable results have been achieved, now we have a deeper understanding of the role of the zero-mass condition. For the problems with general perturbations, the position of the shock wave shifts and diffusion waves are generated, which makes the analysis difficult, cf. \cite{Liu1}. For the uniformly parabolic systems, with partial construction of the fundamental solution, \cite{SX} extended the result of \cite{Liu1} from non-zero mass perturbations with stringent conditions to general ones. 

\smallskip 
The non-zero mass stability is also obtained by the pointwise approach initiated in \cite{Liu1997}, by the construction of Green's function, a $priori$ estimates via Duhamel's Principle, and estimates on wave interactions. And \cite{LZ3} used a significantly different method to construct an approximate Green's function, and their approach was extended to systems with physical viscosity, including compressible Navier-Stokes equations and Magnetohydrodynamics equations, cf. \cite{LZ5}. There are interesting attempts to prove nonlinear stability from linear and spectral stability \cite{MZ2004}. For this approach, we refer to \cite{HLZ2017} and the references therein for the spectral stability. For nonlinear stability derived from spectral stability, we refer to \cite{MZ2004,HRZ2006}. 

\smallskip 
However, for non-zero mass stability of the Navier-Stokes equations by the energy method, due to the diffusion wave, the error terms do not have a sufficient decay rate by the anti-derivative approach. For the non-isentropic case, we refer to \cite{HM}, where the authors construct diffusion waves with a good observation that at the hyperbolic level, a characteristic variable associated with the second characteristic field is decoupled from other characteristic ones up to the second order. They solved this problem in the cases of two shock waves with small strengths of the same order. Recently, for the isentropic Navier-Stokes equations, \cite{WW2022} studied the large-time behavior of shock waves under general $H^2(\R\times\Torus^2)$ perturbations by the a-contraction method with time-dependent shift introduced, where $\Torus:=\R/\mathbb{Z}$. However, it is very difficult to give a good estimate for the shift introduced in their paper, and the position of the shock wave was not identified. According to \cite{P}, one cannot expect stability results under general $H^2(\R\times\Torus^2)$ perturbations and unbounded shifts. Thus, it is very interesting to find such a time-dependent shift to make the profile stable. But in this sense, the shift is essential in their proof, and we would like to use a new technique, different from \cite{WW2022}, to prove the non-zero mass stability of the shock wave at a determined position.

\smallskip 
Although the stability of the shock wave for the compressible Navier-Stokes equations has been extensively studied, the large time stability of the shock wave for the 1-D Boltzmann equation has been open for quite a while, until recently when Liu and Yu \cite{LYs} studied the positivity of the single shock profile to the Boltzmann equation and its nonlinear stability with zero total macroscopic mass condition. Liu and Yu \cite{LYs} and Liu, Yang, and Yu \cite{LYY} introduced the macro-micro decomposition for the Boltzmann equation, which can apply the energy method to study the time-asymptotic stability of nonlinear waves to the Boltzmann equation. Then Liu et al. \cite{LYYZ} proved the time-asymptotic stability of rarefaction waves for the Boltzmann equation. For the viscous contact wave, \cite{HY} obtained its stability with the zero-mass condition, and this result was extended by Huang, Xin, and Yang \cite{HXY} to the perturbations without the zero-mass condition. without the zero-mass condition. Yu \cite{Yu} proved the stability of the single shock profile by the pointwise approach based on the Green function under the general initial perturbation without the total macroscopic zero-mass condition. \cite{WW} proved the stability of the superposition of two shock waves without the zero mass condition. For the fluid dynamic limit of the Riemann solutions of the Euler equations, we refer to \cite{HWY,XZ2,Yu1} for the case of a single wave pattern. For the superposition of multiple wave patterns, we also refer to \cite{HWY1,HWWY}, which justified the limit of the Boltzmann equation to the Euler equations for Riemann problems; their result covers the case of superposition of rarefaction wave, contact wave, and shock wave. We also refer to \cite{DYY1,DYY2,DL} for some interesting wave phenomena in kinetic models.
For the multi-dimensional case, \cite{WW1} proved the stability of the planar rarefaction wave for the 3-D Boltzmann equation and pointed out that the stability of the planar shock wave is completely open.

\smallskip 
For Navier-Stokes equations and viscous conservation laws. There is also a very large body of excellent work attempting to demonstrate the nonlinear stability of multi-dimensional shock waves but with far fewer results compared to the 1-D case. For the viscous shock, by constructing Green's function, the nonlinear stability of planar viscous shocks has been proved by \cite{GM1999,HZ2000,LY2018} in the scalar case with detailed pointwise estimates; in particular, richer wave patterns including Rayleigh-type waves are observed in \cite{LY2018}. We also refer to \cite{KVW2019} for the interesting result about $L^2$ contraction of large planar shock waves. For the cases of systems there are also results deriving nonlinear stability from spectral stability, see for example \cite{HLZ2017} and the references therein.  It is worth mentioning that the nonlinear stability of planar viscous shocks is verified by numerically well-conditioned and analytically justified computations in \cite{HLZ2017} for Navier-Stokes equations, and \cite{FS1,FS2} proved the spectral stability of small amplitude viscous shocks for systems with artificial viscosity. Recently, some interesting progress in rigorous proofs of nonlinear stability for shock waves of Navier-Stokes equations has been made by \cite{Yuan2021,Yuan,WW2022}. As explained earlier, the introduction of a time-dependent shift in \cite{WW2022} makes the location of the shock profile unclear, which is crucial in their theory. The stability of shock waves under periodic perturbations is proved in \cite{Yuan} by introducing the multi-dimensional anti-derivatives, which is the first result to apply this method in the multi-dimensional case to the best of our knowledge. However, a zero-mass type condition is essentially required in \cite{Yuan}. Moreover, how to remove this condition is proposed as an open problem in their paper. 

\subsection{Goal and idea of the proof}
The main result of this paper is to give rigorous proof for nonzero mass stability of the superposition of two shock waves with determined shifts for the 3-D Boltzmann equations under general perturbations. 

\medskip 
Next, we present our strategy of proof. First, based on the result of \cite{LYs}, we construct the superposition of two planar shock profiles, which are traveling wave solutions of the Boltzmann equation \cref{be}. The macro-micro decomposition introduced in \cite{LYY} is also applied. Second, to apply the anti-derivative technique, the initial macroscopic mass should be zero, we introduce the shifts in both a shock wave and a diffusion wave to carry the excess mass as in \cite{Liu1}. However, the presence of a diffusion wave results in bad errors that increase the $L^2$ energy of the antiderivatives. Therefore, a higher order correction is introduced to improve the bad errors inspired by \cite{HM,WW}. Thus, the proper ansatz is composed of two shifted viscous shock waves, a diffusion wave and higher-order derivatives.

\smallskip 
For the a priori estimate, we perform the energy estimates for anti-derivatives. Since the ansatz consists of two shock waves, a diffusion wave, and higher-order corrections, we should study the wave interaction carefully. Moreover, due to the introduction of the diffusion wave and the complexity of the viscosity, some terms cannot be controlled by the compressibility, see \cref{be2}. We further study the diagonalized system to make full use of the compressibility of the shock profile, but some critical terms appear. To deal with these terms, we introduce a Poincar\'e type inequality based on a key cancellation and an estimate of the heat kernel introduced by \cite{HLM}, see \cref{pte}. Note that the anti-derivative technique can only be used to estimate the terms along the shock propagation direction. For the transverse direction, we decompose the perturbation into zero and non-zero modes, the zero mode can be controlled by the energy estimate of anti-derivatives and higher order derivatives, and the Poincar\'e inequality is available for non-zero modes, so they can be controlled by higher order derivatives. Finally, estimates on the microscopic part and an estimate on the highest order of derivatives by the original equation \cref{be} help us to close the a priori estimate. The global existence and stability of the solution of the Boltzmann equation can be derived from a local existence result, the a priori estimates, and a standard continuity argument. 

\smallskip 
The rest of this paper is organized as follows. In \cref{Sec2}, we introduce the construction of the ansatz and the main theorem. In \cref{Sec3}, we formulate the problem, and in \cref{Sec4}, the a priori estimates are performed.

\section{Ansatz and main result}\label{Sec2}
In this Section, we focus on the construction of the ansatz and the mathematical description of the main theorem. Firstly, we apply the macro-micro decomposition introduced in \cite{LYY,LYs}. For the macroscopic part, we construct two viscous shock waves that are traveling wave solutions to the Boltzmann equation. Secondly, to apply the anti-derivative technique, the initial macroscopic mass is required to be zero. We determine the shifts of two shock profiles and introduce the diffusion wave and the coupled diffusion wave to carry the excess initial mass. Then we introduce a decomposition that decomposes the required quantities into zero and non-zero modes in Fourier space and give some properties of the decomposition. Finally, we state the main theorem. 

\subsection{Notations}
We begin with the basic notations. For convenience, let $\nabla=\nabla_x$ be the gradient in $x$, $\operatorname{div}_xu=\nabla_x\cdot u$ be the divergence in $x$, $\alpha=(\alpha_0,\alpha_1,\alpha_2,\alpha_3)$ be the multi-index, and denote
$$
\partial^\alpha=\p^{\alpha_0}_t\partial^{\alpha_1}_{x_1}\partial^{\alpha_2}_{x_2}\partial^{\alpha_3}_{x_3}.
$$
Given the mass density $\rho(t,x)$, the fluid velocity $u(t,x)$, the temperature $\theta(t,x)$, the gas constant $R$, and
$$
e:=\frac{R}{\gamma-1}\theta,\qquad \qquad p:=R\rho\theta
$$
denotes the internal energy and the pressure function, respectively. We take $R=\frac{2}{3}$ in this paper for convenience.
The local Maxwellian $\mathbf{M}$ is given by 
\begin{align}\label{M}
	\M := \mathbf{M} _{[\rho, u, \theta]}(t, x, v) := \frac{\rho(t, x)}{\sqrt{(2 \pi R \theta(t, x))^{3}}} \exp \left(-\frac{|v-u(t, x)|^{2}}{2 R \theta(t, x)}\right). 
\end{align}
We define the $L^2$ inner products (on $\mathbb{R}^3_v$ and $\mathbb{D}^{}\times\mathbb{R}^3_v$) with respect to a given Maxwellian $ \widetilde{\M} $ as 
\begin{align}
	\langle h,g\rangle_{\widetilde{\M}}& =\int_{\mathbb{R}^3}\frac{h(v)g(v)}{\widetilde{\M}}\,dv \qquad \text{with norm}\  |h|_{L_v^2(\frac{1}{\sqrt{\widetilde{\M}_*}})}^2=\int_{\mathbb{R}^3}\frac{h^2(v)}{\widetilde{\M}}\,dv. \label{in-prod1}
\end{align}
Moreover, the ${H^k_x}$ is the standard Sobolev space on $\mathbb{D}$ and the mixed space $H^k_xL_v^2(\frac{1}{\sqrt{\widetilde{\M}_*}})$ on $\mathbb{D}\times\R^3_v$ are equipped with norms 
\begin{align*}
    \|f\|_{H^k_x}^2:=\sum_{0\le j\le k}\int_{\mathbb{D}}|\nabla_x^jf|^2\,dx,\quad \text{ and }\quad \|f\|^2_{H^k_xL_v^2(\frac{1}{\sqrt{\widetilde{\M}_*}})}:=\big\|\|f\|_{L_v^2(\frac{1}{\sqrt{\widetilde{\M}_*}})}\big\|_{H^k_x}^2. 
\end{align*}
For simplicity, we write $\|f\|:=\|f\|_{L^2_x}$. 
We will write the macroscopic quantities as 
\begin{align}
	U:=(\rho,m,E)^t,\quad \bar{U}_\pm:=(\rhob_\pm,\mbb_\pm,\Eb_\pm)^t,\quad\tilde{U}_\pm:=(\rhot_\pm,\mtb_\pm,\Et_\pm)^t. 
\end{align}
where $(\cdot)^t$ denotes the transpose of the vector $(\cdot)$.

\subsection{Macro-micro decomposition}
We will use the energy method for the Boltzmann equation and summarize it as follows; see \cite{LYY,LYYZ,LY2}
for more details. 
The macro-micro decomposition of Boltzmann equation \cref{be} is related to decomposing the solution $ f $ with respect to a local Maxwellian $\M(t,x,v)$, that is 
\begin{align}
	f(t,x,v)=\mathbf{M}(t,x,v)+\mathbf{G}(t,x,v).
\end{align}
The local Maxwellian $\M$ represents the macroscopic fluid part, which is defined by five conserved quantities, the mass density $\rho(t,x)$, the momentum $m(t, x) = \rho(t, x) u(t, x)$, and the energy $ E=\rho\big(e(t, x) + \frac{1}{2}|u(t, x)|^2\big)$. More precisely, they are given by 
\begin{align}\label{fluid}
\left\{
	\begin{aligned}
		&\rho(t, x) \equiv \int_{\mathbb{R}^{3}}\Xi_{0} f(t, x, v) \,\mathrm{d} v \\
		&m_{i}(t, x) \equiv \int_{\mathbb{R}^{3}} \Xi_{i}(v) f(t, x, v) \,\mathrm{d} v \qquad \text { for } i=1,2,3 \\
		&\big(\rho e+\frac{1}{2}\rho|u|^{2}\big)(t, x) \equiv \int_{\mathbb{R}^{3}} \Xi_{4}(v) f(t, x, v) \,\mathrm{d} v,
	\end{aligned}\right.
\end{align}
where $\Xi_i$, $(i=0,1,..,4)$ are the collision invariants given by
\begin{align}\label{varphi}
	\left\{\begin{aligned}
		&\Xi_{0}(v) \equiv 1 \\
		&\Xi_{i}(v) \equiv v_{i} \qquad \text { for } i=1,2,3 \\
		&\Xi_{4}(v) \equiv \frac{1}{2}|v|^{2}
	\end{aligned}\right.
\end{align}
and satisfy
\[
\int_{\mathbb{R}^{3}} \Xi_{i}(v) Q(f, g) \,\mathrm{d} v=0 \quad \text { for } i=0,1,2,3,4.
\]
Correspondingly, $\mathbf{G}:=f-\M$ is the microscopic non-fluid part. Applying the orthogonalized method to \cref{varphi}, we have collision invariants 
\begin{align}\label{chi}
	\left\{\begin{aligned}
		&\chi_{0}(v) \equiv \frac{1}{\sqrt{\rho}} \mathbf{M}, \\
		&\chi_{i}(v) \equiv \frac{v_{i}-u_{i}}{\sqrt{R \theta \rho}} \mathbf{M}, \quad \text { for } i=1,2,3, \\
		&\chi_{4}(v) \equiv \frac{1}{\sqrt{6 \rho}} \big(\frac{|v-u|^{2}}{R \theta}-3\big) \mathbf{M},
	\end{aligned}\right.
\end{align}
that are pairwise orthogonal with respect to the inner product \cref{in-prod1}, i.e.
\begin{align}
	\left\langle\chi_{i}, \chi_{j}\right\rangle_\M =\delta_{i j}, \quad i, j=0,1,2,3,4.
\end{align}
Naturally, we can define the projection operators by the pairwise orthogonal function \cref{chi}: 
\begin{align}\label{pro}
	\Pb_{0} h \equiv \sum_{i=0}^{4}\left\langle h, \chi_{i}\right\rangle_\M \chi_{i}, \quad \Pb_{1} h \equiv h-\Pb_{0} h,
\end{align}
which are
called the macroscopic projection and the microscopic projection, respectively.
Direct calculation yields 
\begin{align*}
	\Pb_{0} \Pb_{0}=\Pb_{0}, \quad \Pb_{1} \Pb_{1}=\Pb_{1}, \quad \Pb_{1} \Pb_{0}=\Pb_{0} \Pb_{1}=0.
\end{align*}
Then by \cref{fluid} and \cref{M}, we have
\[
\Pb_{0}f=\mathbf{M}, \ \ \ \ \ \ \ \ \ \Pb_{1}f=\mathbf{G}.
\]
Thus, the Boltzmann equations \cref{be} can be rewritten as 
\begin{align}\label{MG}
	(\mathbf{M}+\mathbf{G})_{t}+v\cdot\nabla_x\left(\mathbf{M}+\mathbf{G}\right)=Q(\mathbf{G}, \mathbf{M})+Q(\mathbf{M}, \mathbf{G})+Q(\mathbf{G}, \mathbf{G}).
\end{align}
Multiplying \cref{MG} by the collision invariant $\Xi_i$ ($i=0,1,..,4$) given by \cref{varphi}, we have the fluid-type system for the macroscopic parts: 
\begin{align}\label{nsg}
	\left\{\begin{array}{l}
		\rho_{t}+\operatorname{div}_x(\rho u)=0 \\[2mm]
		(\rho u)_{t}+\operatorname{div}(\rho u \otimes u)+\nabla_x p=-\int v \otimes v \cdot \nabla_x \mathbf{G}\, d v \\[2mm]
		{\big[\rho\big(e+\frac{|u|^{2}}{2}\big)\big]_{t}+\operatorname{div}_x\big[\rho u\big(e+\frac{|u|^{2}}{2}\big)+p u\big]=-\int \frac{1}{2}|v|^{2} v \cdot \nabla_x \mathbf{G}\, d v},
	\end{array}\right.
\end{align}
where $p := \frac{2}{3}
\rho e = R\rho \theta$ is the pressure for the mono-atomic gas. Applying $\Pb_1$ to \cref{MG}, the microscopic part $\mathbf{G}$ satisfies 
\begin{align}\label{g}
	\mathbf{G}_{t}+\Pb_{1}\left(v\cdot \nabla_x \mathbf{G}\right)+\Pb_{1}\left(v\cdot \nabla_x \mathbf{M}\right)=\LL_{\mathbf{M}} \mathbf{G}+Q(\mathbf{G}, \mathbf{G}),
\end{align}
where $\LL_\mathbf{M}$ is the linearized operator around the local Maxwellian $\M$, that is
\begin{align}\label{Lm}
	\LL_{\mathbf{M}} g=Q(\mathbf{M}+g, \mathbf{M}+g)-Q(\mathbf{G}, \mathbf{G})=2 Q(\mathbf{M}, g).
\end{align}
Then the microscopic part $\mathbf{G}$ takes the following form
\begin{align}\label{GL}
	\mathbf{G}=\LL_{\mathbf{M}}^{-1}\left(\Pb_{1}\left(v\cdot \nabla_x \mathbf{M}\right)\right)+\Pi,
\end{align}
where
\begin{align}\label{the}
	\Pi=\LL_{\mathbf{M}}^{-1}\Big[\mathbf{G}_{t}+\Pb_{1}\left(v\cdot \nabla_x \mathbf{G}\right)-Q(\mathbf{G}, \mathbf{G})\Big].
\end{align}
Direct calculation implies that the null space of $\LL_\mathbf{M}$, denoted by $\mathfrak{N}$, is spanned by $\{\chi_\alpha\}_{0\le\alpha\le 4}$ given by \cref{chi}. Moreover, we can rewrite the expression of the linearized collision operator $\LL_\mathbf{M}$ as
\begin{align}\label{lp}
	\left(\LL_{\mathbf{M}} h\right)(v)=-\nu_{\mathbf{M}}(v) h(v)+K_{\mathbf{M}} h(v),
\end{align}
where $\nu_{\mathbf{M}}(v)$ is the collisional frequency given by
\begin{align}
	\nu_{\mathbf{M}}(v)=\int_{\mathbb{R}^{3} \times \mathbb{S}^{2}}\left|v-v_-\right|^{\gamma} \mathbf{M}\left(v_-\right) B(\vartheta) \,\mathrm{d} v_- \mathrm{d} \omega=c \int_{\mathbb{R}^{3}}\left|v-v_-\right|^{\gamma} \mathbf{M}\left(v_-\right) \,\mathrm{d} v_-,
\end{align}
for some constant $c>0$, and $K_{\M}:=K_{2 \M}-K_{1 \M}$ is given by 
\[
\begin{array}{l}
	K_{1 \mathbf{M}} g=\int_{\mathbb{R}^{3} \times \mathbb{S}^{2}}\left|v-v_-\right|^{\gamma} \mathbf{M}(v) g\left(v_-\right) B(\vartheta) \mathrm{d} v_- \mathrm{d} \omega, \\[2mm]
	K_{2 \mathbf{M}} g=\int_{\mathbb{R}^{3} \times \mathbb{S}^{2}}\left|v-v_-\right|^{\gamma}\left\{\mathbf{M}\left(v^{\prime}\right) g\left(v_-^{\prime}\right)+g\left(v^{\prime}\right) \mathbf{M}\left(v_-^{\prime}\right)\right\} B(\vartheta) \mathrm{d} v_- \mathrm{d} \omega.
\end{array}
\]
Utilizing the exponential decay in $\M$, the collision frequency $\nu_\M$ satisfies 
\begin{align}\label{cf}
	\nu_{0}\left(1+|v-u|^{2}\right)^{\gamma / 2} \leq \nu_{\M}(v) \leq \nu_{1}\left(1+|v-u|^{2}\right)^{\gamma / 2},
\end{align}
for some positive constants $\nu_0$ and $\nu_1$ which depend on the lower and upper bounds of $(\rho,u,\theta)$ (which are universal constants by assuming the small perturbation below). 
Moreover, $\LL_\mathbf{M}$ is dissipative, i.e., there exists a positive constant $\sigma_0$, such that for any $h \in \mathfrak{N}^{\perp}$, we have 
\begin{align}\label{dp}
	\langle h,\LL_\mathbf{M}h\rangle_{\mathbf{M}}\leq-\sigma_0\langle \nu_\mathbf{M}(v)h,h\rangle_{\mathbf{M}}.
\end{align}
Substituting \cref{GL} into \cref{nsg}, we obtain the fluid-type system 
\begin{equation}\label{NSG}
	\begin{cases}
		\rho_{t}+\operatorname{div}_{x}(\rho u)=0, \\[1mm]
		(\rho u)_{t}+\operatorname{div}_{x}(\rho u \otimes u)+\nabla_{x} p \\
		\qquad\qquad\qquad =-\int v \otimes v \cdot \nabla_{x}\left(\LL_{\mathbf{M}}^{-1}\left[\mathbf{P}_{1}\left(v \cdot \nabla_{x} \mathbf{M}\right)\right]\right) \,d v-\int v \otimes v \cdot \nabla_{x} \Pi\, d v, \\[1mm]
		{\big[\rho\big(e+\frac{|u|^{2}}{2}\big)\big]_{t}+\operatorname{div}_{x}\big[\rho u\big(e+\frac{|u|^{2}}{2}\big)+p u\big]} \\
		\qquad\qquad\qquad =-\int \frac{1}{2}|v|^{2} v \cdot \nabla_{x}\left(\LL_{\mathbf{M}}^{-1}\left[\mathbf{P}_{1}\left(v \cdot \nabla_{x} \mathbf{M}\right)\right]\right)\,\, d v-\int \frac{1}{2}|v|^{2} v \cdot \nabla_{x} \Pi\, d v.
	\end{cases}
\end{equation}
Further, a direct calculation yields (cf. \cite{XYY})
\begin{align}\label{diff}
	&-\int v_{i} v \cdot \nabla_{x}\left(\LL_{\mathbf{M}}^{-1}\left[\mathbf{P}_{1}\left(v \cdot \nabla_{x} \mathbf{M}\right)\right]\right) d v=\sum_{j=1}^{3}\left[\mu(\theta)\left(u_{i x_{j}}+u_{j x_{i}}-\frac{2}{3} \delta_{i j} \operatorname{div}_{x} u\right)\right]_{x_{j}}=:\sum_{j=1}^3[S_{ij}]_{x_j}, \\
	&-\int \frac{1}{2}|v|^{2} v \cdot \nabla_{x}\left(\LL_{\mathbf{M}}^{-1}\left[\mathbf{P}_{1}\left(v \cdot \nabla_{x} \mathbf{M}\right)\right]\right) d v 
	=\dv\left(\kappa(\theta) \nabla\theta\right)+\dv\left\{\mu(\theta) \uf\cdot \mathbf{S}\right\},\nonumber
\end{align} 
where $\mathbf{S}:=(S_{ij})\in\R^{3\times3}$, $i,j=1,2,3$ is given by \cref{diff}. 
We use the following notations to describe the macroscopic parts of the solution $ f $ to \cref{be}. Write the nonlinear part as 
\begin{align}\label{calG}
	\begin{aligned}
		\mathcal{G}_2(U):=&\frac{m\otimes m}{\rho}+\frac{2}{3}(E-\frac{|m|^2}{2\rho})\mathbb{I}-\mathbf{S},\\\quad \mathcal{G}_3(U):=&\frac{m E}{\rho}+pu-\kappa(\T)\nabla_x(\frac{E}{\rho}-\frac{|m|^2}{2\rho^2})-\frac{m}{\rho}\cdot\mathbf{S},
	\end{aligned}
\end{align}
where $\mathbb{I}$ denotes the identity matrix.
Using these notations, \cref{NSG} can be rewritten as
\begin{align}\label{NS-re}
\left\{
	\begin{aligned}
		&\p_t\rho+\dv_x m=0,\\
		&\pt m+\dv_x\mathcal{G}_2(U)=-\int v \otimes v \cdot \nabla_{x} \Pi\, d v,\\
		&\pt E+\dv_x \mathcal{G}_3(U)=-\int \frac{1}{2}|v|^{2} v \cdot \nabla_{x} \Pi\, d v.
	\end{aligned}\right. 
\end{align}


\subsection{1-D Shock profile of the Boltzmann equation}
In this Subsection, we introduce the shock wave profiles, the traveling wave solutions, of the Boltzmann equation \cref{be}
with the following kind of initial data
\begin{align}\label{ini}
	f(0, x, v)=f_0(x, v) \rightarrow \mathbf{M}_{\left[\rho_{\pm}, u_{\pm}, \theta_{\pm}\right]}(v) \quad \text{as}\quad x_1 \rightarrow \pm \infty,
\end{align}
where $\rho_\pm>0$, $u_{\pm}=\left(u_{1 \pm}, 0,0\right)^t$ and $\theta_\pm>0$ are given states. With the macro-micro decomposition, one can deduce that the macroscopic initial data satisfy 
\begin{align}
	(\rho, u, \theta)(x, 0)=\left(\rho_0, u_0, \theta_0\right)(x), \quad \lim _{x \rightarrow \pm \infty}\left(\rho_0, u_0, \theta_0\right)(x)=\left(\rho_{\pm}, u_{\pm}, \theta_{\pm}\right).
\end{align}
In this paper, we are interested in the stability of shock waves for the Boltzmann equation, which is closely related to  the Riemann problem for the compressible Euler equation 
\begin{equation}\label{EE1}
	\begin{cases}
		\rho_{t}+\operatorname{div}_{x}(\rho u)=0, \\
		(\rho u)_{t}+\operatorname{div}_{x}(\rho u \otimes u)+\nabla_{x} p =0, \\
		{\big[\rho\big(e+\frac{|u|^{2}}{2}\big)\big]_{t}+\operatorname{div}_{x}\big[\rho u\big(e+\frac{|u|^{2}}{2}\big)+p u\big]} =0.
	\end{cases}
\end{equation}
with the Riemann initial data
\begin{align}\label{initial-Euler}
(\rho, m, E)(x, 0)=\left\{\begin{array}{l}
	\left(\rho_{-}, m_{-}, E_{-}\right), x<0, \\
	\left(\rho_{+}, m_{+}, E_{+}\right), x>0,
\end{array}\right.
\end{align}
where $m_{\pm}=\rho_{\pm} u_{\pm}, E_{\pm}=\rho_{\pm}\left(\theta_{\pm}+\frac{\left|u_{\pm}\right|^2}{2}\right)$. There are three eigenvalues for the system \cref{EE1}:
\begin{align}
	\lambda_1=u_1-\frac{\sqrt{10 \theta}}{3},\quad \lambda_2=u_1,\quad \lambda_3=u_1+\frac{\sqrt{10 \theta}}{3},
\end{align} 
where the second characteristic field is linear degenerate and the others are genuinely nonlinear. 

\smallskip 
For the case of the composite wave consisting of two shock waves, the Riemann solution can be described as follows (see for example \cite{smoller}). There exists an intermediate constant state between $+$ and $-$, 
\begin{align}\label{Ex}
	\Big(\rho_{\#}, \quad m_{\#}=\rho_{\#} u_{\#},\quad E_{\#}=\rho_{\#}\big(\theta_{\#}+\frac{\left|u_{\#}\right|^2}{2}\big)\Big),
\end{align} 
such that Rankine-Hugoniot conditions
\begin{align}\label{Rh1}
\left\{\begin{array}{l}
	-s_3\left(\rho_{+}-\rho_{\#}\right)+(\rho_{+} u_{1+}-\rho_{\#} u_{1 \#})=0, \\ -s_3(\rho_{+} u_{1+}-\rho_{\#} u_{1 \#})+\big(\rho_{+} u_{1+}^2+p_{+}-\rho_{\#} u_{1 \#}^2-p_{\#}\big)=0, \\ -s_3(\rho_{+} E_{+}-\rho_{\#} E_{\#})+\big(\rho_{+} u_{1+} E_{+}+p_{+} u_{1+}-\rho_{\#} u_{1 \#} E_{\#}-p_{\#} u_{1 \#}\big)=0,
\end{array}\right.\\
\left\{\begin{array}{l}\label{Rh2}
	-s_1(\rho_{\#}-\rho_{-})+\left(\rho_{\#} u_{1 \#}-\rho_{-} u_{1-}\right)=0, \\ -s_1(\rho_{\#} u_{1 \#}-\rho_{-} u_{1-})+\big(\rho_{\#} u_{1 \#}^2+p_{\#}-\rho_{-} u_{1-}^2-p_{-}\big)=0, \\ -s_1(\rho_{\#} E_{\#}-\rho_{-} E_{-})+\big(\rho_{\#} u_{1 \#} E_{\#}+p_{\#} u_{1 \#}-\rho_{-} u_{1-} E_{-}-p_{-} u_{1-}\big)=0
\end{array}\right.
\end{align}
and Lax entropy conditions
\begin{align}\label{lax}
	\lambda_{3+}<s_3<\lambda_{3 \#}, \quad \lambda_{1 \#}<s_1<\lambda_{1-}
\end{align}
holds. Here, $s_i$ $(i=1,3)$ are the $i$-shock wave speeds and 
\begin{align}
\lambda_{1-}=u_{1-}-\frac{\sqrt{10 \theta_{-}}}{3},\quad \lambda_{i \#}=u_{1 \#}+(-1)^{\frac{i+1}{2}} \frac{\sqrt{10 \theta_{\#}}}{3}\quad \text{and}\quad \lambda_{3+}=u_{1+}+\frac{\sqrt{10 \theta_{+}}}{3}.
\end{align}
To describe the strengths of the shock waves for later use, we set
\begin{align}\label{wave strength}
	\delta^{s_1}=\left|\rho_{\#}-\rho_{-}\right|+\left|m_{\#}-m_{-}\right|+\left|E_{\#}-E_{-}\right|, \quad \delta^{s_3}=\left|\rho_{\#}-\rho_{+}\right|+\left|m_{\#}-m_{+}\right|+\left|E_{\#}-E_{+}\right|,
\end{align}  and 
\begin{align}\label{deltaDef}
    \delta=\min \left\{\delta^{s_1}, \delta^{s_3}\right\}.
\end{align}
Now we fix $\left(\rho_{-}, m_{-}, E_{-}\right)$ and choose $\left|\left(\rho_{+}-\rho_{-}, m_{+}-m_{-}, E_{+}-E_{-}\right)\right|$ sufficiently small, then 
$$
\delta^{s_1}+\delta^{s_3} \leq C\left|\left(\rho_{+}-\rho_{-}, m_{+}-m_{-}, E_{+}-E_{-}\right)\right|,
$$
where $C$ is a positive constant depending only on $\left(\rho_{-}, m_{-}, E_{-}\right)$. Then, following \cite{HM}, if it also holds that
\begin{align}\label{wave strength-same order}
\delta^{s_1}+\delta^{s_3} \leq C \delta \quad \text { as } \quad \delta^{s_1}+\delta^{s_3} \rightarrow 0,
\end{align}
for a positive constant $C$, we say the strengths of the shock waves are "small with the same order." In the following, we will assume that \cref{wave strength-same order} holds.

\smallskip Next, we recall some basic theory of the one-dimensional $i$-shock profile $F^{s_i}\left(x_1-s_i t, v\right)$ $(i=1,3)$ of the Boltzmann equation \cref{be} in Eulerian coordinates; the global existence and properties can be found in \cite{CN,LYs,LY2}. The $i$-shock profile $F^{s_i}\left(x_1-s_i t, v\right)$ are travelling wave solutions satisfying
\begin{align}\label{Fsi}
\begin{cases}
	-s_i\left(F^{s_i}\right)^{\prime}+v_1\left(F^{s_i}\right)^{\prime}=Q\left(F^{s_i}, F^{s_i}\right), \quad i=1,3, \\ F^{s_1}(-\infty, v)=\mathbf{M}_{\left[\rho_{-}, u_{-}, \theta_{-}\right]}(v), \qquad F^{s_3}(-\infty, v)=\mathbf{M}_{\left[\rho_{\#}, u_{\#}, \theta_{\#}\right]}(v), \\ F^{s_1}(+\infty, v)=\mathbf{M}_{\left[\rho_{\#}, u_{\#}, \theta_{\#}\right]}(v), \qquad F^{s_3}(+\infty, v)=\mathbf{M}_{\left[\rho_{+}, u_{+}, \theta_{+}\right]}(v),
\end{cases}
\end{align}
where $(\cdot)^{\prime}:=\frac{d}{d \vartheta_i}(\cdot),$ and $\vartheta_i:=x_1-s_i t$.
By the micro-macro decomposition around the local Maxwellian $\mathbf{M}^{s_i}(i=1,3)$, set
$$
F^{s_i}\left(x_1-s_i t, v\right)=\mathbf{M}^{s_i}\left(x_1-s_i t, v\right)+\mathbf{G}^{s_i}\left(x_1-s_i t, v\right),
$$
where 
\begin{align}\label{local maxwellian}
	\mathbf{M}^{s_i}\left(x_1-s_i t, v\right)=\mathbf{M}_{\left[\rho^{s_i}, u^{s_i}, \theta^{s_i}\right]}\left(x_1-s_i t, v\right)=\frac{\rho^{s_i}\left(x_1-s_i t\right)}{\sqrt{\left(2 \pi R \theta^{s_i}\left(x_1-s_i t\right)\right)^3}} e^{-\frac{| v-u^{s_i}(x_1-s_i t)|^2}{2 R \theta^{s_i}\left(x_1-s_i t\right)}}
\end{align}
with
\begin{align*}
    \begin{pmatrix}
	\rho^{s_i} \\
	\rho^{s_i} u_j^{s_i} \\[1mm]
\rho^{s_i}\left(\theta^{s_i}+\frac{\left|u^{s_i}\right|^2}{2}\right)
\end{pmatrix}
:=\int_{\mathbb{R}^3}
\begin{pmatrix}
	1 \\
	v_j \\[1mm]
	\frac{|v|^2}{2}
\end{pmatrix} F^{s_i}\left(x_1-s_i t, v\right) d v, \quad j=1,2,3 .
\end{align*}
The corresponding macroscopic projection $\mathbf{P}_0^{s_i}$ and microscopic projection $\mathbf{P}_1^{s_i}$ are given by
$$
\mathbf{P}_0^{s_i} g=\sum_{j=0}^4\left\langle g, \chi_j^{s_i}\right\rangle_{\mathbf{M}^{s_i}} \chi_j^{s_i}, \quad \mathbf{P}_1^{s_i} g=g-\mathbf{P}_0^{s_i} g,
$$
where $\chi_j^{s_i}$ $(0\le j\le 4)$ and $\langle\cdot, \cdot\rangle_{\mathbf{M}^{s_i}}$ are defined in \cref{dp,chi} respectively with respect to $\mathbf{M}^{s_i}$. Then rewritting $F^{s_i}=\mathbf{M}^{s_i}+\mathbf{G}^{s_i}$, the Boltzmann equation \cref{Fsi} can be rewritten as 
$$
\left(\mathbf{M}^{s_i}+\mathbf{G}^{s_i}\right)_t+v_1\left(\mathbf{M}^{s_i}+\mathbf{G}^{s_i}\right)_{x_1}=2 Q\left(\mathbf{M}^{s_i}, \mathbf{G}^{s_i}\right)+Q\left(\mathbf{G}^{s_i}, \mathbf{G}^{s_i}\right) .
$$
Correspondingly, the system for fluid components of the shock profile is 
\begin{align}\label{sw-fluid}
\left\{\begin{aligned}
	&\pt\rho^{s_i}+\px\left(\rho^{s_i} u_1^{s_i}\right)=0, \\
	&\pt\left(\rho^{s_i} u_1^{s_i}\right)+\px\left[\rho^{s_i}\left(u_1^{s_i}\right)^2+p^{s_i}\right]=\frac{4}{3}\px\left(\mu\left(\theta^{s_i}\right) \px u_{1 }^{s_i}\right)-\int v_1^2 \px\Pi^{s_i} d v, \\
	&\pt\left(\rho^{s_i} u_j^{s_i}\right)+\px\left(\rho^{s_i} u_1^{s_i} u_j^{s_i}\right)=\px\left(\mu\left(\theta^{s_i}\right) \px u_{j}^{s_i}\right)-\int v_1 v_j \px\Pi^{s_i} d v,\ \  j=2,3, \\
	&\pt\left[\rho^{s_i}\left(\theta^{s_i}+\frac{\left|u^{s_i}\right|^2}{2}\right)\right]+\px\left[\rho^{s_i} u_1^{s_i}\left(\theta^{s_i}+\frac{\left|u^{s_i}\right|^2}{2}\right)+p^{s_i} u_1^{s_i}\right]=\px\left(\kappa\left(\theta^{s_i}\right) \px\theta^{s_i}\right) \\
	&\quad+\frac{4}{3}\px\big(\mu\left(\theta^{s_i}\right) u_1^{s_i} \px u_{1}^{s_i}\big)+\sum_{j=2}^3\px\left(\mu\left(\theta^{s_i}\right) u_j^{s_i} \px u_{j}^{s_i}\right)-\int \frac{1}{2} v_1|v|^2 \px\Pi^{s_i} d v, 
\end{aligned}\right.
\end{align}
while the equation for the nonfluid component $\mathbf{G}^{s_i}$ $(i=1,3)$ is
\begin{align}\label{Gsi}
\pt\mathbf{G}^{s_i}+\mathbf{P}_1^{s_i}\left(v\cdot\nabla_x  \mathbf{M}^{s_i}\right)+\mathbf{P}_1^{s_i}\left(v\cdot\nabla_x\mathbf{G}^{s_i}\right)=\LL_{\mathbf{M}^{s_i}} \mathbf{G}^{s_i}+Q\left(\mathbf{G}^{s_i}, \mathbf{G}^{s_i}\right) .
\end{align}
where $\LL_{\mathbf{M}^{s_i}}$ is the linearized collision operator of $Q\left(F^{s_i}, F^{s_i}\right)$ with respect to the local Maxwellian $\mathbf{M}^{s_i}$:
$$
\LL_{\mathbf{M}^{s_i}} g:=2 Q\left(\mathbf{M}^{s_i}, g\right)=Q\left(\mathbf{M}^{s_i}, g\right)+Q\left(g, \mathbf{M}^{s_i}\right) .
$$
Further,
\begin{align}
\label{the1}
\begin{aligned}
	& \mathbf{G}^{s_i}:=\LL_{\mathbf{M}^{s_i}}^{-1}\left[\mathbf{P}_1^{s_i}\left(v\cdot \nabla_x\mathbf{M}^{s_i}\right)\right]+\Pi^{s_i}, \\
	& \Pi^{s_i}:=\LL_{\mathbf{M}^{s_i}}^{-1}\left[\pt\mathbf{G}^{s_i}+\mathbf{P}_1^{s_i}\left(v\cdot\nabla_x \mathbf{G}^{s_i}\right)-Q\left(\mathbf{G}^{s_i}, \mathbf{G}^{s_i}\right)\right].
\end{aligned}
\end{align}
Now we can state the important Lemma given by \cite{LY2}:
\begin{Lem}\label{sw}
Let $(\rho_{\pm}, u_{\pm}, \theta_{\pm})$ and $(\rho_{\#},u_{\#},\theta_{\#})$ be the given states satisfying \cref{Rh1,Rh2,lax}, and $\delta^{s_i}$ be the shock wave strength defined in \cref{wave strength}.  If $\delta^{s_i}$ is small enough, then the Cauchy problem of Boltzmann equation \cref{Fsi} admits a unique traveling wave solution $F^{s_i}\left(x_1-s_i t, v\right)$  up to a shift, which satisfies the following properties:
\begin{itemize}

\item It holds that
\begin{align}\label{decayshock}
\left\{\begin{aligned}
&	\left|\left(\rho^{s_1}-\rho_{-}, u_1^{s_1}-u_{1-}, \theta^{s_1}-\theta_{-}\right)\right| \leq C \delta^{s_1} e^{-c \delta^{s_1}\left|\vartheta_1\right|} \quad \text { as } \quad \vartheta_1<0, \\[2mm]
&	\left|\left(\rho^{s_1}-\rho_{\#}, u_1^{s_1}-u_{1 \#}, \theta^{s_1}-\theta_{\#}\right)\right| \leq C \delta^{s_1} e^{-c \delta^{s_1}\left|\vartheta_1\right|} \quad \text { as } \quad \vartheta_1>0, \\[2mm]
&	\left|\left(\rho^{s_3}-\rho_{+}, u_1^{s_3}-u_{1+}, \theta^{s_3}-\theta_{+}\right)\right| \leq C \delta^{s_3} e^{-c \delta^{s_3}\left|\vartheta_3\right|} \quad \text { as } \quad \vartheta_3>0, \\[2mm]
&	\left|\left(\rho^{s_3}-\rho_{\#}, u_1^{s_3}-u_{1 \#}, \theta^{s_3}-\theta_{\#}\right)\right| \leq C \delta^{s_3} e^{-c \delta^{s_3}\left|\vartheta_3\right|} \quad \text { as } \quad \vartheta_3<0, \\[2mm]
&	\left(\int \frac{\nu(|v|)\left|\mathbf{G}^{s_i}\right|^2}{\mathbf{M}_0} d v\right)^{\frac{1}{2}} \leq C\left(\delta^{s_i}\right)^2 e^{-c \delta^{s_i}\left|\vartheta_i\right|}, \quad i=1,3,
\end{aligned}\right.
\end{align}
where $\mathbf{M}_0$ is the global Maxwellian close to the shock profile and $\vartheta_i:=x_1-s_i t$. For the precise definition, we refer to \cite{LY2}.

\smallskip 
\item Compressibility of $i$-shock profile:
\begin{align}\label{compressibility}
\left(\lambda_i^{s_i}\right)_{\vartheta_i}<0, \quad \lambda_i^{s_i}=u_1^{s_i}+(-1)^{\frac{i+1}{2}} \frac{\sqrt{10 \theta^{s_i}}}{3} .
\end{align}

\item It holds that
\begin{align}\label{derismall}
\rho_{\vartheta_i}^{s_i} \sim u_{1 \vartheta_i}^{s_i} \sim \theta_{\vartheta_i}^{s_i} \sim\left(\lambda_i^{s_i}\right)_{\vartheta_i} \sim\left(\int \frac{\nu(|v|)\left|\mathbf{G}^{s_i}\right|^2}{\mathbf{M}_0} d v\right)^{\frac{1}{2}}\leq C\delta^2,
\end{align}
where $A \sim B$ denotes the equivalence of the quantities $A$ and $B$, and
\begin{align}\label{derismall2}
\left\{\begin{aligned}
	&u_j^{s_i} \equiv 0, \quad \int v_j \Pi^{s_i} d v \equiv 0, \quad j=2,3, \\[2mm]
&\left|\partial_{\vartheta_i}^k\left(\rho^{s_i}, u_1^{s_i}, \theta^{s_i}\right)\right| \leq C\left(\delta^{s_i}\right)^{k-1}\left|\left(\rho_{\vartheta_i}^{s_i}, u_{1 \vartheta_i}^{s_i}, \theta_{\vartheta_i}^{s_i}\right)\right|, \quad k \geq 2, \\[2mm]
	&\left(\int \frac{\nu(|v|)\left|\partial_{\vartheta_i}^k \mathbf{G}^{s_i}\right|^2}{\mathbf{M}_0} d v\right)^{\frac{1}{2}} \leq C\left(\delta^{s_i}\right)^k\left(\int \frac{\nu(|v|)\left|\mathbf{G}^{s_i}\right|^2}{\mathbf{M}_0} d v\right)^{\frac{1}{2}}, \quad k \geq 1, \\[2mm]
	&\left|\int v_1 \Xi_j(v) \Pi_{\vartheta_i}^{s_i} d v\right| \leq C \delta^{s_i}\left|u_{1 \vartheta_i}^{s_i}\right|, \quad j=1,2,3,4,
\end{aligned}\right.
\end{align}
where $\Xi_j(v)$ are the collision invariants defined in \cref{varphi}.
\end{itemize}
\end{Lem}
In this work, we consider the linear superposition of 1-shock and 3-shock profiles
$\bar{U}(x, t)=(\bar{\rho}, \bar{m}, \bar{E})^t$ given by 
\begin{align}\label{shock wave}
\left\{\begin{array}{l}
\bar{\rho}=\rho^{s_1}\left(x_1-s_1t\right)+\rho^{s_3}\left(x_1-s_3t\right)-\rho_{\#},\\	\bar{m}_1=m_1^{s_1}\left(x_1-s_1t\right)+m_1^{s_3}\left(x_1-s_3t\right)-m_{1 \#}, \quad \bar{m}_i=0, \quad i=2,3,\\
 \bar{E}=E^{s_1}\left(x_1-s_1t\right)+E^{s_3}\left(x_1-s_3t\right)-E_{\#}, 
\end{array}\right.
\end{align}
and consider the stability of this profile $\bar{U}(x,v)$. 

\subsection{Construction of the ansatz with decoupled diffusion wave}
Now we construct the ansatz for our perturbation theory. 
In this work, we consider the general initial perturbation such that the initial perturbed macroscopic quantities may not be zero: 
\begin{align}\label{inimass}
    \int_{\mathbb{D}}(U-\bar{U})(x, 0) \,d x=c_0\neq0; 
\end{align}
but we will assume that $c_0$ is small enough. 
 To apply the anti-derivative technique, we need to find an ansatz $\tilde{U}$ such that 
 \begin{align}\label{zm}
 	\int_{\mathbb{D}} \big( U(x, t)-\tilde{U}(x, t)\big)\,d x=0 ,\quad \text{ and }\quad |\tilde{U}-\bar{U}|\rightarrow 0 \quad \text{as} \quad t \rightarrow \infty.
 \end{align}
Inspired by \cite{Liu1,HM}, we will construct the decoupled diffusion wave and determine two shifts on the shock profile to obtain the desired ansatz $\tilde{U}$. However, the interaction between transveral direction $(m_2,m_3)$ and the principal direction $m_1$ is very difficult, which makes the energy of anti-derivative increase. The previous results often assumed the initial mass of the transversal direction to be zero. We then decompose the profile into two parts.
\begin{align}
	U:=(U_1,U_2),\qquad \text{where}\qquad U_1:=(\rho,m_1,E),\quad U_2:=(m_2,m_3).
\end{align} 
The extra mass of $U_1$ will be removed and we won't take anti-derivatives on $U_2$. Instead of which, we will obtain the energy estimate of $U_2$ with the help of the entropy-entropy pair; see \cref{EEP}.
We start studying
$$
A\left(\rho, m_1, E\right)=\left(\begin{array}{ccc}
	0 & 1 & 0 \\
	-\frac{m_1^2}{\rho^2}+\frac{m_1^2}{3 \rho^2} & \frac{4 m_1}{3 \rho} & \frac{2}{3} \\[1mm]
	-\frac{5 m_1 E}{3 \rho^2}+\frac{2 m_1^2 m_1}{3 \rho^3} & \frac{5 E}{3 \rho}-\frac{2 m_1^2}{3 \rho^2}-\frac{m_1^2}{3 \rho^2} & \frac{5 m_1}{3 \rho}
\end{array}\right),
$$
which is the Jacobi matrix of the flux $\big(m_1, \frac{2}{3} E+\frac{m_1^2}{\rho}-\frac{m^2}{3 \rho}, \frac{5 m_1 E}{3 \rho}-\frac{m_1 m^2}{3 \rho^2}\big)^t$ of the Euler system \cref{EE1} with respect to $\left(\rho, m_1, E\right)$. Since the wave strengths $\delta^{s_1}$ and $\delta^{s_3}$ (given in \cref{wave strength}) are small, the following three vectors are linearly independent (cf. \cite{Liu1,SX}):
\begin{align*}
   & r_1=\left(\rho_{\#}-\rho_{-}, m_{1\#}-\right.\left.m_{1-}, E_{\#}-E_{-}\right)^t,\qquad r_2=\big(1, u_{1 \#}, \frac{u_{1 \#}^2}{2}\big)^t , \\
   &r_3=\left(\rho_{+}-\rho_{\#}, m_{1+}-m_{1\#}, E_{+}-E_{\#}\right)^t,
\end{align*}
where $r_2$ is the second right eigenvector of the $\operatorname{matrix} A\left(\rho, m_1, E\right)$ at the point $\left(\rho_{\#}, m_{1 \#}, E_{\#}\right)$. 
Next, one can distribute the initial mass along  $r_1, r_2$, and $r_3$, that is, there exist three constants that depend only on the initial perturbation satisfying
\begin{align}\label{initialal}
\int_{\mathbb{D}}(U_1(x, 0)-\bar{U}_1(x, 0))\, d x=\sum_{i=1}^3 \alpha_i r_i.
\end{align}
As in \cite{Liu1,SX,HM}, 
{\cite[pp. 1254, Eq. (2.a25)]{HXY}}, 
the excess masses $\alpha_1r_1$ and $\alpha_3r_3$ can be removed by translating the 1-shock profile and 3-shock profile with the shifts $\alpha_1$ and $\alpha_3$, i.e. $F^{s_1}(x-s_1t+\alpha_1)$ and  $F^{s_3}(x-s_3t+\alpha_3)$, respectively. Whlie for $\alpha_2r_2$, we introduce a linear diffusion wave $\Theta$:
\begin{align}\label{diffusionwave}
	\p_t\Theta+u_{1 \#}\p_1 \Theta=a \px^2\Theta,\left.\quad \Theta\right|_{t=-1}=\alpha_2 \delta(x_1), \quad \int_{-\infty}^{\infty} \Theta(x_1, t) d x_1=\alpha_2,
\end{align}
where $a>0$ is a constant that will be determined later. 
Then one has
\begin{align}
	\Theta(x, t)=\frac{\alpha_2}{\sqrt{4 \pi a(1+t)}} e^{-\frac{\left(x_1-u_{1 \#} t\right)^2}{4 a(1+t)}},
\end{align}

Then the ansatz satisfying \cref{zm} can be written as
\begin{align}
	(\hat{\rho},\hat{m}_1,\hat{E}):=(\bar{\rho},\bar{m}_1,\bar{E})+\Theta r_2	.
\end{align}
However, due to the diffusion wave and the interactions of two shock profiles, $(\hat{\rho},\hat{m}_1,\hat{E})$ is not the exact solution of a fluid-type equation like \cref{NSG} or \cref{sw-fluid}, thus the estimates of such errors are necessary and important. A direct calculation yields
\begin{align}
	\p_t \hat{\rho}+\px\hat{m}_{1}=a\px^2\Theta.
\end{align}
The error term is  $a\px^2\Theta$ with a poor decay rate that makes the energy of anti-derivatives increase. Thus, we shall find higher-order corrections. One immediate idea is to set
\begin{align}\label{tirhomE1}	\big(\tilde{\rho},\tilde{m}_1,\tilde{E}\big):=\big(\hat{\rho},\hat{m}_1,\hat{E}\big)+(g_1,g_2,g_3),
\end{align}
where $g_i$ ($i=1,2,3$) don't carry any mass, i.e.
$$
\int_{\mathbb{D}}g_idx=0.
$$ Then one has
\begin{align}
		\p_t \tilde{\rho}+\px\tilde{m}_{1}=\p_t\Theta+u_{1\#}\px\Theta+\p_tg_1+\px g_2.
\end{align}
Letting $g_1=0$, $g_2=-a\px\Theta$,  one has 
\begin{align}
	\p_t \tilde{\rho}+\px\tilde{m}_{1}=0.
\end{align}
To deal with the error terms in the nonlinear equations involving $\pt\mt_1$ and $\pt\Et$, we should first study the interaction between two shock profiles and the diffusion wave. In fact, by \cref{sw}, one has
\begin{align}\label{1}\notag
|\rho^{s_1}-\rho_{\#}|| \rho^{s_3}-\rho_{\#}| & \leq C \delta^{s_1} \delta^{s_3}\left(e^{-c \delta^{s_1}(|x_1|+t)+c \delta^{s_1}\left|\alpha_1\right|}+e^{-c \delta^{s_3}(|x_1|+t)+c \delta^{s_3}\left|\alpha_3\right|}\right) \\
&\leq C \delta^2 e^{-c \delta(|x|+t)},
\end{align}
and for $i=1,3$,
\begin{align}\label{2}
\left|\rho^{s_i}-\rho_{\#}\right| |\Theta| \leq C\left|\alpha_2\right| \delta^{\frac{3}{2}} e^{-c \delta(|x_1|+t)}+C \frac{\left|\alpha_2\right|}{(1+t)^{\frac{3}{2}}} e^{-\frac{c\left(x_1-u_{1 \#} t\right)^2}{1+t}}+C\left(\delta+\left|\alpha_2\right|\right) e^{-c(|x_1|+t)},
\end{align}
where we split the cases as in \cref{decayshock} and used the fact that $\delta^{s_i},\alpha_2$ are suitably small from \cref{wave strength,inimass}. Note that the right-hand side of \cref{1,2} are all good terms, so for convenience, we adopt a universal notation to denote such good terms throughout the remaining content of this work, which serves as the {\bf remainder terms} in fluid-type equations. Indeed, we denote the function $q$ to belong to the set $Q$, which may vary depending on the position, where the set $Q$ is given by  
\begin{align}\label{DefQ}
	Q=\Big\{q(t, x) :|q| \leq C\big(\delta^2+\left|\alpha_2\right| \delta^{\frac{3}{2}}\big) e^{-c(|x_1|+t)}+C \frac{\left|\alpha_2\right|}{(1+t)^{\frac{3}{2}}} e^{-\frac{\mu\left(x_1-u_{1\#} t\right)^2}{1+t}}+C\left(\delta+\left|\alpha_2\right|\right) e^{-c(|x_1|+t)}\Big\} . 
\end{align}
 By \cref{1,2} and direct calculations, one has from \cref{tirhomE1} that 
\begin{align}\label{u1}\notag
	\tilde{u}_1:=\frac{\tilde{m}_1}{\tilde{\rho}}&=  u_1^{s_1}+u_1^{s_1}-u_{1 \#}+\frac{1}{\tilde{\rho}}\big[-a \px\Theta+\left(\rho_{\#}-\rho^{s_3}-\Theta\right)\left(u_1^{s_1}-u_{1 \#}\right)\\
	&\notag\quad +\left(\rho_{\#}-\rho^{s_1}-\Theta\right)\left(u_1^{s_3}-u_{1 \#}\right)\big]\\
 &=u_1^{s_1}+u_1^{s_3}-u_{1 \#}-\frac{a}{\rho_{\#}} \px\Theta+q,
\end{align}
and hence, 
\begin{align}\label{m1}
 \notag
\frac{\tilde{m}_1^2}{\tilde{\rho}}&=\big(u_1^{s_1}+u_1^{s_3}-u_{1 \#}-\frac{a}{\rho_{\#}} \px\Theta+q\big)(m^{s_1}_1+m^{s_3}_1-m_{1\#}+u_{1\#}\Theta-a\px\Theta)\\
&=u_1^{s_1}m_1^{s_1}+u_1^{s_3}m_1^{s_3}+u_{1 \#}^2 \Theta-2 u_{1 \#} a \px\Theta+q.
\end{align}
Moreover, for terms involving $\tilde{E}$, one has,
\begin{align}
	\begin{aligned}
	\frac{\tilde{E}}{\tilde{\rho}}=\frac{E^{s_1}}{\rho^{s_1}}+\frac{E^{s_3}}{\rho^{s_3}}-\frac{E_\#}{\rho_{\#}}-\frac{\theta_{\#}}{\rho_{\#}}\Theta+\frac{\theta_{\#}}{\rho_{\#}^2}\Theta^2+\frac{g_3}{{\rho_{\#}}}+q.
\end{aligned}
\end{align}
Similarly,
\begin{align}\label{Q1}
	\begin{aligned}
	\frac{\tilde{E}\tilde{m}_1}{\tilde{\rho}}&=E^{s_1}u^{s_1}_1+E^{s_3}u^{s_3}_1-E_{\#}u_{1\#}+\frac{1}{2}u^3_{1\#}\Theta-\frac{E_{\#}a}{\rho_{\#}}\px\Theta+u_{1\#}g_3+q,\\
	\tilde{u}_1^2&:=\frac{\tilde{m}_1^2}{\tilde{\rho}^2}=(u^{s_1}_1)^2+(u^{s_3}_1)^2-(u_{1\#})^2-2\frac{u_{1\#}a}{\rho_{\#}}\px\Theta+q,\\
  	\tilde{\theta}&:=\frac{\tilde{E}}{\tilde{\rho}}-\frac{1}{2}|\tilde{u}|^2=\theta^{s_1}+\theta^{s_3}-\theta_{\#}-\frac{\theta_{\#}}{\rho_{\#}}\Theta+\frac{\theta_{\#}}{\rho_{\#}^2}\Theta^2+\frac{g_3}{\rho_{\#}}+\frac{u_{1\#}a}{\rho_{\#}}\px\Theta+q,\\
	\tilde{p}&:=\frac{2}{3}\tilde{\rho}\tilde{\theta}=\frac{2}{3}\big(\tilde{E}-\frac{1}{2}\rhot\abs{\tilde{u}}^2\big)=p^{s_1}+p^{s_3}-p_{\#}+\frac{2}{3}\big(g_3+au_{1\#}\px\Theta\big)+q,\\
	\frac{\tilde{p}\tilde{m}_1}{\tilde{\rho}}&=\tilde{p}\tilde{u}_1=p^{s_1}u_1^{s_1}+p^{s_3}u_1^{s_3}-p_{\#}u_{1\#}+\frac{2a}{3}\theta_{\#}\px\Theta+\frac{2u_{1\#}}{3}\big(g_3+au_{1\#}\px\Theta\big)+q,\\
	\kappa(\tilde{\T})\px\Tt&=\kappa(\T^{s_1})\px\T^{s_1}+\kappa(\T^{s_3})\px\T^{s_3}-\frac{\kappa(\theta_{\#})\theta_{\#}}{\rho_{\#}}\px\Theta+q.
\end{aligned}
\end{align}
Consequently, by directly substituting \cref{u1}-\cref{Q1} into the fluid-type equation for shock profiles \cref{sw-fluid}, one has 
\begin{align}\label{eq-ansat}
\left\{\begin{aligned}
	&\p_t \tilde{\rho}+\px\tilde{m}_{1}=0, \\
	& \p_t\tilde{m}_{1}+\px\left(\frac{\tilde{m}_1^2}{\tilde{\rho}}+\tilde{p}\right)=\px\left(\frac{4}{3}\mu(\tilde{\theta}) \px\tilde{u}_{1}\right)-\int v_1^2\left(\px\Pi^{s_1}+\px\Pi^{s_3}\right) d v+\px \tilde{Q}_{1}, \\
	& \p_t\tilde{E}+\px\left(\frac{\tilde{E} \tilde{m}_1}{\tilde{\rho}}+\frac{\tilde{p} \tilde{m}_1}{\tilde{\rho}}\right)=\px\left(\kappa(\tilde{\theta}) \px\tilde{\theta}\right)+\frac{4}{3}\px\left(\mu(\tilde{\theta}) \tilde{u}_1 \px\tilde{u}_{1}\right) \\
	& \quad\qquad-\int v_1 \frac{|v|^2}{2}\left(\px\Pi^{s_1}+\px\Pi^{s_3}\right) d v+\px \tilde{Q}_{2},
\end{aligned}\right.
\end{align}
where
\begin{align*}
	\tilde{Q}_1&=  \Big(\frac{\tilde{m}_1^2}{\tilde{\rho}}-\frac{\left(m_1^{s_1}\right)^2}{\rho^{s_1}}-\frac{\left(m_1^{s_3}\right)^2}{\rho^{s_3}}+\frac{m_{1 \#}^2}{\rho_{\#}}\Big)+\Big(\tilde{p}-p^{s_1}-p^{s_3}+p_{\#}\Big) \\
	&\quad -\frac{4}{3}\Big(\mu(\tilde{\theta}) \px\tilde{u}_{1}-\mu\left(\theta^{s_1}\right) \px u_{1}^{s_1}-\mu\left(\theta^{s_3}\right) \px u_{1}^{s_3}\Big) \\
	&\quad +2 u_{1 \#} a \px\Theta-u_{1 \#}^2 \Theta-a^2\px^2 \Theta,\\
	\tilde{Q}_2&=  \Big(\frac{\tilde{m}_1 \tilde{E}}{\tilde{\rho}}-\frac{m_1^{s_1} E^{s_1}}{\rho^{s_1}}-\frac{m_1^{s_3} E^{s_3}}{\rho^{s_3}}+\frac{m_{1 \#} E_{\#}}{\rho_{\#}}\Big) \\
	&\quad +\Big(\frac{\tilde{m}_1 \tilde{p}}{\tilde{\rho}}-\frac{m_1^{s_1} p^{s_1}}{\rho^{s_1}}-\frac{m_1^{s_3} p^{s_3}}{\rho^{s_3}}+\frac{m_{1 \#} p_{\#}}{\rho_{\#}}\Big) \\
	&\quad -\left(\kappa(\tilde{\theta}) \px\tilde{\theta}-\kappa\left(\theta^{s_1}\right) \px\theta^{s_1}-\kappa\left(\theta^{s_3}\right) \px\theta^{s_3}\right) \\
	&\quad -\frac{4}{3}\left(\mu(\tilde{\theta}) \tilde{u}_1 \px\tilde{u}_{1}-\mu\left(\theta^{s_1}\right) u^{s_1}\px u_{1}^{s_1}-\mu\left(\theta^{s_3}\right) u^{s_3} \px u_{1}^{s_3}\right) \\
	&\quad +\frac{1}{2} u_{1 \#}^2 a \px\Theta-\frac{1}{2} u_{1 \#}^3 \Theta+\p_tg_3.
\end{align*}
Similar to \cref{1,2}, using \cref{m1,Q1}, the
direct calculation yields $\tilde{Q}_1\in Q$. If we want to obtain $\tilde{Q}_2\in Q$, then by the expansion in \eqref{Q1}, the last term in $\tilde{Q}_2$ necessarily satisfies
\begin{align}
-\frac{E_{\#}a}{\rho_{\#}}\px\Theta+u_{1\#}g_3-\frac{2a}{3}\theta_{\#}\px\Theta+\frac{2}{3}\big(g_3+au_{1\#}\px\Theta\big)+\frac{1}{2} u_{1 \#}^2 a \px\Theta+\p_tg_3+\frac{\kappa(\theta_{\#})\theta_{\#}}{\rho_{\#}}\px\Theta=q.
\end{align}
Letting $g_3=-au_{1\#}\px\Theta$, one has
\begin{align}	-\big(\frac{E_{\#}}{\rho_{\#}}-\frac{1}{2}u_{1\#}^2\big)a\px\Theta-\frac{2a}{3}\theta_{\#}\px\Theta+\frac{\kappa(\theta_{\#})\theta_{\#}}{\rho_{\#}}\px\Theta=q,
\end{align}
and hence, we should take $a=\frac{3\kappa(\theta_{\#})}{5\rho_{\#}}$. Thus, the suitable ansatz is
\begin{align}\label{tirhomE}\left\{
\begin{aligned}
 \tilde{\rho}(x_1,t)&={\rho}^{s_1}(x_1-s_1t+\alpha_1)+{\rho}^{s_3}(x_1-s_3t+\alpha_3)-\rho_{\#}+\Theta(x_1, t),\\
 \tilde{m}_{1}(x_1, t)&={m}_1^{s_1}(x_1-s_1t+\alpha_1)+{m}_1^{s_3}(x_1-s_3t+\alpha_3)-m_{1\#}\\
&\quad+u_{1\#} \Theta(x_1, t)-a \px\Theta(x_1, t),\\
\tilde{m}_2&=\tilde{m}_3=0,\\
\tilde{E}(x_1, t)&={E}^{s_1}(x_1-s_1t+\alpha_1)+{E}^{s_3}(x_1-s_3t+\alpha_3)-E_{\#}\\
&\quad+\frac{1}{2} u_{1 \#}^{2} \Theta(x_1, t)-a u_{1 \#} \px\Theta(x_1, t).
\end{aligned}\right.
\end{align}
Then by a direct calculation as \cref{Q1}, one has
\begin{align}
    \tilde{\theta}&:=\frac{\tilde{E}}{\tilde{\rho}}-\frac{1}{2}|\tilde{u}|^2=\theta^{s_1}+\theta^{s_3}-\theta_{\#}-\frac{\theta_{\#}}{\rho_{\#}}\Theta+\frac{\theta_{\#}}{\rho_{\#}^2}\Theta^2+q.
\end{align}
Using \cref{initialal}, \cref{tirhomE1}, and the construction of $g_i$ above, we deduce that $\big(\tilde{\rho},\tilde{m}_1,\tilde{E}\big)$ satisfies \cref{zm}, and one also has the fluid-type equation: 
\begin{align}
    \label{tilderho}
 \left\{\begin{aligned}
&\tilde{\rho}_t+\px\left(\tilde{\rho} \tilde{u}_1\right)=0 \\[2mm]
&\tilde{\rho} \tilde{u}_{1t}+\tilde{\rho} \tilde{u}_1\px\tilde{u}_1+\frac{2}{3} \px\big(\tilde{\rho} \tilde{\theta}\big)=\px(\frac{4}{3} \mu(\Tt)\px\tilde{u}_{1})-\int v_1^2\left(\px\Pi^{s_1}+\px\Pi^{s_3}\right) d v+\px \tilde{Q}_{1}, \\[2mm]
&\tilde{\rho} \tilde{\theta}_t+\tilde{\rho} \tilde{u}_1 \px\tilde{\theta}+\frac{2}{3} \tilde{\rho} \tilde{\theta} \px\tilde{u}_{1}=\px\left(\kappa(\tilde{\theta}) \px\tilde{\theta}\right)+\frac{4}{3}\px\left(\mu(\tilde{\theta}) \tilde{u}_1 \px\tilde{u}_{1}\right) \\
	& \quad\qquad\quad\qquad-\int v_1 \frac{|v|^2}{2}\left(\px\Pi^{s_1}+\px\Pi^{s_3}\right) d v+\px \tilde{Q}_{2}-\tilde{u}_1\px\tilde{Q}_1.
\end{aligned}\right.
\end{align}

\subsection{Main theorem}
To present the main result, we first denote the macroscopic, microscopic, and full perturbations around the ansatz by
\begin{align}\label{perturbation}
\left\{\begin{aligned}
&(\phi, \psi, \omega)(t, x)=(\rho-\tilde{\rho}, m-\tilde{m}, E-\tilde{E})(t, x),\\
	& \tilde{\mathbf{G}}(t, x, v)=\mathbf{G}(t, x, v)-\mathbf{G}^{s_1}\left(x_1-s_1t+\alpha_1, v\right)-\mathbf{G}^{s_3}\left(x_1-s_3t+\alpha_3, v\right), \\
	& \tilde{f}(t, x, v)=f(t, x, v)-F^{s_1}\left(x_1-s_1t+\alpha_1, v\right)-F^{s_3}\left(x_1-s_3t+\alpha_3, v\right)+\mathbf{M}_{\#},
\end{aligned}\right.
\end{align}
where $(\tilde{\rho},\tilde{m},\tilde{E})$ is given by \cref{tilderho} and $\mathbf{M}_{\#}=\frac{\rho_{\#}}{\sqrt{(2 \pi R \theta_{\#})^3}} e^{-\frac{|v-u_{\#}|^2}{2 R \theta_{\#}}}$. The anti-derivatives in the $x_1$-axis are also well defined by first integrating the transverse direction $\Torus^2$: 
\begin{align}\label{Phi}
	\big(\Phi,\Psi_1,W\big)(x_1,t):=\int_{-\infty}^{x_1}\int_{\Torus^2}(\phi, \psi_1, \omega)(t, y_1,x_2,x_3)\,dx_2dx_3dy_1.
\end{align}
To capture the viscous effect of velocity and temperature, following \cite{LYYZ}, we set $(\tilde{\Psi}_1,\tilde{W})$ by 
\begin{align}\label{tiPsi}\begin{aligned}
	&\Psi_1=  \tilde{\rho} \tilde{\Psi}_1+ \tilde{u} \Phi,\\
	&W=  \tilde{\rho} \tilde{W}+ \tilde{u}_1  \Psi_1+\Big(\tilde{\theta}-\frac{| \tilde{u}|^{2}}{2}\Big) \Phi=  \tilde{\rho} \tilde{W}+  \tilde{\rho}  \tilde{u}_1 \tilde{\Psi}_1+\Big(\tilde{\theta}+\frac{| \tilde{u}|^{2}}{2}\Big) \Phi.
\end{aligned}
\end{align}
Then we denote the ``zeroth-derivative'' perturbations as 
\begin{align}\label{tipsi}
(\tilde{\psi}_1, \tilde{\omega})(t, x_1)=\big(\px\tilde{\Psi}_1, \px\tilde{W}\big)(t, x_1),\qquad (\varphi,\zeta):=(u-\tilde{u},\theta-\Tt).
\end{align}
To state the \emph{a priori} estimates more clearly, we introduce the instant energy functional $\mathcal{E}(t)$ and the dissipation energy functional $\mathcal{D}(t)$ as 
\begin{align}
\begin{aligned}
\mathcal{E}(t):=&\big\|\big(\Phi, \tilde{\Psi}_1, \tilde{W}\big)(t, \cdot)\big\|_{L^2_x}^2+\norm{(\phi, \psi, \zeta)(t, \cdot)}_{H^2_x}^2\\
 &+\sum_{0\leq \left|\alpha\right|\leq2} \iint \frac{\big|\partial^{\alpha} \tilde{\mathbf{G}}\big|^2}{\mathbf{M}_*} d v d x+\sum_{|\alpha|=3} \iint \frac{\big|\partial^\alpha \tilde{f}\big|^2}{\mathbf{M}_*} d v d x,\\
\mathcal{D}(t):=&\sum_{|\alpha|=1}\big\|\p^{\alpha}\big({\Phi},\tilde{\Psi}_1, \tilde{W}\big)\big\|_{L_x^2}^2+\big\|\sqrt{\left|\px u^{s_1}\right|+\left|\px\Theta\right|}(\Phi, \tilde{\Psi}_1, \tilde{W})\big\|_{L_x^2}^2 \\
& +\sum_{1 \leq\left|\alpha\right| \leq 3}\left\|\partial^{\alpha}\left(\phi, \psi, \omega\right)\right\|_{L_x^2}^2+ \sum_{1\leq |\beta'|\leq 3}\iint \frac{\nu(v)\big|\partial^{\beta} \tilde{\mathbf{G}}\big|^2}{\mathbf{M}_*}\, dvdx, 
		\end{aligned}
\end{align}
with some global Maxwellian $\M_*=\M_{[\rho_*,u_*,\theta_*]}$. 
Then we give the \emph{a priori} assumptions: 
\begin{align}\label{priori}
{\sup _{0 \leq t \leq T}\mathcal{E}(t)} \leq \chi^2,
\end{align}
where $\chi>0$ is a small positive constant depending on the initial data but independent of the time $T$. 
If we denote 
\begin{align}
    \mathcal{I}(0):=\big\|f_0(x, v)-\mathbf{M}_{\left[\tilde{\rho}\left(0, x_1\right), \tilde{u}\left(0, x_1\right), \tilde{\theta}\left(0, x_1\right)\right]}\big\|_{H_x^3(L
_v^2(\frac{1}{\sqrt{M_*}}))}, 
\end{align}
then by \cref{sys-perturb} and \cref{sobolev} below, one can easily verify the initial ``equivalence"
\begin{align*}
\mathcal{E}(0)=O(1)(\mathcal{I}(0)+\delta_0)\quad \text { and }\quad \mathcal{I}(0)=O(1)(\mathcal{E}(0)+\delta_0),
\end{align*} 
and the estimate on time derivative:
$$
\sup _{(\tau, x) \in[0, t] \times \mathbb{D}} \sum_{0 \leq\left|\alpha\right| \leq 1}\Big\{\big|\partial^{\alpha}(\phi, \psi, \zeta)\big|+\Big(\int \frac{\big|\partial^{\alpha} \mathbf{G}\big|^2}{\mathbf{M}_*} d v\Big)^{\frac{1}{2}}\Big\} \leq C(\chi+\delta_0).
$$

\smallskip 
Now we are ready to state the main result.
\begin{Thm}\label{mt}
Let $\left(\rho_{\pm}, m_{\pm}, E_{\pm}\right)$ be any two constant states satisfying \cref{Ex,Rh1,Rh2,lax,wave strength-same order}, and $\bar U=(\bar{\rho}, \bar{m}, \bar{E})^t$ be the corresponding superposition shock profiles given by \cref{shock wave}. Then there exist positive constants $\delta_0$, $\varepsilon_0$, $\eta_0$ and a global Maxwellian $\M_*=\M_{[\rho_*,u_*,\theta_*]}$ such that if 
\begin{align}
	\begin{aligned}
	&\abs{\left(\rho_{+}-\rho_{-}, m_{+}-m_{-}, E_{+}-E_{-}\right)}\leq \delta_0,\\
	&\ \ \rho_*>0,\ \ \frac{1}{2} \theta(t, x)<\theta_*<\theta(t, x), \\
	&|\rho(x, t)-\rho_*|+|u(x, t)-u_*|+|\theta(x, t)-\theta_*|<\eta_0,
		\end{aligned}
\end{align}
and the initial data satisfies that
\begin{align}\label{initialdata}
\begin{aligned}
	& \Big\{\|(\Phi, \Psi, W)\|_{L_x^2}^2+\|\tilde{f}\|_{H_x^3(L_v^2(\frac{1}{\sqrt{\mathbf{M}_*}}))}+\big\|\tilde{\mathbf{G}}\big\|_{H_x^2(L_v^2(\frac{1}{\sqrt{\mathbf{M}_*}}))}\Big\}\Big|_{t=0}+\Big|\int_{\mathbb{D}}(U(x, 0)-\bar{U}(x, 0))\,dx\Big| \leq \varepsilon_0.
\end{aligned}
\end{align}
Then the Cauchy problem \cref{be,ini} admits a unique global-in-time solution $f(t, x, v)$ satisfying energy estimate
	\begin{align}\label{priesti}
		\sup_{0\leq t\leq T}\mathcal{E}(t)+\int_0^T\mathcal{D}(t)dt\leq \mathcal{E}(0)+\delta_0^{\frac{1}{2}}, 
	\end{align}  
 for any $T>0$, 
 and asymptotic behavior 
\begin{align}\label{asympt}
\lim _{t \rightarrow \infty}\left\|f(t, x, v)-\left[F^{s_1}\left(x_1-s_1t+\alpha_1, v\right)+F^{s_3}\left(x_1-s_3t+\alpha_3, v\right)-\mathbf{M}_{\#}\right]\right\|_{L_x^{\infty} L_v^2(\frac{1}{\mathbf{M}_*})}=0. 
\end{align}
\end{Thm}
The proof of the asymptotic behavior \eqref{asympt} will be given in \cref{subproof} based on the energy estimate \cref{priesti}, whose proof is given in \cref{Sec4}. 

\begin{Rem}
	Our result is the first one concerning the time-asymptotic stability of a composite wave of two planar viscous shock waves to the Boltzmann equation with general 3-D initial perturbation without the macroscopic zero-mass condition in Eulerian coordinates. That is, we do not need the zero-mass type condition
 \begin{align}
    \int_{\mathbb{D}}(U-\bar{U})(x, 0) \,d x=0,
\end{align}
 as in Yuan \cite{Yuan}. Moreover, we consider the non-isentropic case and constant shifts, which is essentially different from Wang-Wang \cite{WW}. In addition, the stability of the shock profile for the multi-dimensional Boltzmann equation is proposed as an open problem by Wang-Wang \cite{WW1}. 
\end{Rem}

\section{Preliminaries and formulation of problem}\label{Sec3}
In this Section, we present some useful lemmas and properties for the decomposition of zero and non-zero modes. Then we introduce the perturbed system and perform a transformation to formulate the problem.

\subsection{Useful Lemmas}
Based on the celebrated H-theorem, we have the following lemmas about linearized collision operator in weighted $L^2$ space, which can be found in  \cite{LYY,LYYZ,LYs}. 
\begin{Lem}[{\cite[Lemma B.1]{LYs}}]
\label{fg}
There exists a positive constant $C$ such that
$$
\int \frac{\nu(v)^{-1} Q(f, g)^2}{\widetilde{\M}} d v \leq C\left\{\int \frac{\nu(v) f^2}{\widetilde{\M}} d v \cdot \int \frac{g^2}{\widetilde{\M}} d v+\int \frac{f^2}{\widetilde{\M}} d v \cdot \int \frac{\nu(v) g^2}{\widetilde{\M}} d v\right\},
$$
where $\widetilde{\M}$ can be any Maxwellian so that the above integrals are well-defined.
\end{Lem} 
\begin{Lem}[{\cite[Lemma 4.2--Lemma 4.4]{LYY}}]
\label{linear}
 If $\theta / 2<\theta_*<\theta$, then there exist two positive constants $\tilde{\sigma}=\tilde{\sigma}\left(\rho, u, \theta ; \rho_*, u_*, \theta_*\right)$ and $\eta_0=\eta_0\left(\rho, u, \theta ; \rho_*, u_*, \theta_*\right)$ such that if $\left|\rho-\rho_*\right|+\left|u-u_*\right|+\left|\theta-\theta_*\right|<\eta_0$, $\tilde{\sigma}>\tilde{c}>0$, where $\tilde{c}$ is a constant, then for $g(v) \in \mathfrak{N}^{\perp}$ where $\mathfrak{N}$ is the null space of $\LL_{\M}$, one has the following estimates
\begin{align*}
\text{1)}&\qquad -\int \frac{g \LL_{\mathbf{M}} g}{\mathbf{M}_*} d v \geq \tilde{\sigma} \int \frac{\nu(v) g^2}{\mathbf{M}_*} d v,\\
\text{2)}&\qquad \int \frac{\nu(v)}{\mathbf{M}}\left|\LL_{\mathbf{M}}^{-1} g\right|^2 d v \leq \tilde{\sigma}^{-2} \int \frac{\nu(v)^{-1} g^2}{\mathbf{M}} d v,\\
\text{3)}&\qquad \int \frac{\nu(v)}{\mathbf{M}_*}\left|\LL_{\mathbf{M}}^{-1} g\right|^2 d v \leq \tilde{\sigma}^{-2} \int \frac{\nu(v)^{-1} g^2}{\mathbf{M}_*} d v .\\
\end{align*}
Furthermore, for any positive constants $k$ and $\lambda$, it holds that
$$
\text{4)}\qquad\left|\int \frac{g_1 \mathbf{P}_1\left(|v|^k g_2\right)}{\mathbf{M}_*} d v-\int \frac{g_1|v|^k g_2}{\mathbf{M}_*} d v\right| \leq C_k \int \frac{\lambda\left|g_1\right|^2+\lambda^{-1}\left|g_2\right|^2}{\mathbf{M}_*} d v,
$$
where the constant $C_k$ only depends on $k$ and $\mathbf{P}_1$ is given by \eqref{pro}.
\end{Lem}
\begin{Lem}[Sobolev inequality, {\cite[Lemma 2.8]{WW1}}]\label{sobolev}
There exists some positive constant $C$ such that for $g \in H^2(\mathbb{D})$ with $\mathbb{D}:=\mathbb{R} \times \mathbb{T}^2$, it holds that
$$
\|g\|_{L^{\infty}(\mathbb{D})}^2 \leq C\left[\|g\|_{L^2(\mathbb{D})}\|\nabla g\|_{L^2(\mathbb{D})}+\|\nabla g\|_{L^2(\mathbb{D})}\left\|\nabla^2 g\right\|_{L^2(\mathbb{D})}\right] .
$$
\end{Lem} 
\subsection{The decomposition for zero and non-zero modes.}

In this Subsection, we will decompose the solution into the principal and transversal parts, which correspond to the zero and non-zero modes in Fourier space. Recalling the decomposition of $\mathbf{D}_0$ and $\mathbf{D}_{\neq}$,
\begin{align}\label{def-decom}
	\D_{0} h:=\mathring{h}:=\int_{\Torus^2}h\,dx_2dx_3,\qquad \D_{\neq} h:=\acute{h}:=h-\mathring{h},
\end{align}
for an arbitrary function $h$ which is integrable on $\Torus^2$.
By simple analysis, the following propositions of $\D_0$ and $\D_{\neq}$ hold for any suitable integrable function $h$, whose proof is basic and we omit it. 
\begin{Prop}\label{prop-decom}
	For the projections $\Do$ and $\Dn$ defined in \cref{def-decom}, the following holds,
\begin{itemize}
    \item  $\Do\Dn h=\Dn\Do h=0$.
     \item  For any nonlinear function $F\in C^2$, one has
	\begin{align}
		\Do F(U)-F(\Do U)=O(1)\D_0 (\Dn U)^2,
	\end{align}
	and similar results hold for $\tilde{U}$, $\bar{U}$, etc.
	
	\smallskip\item  $\|h\|^2_{L^2_x}=\|\Do h\|^2_{L^2_x}+\|\Dn h\|^2_{L^2_x}.$

 \smallskip\item The Poincar\'e inequality (with respect to $\Torus^2$) is valid for $\Dn h$:
 \begin{align}
     \norm{\Dn h}_{L^p_x}\leq C_p\norm{\nabla \Dn h}_{L^p_x},\quad \norm{\Dn h}_{L^p_x}\leq C_p\norm{\nabla  h}_{L^p_x}, \quad \text{for\ } p\geq 2,
 \end{align}
 where $C_p>0$ is a constant.
 \end{itemize}
\end{Prop}

\subsection{Formulation of problem}
In this Subsection, we give the formulation of some key quantities that play a crucial role in closing the a priori estimate. The most important of these is the estimate of the so-called anti-derivatives in the multi-dimensional case. 
Denote $(\Phi,\Psi,W)$ by \eqref{Phi}. 
Then by \cref{eq-ansat,NSG}, we have
\begin{align}\label{eq-pert1}
\left\{\begin{aligned}
&\pt\Phi+\px\Psi_{1}=0, \\ 
&\pt\Psi_{1}+\frac{2}{3}\left( \Em- \frac{| \mm|^{2}}{2 \rhom}-  {\tilde{E}}+ \frac{|{\tilde{m}}|^{2}}{2   \tilde{\rho}}\right)+ \frac{\mm_{1}^2}{ \rhom}- \frac{ {\tilde{m}}_{1}^2}{ \tilde{\rho}}=\frac{4}{3}\left(\mu(\Tm)\px\um_{1}-  \mu(\Tt)\px\tilde{u}_{1}\right)\\
&\qquad\qquad\qquad-\int v_1^2\left(\mr{\Pi}-\tilde{\Pi}\right) d v-Q_1, \\ 
&\pt W+\frac{5  \mm_1   \Em}{  3\rhom}-\frac{5 {\tilde{m}}_1   {\tilde{E}}}{  3\tilde{\rho}}-\frac{ |\mm|^{2}  \mm_1}{3 \rhom^2}+\frac{ {|{\tilde{m}}}|^{2}  {\tilde{m}}_1}{3\tilde{\rho}^2}= \left(\kappa(\mr{\theta}) \px\mr{\theta}-\kappa(\tilde{\theta}) \px\tilde{\theta}\right) \\ &\qquad+\frac{4}{3}\left(\mu(\mr{\theta}) \mr{u}_1 \px\mr{u}_{1}-\mu(\tilde{\theta}) \tilde{u}_1 \px\tilde{u}_{1}\right)+\sum_{i=2}^3 \mu(\mr\theta) \mr u_i \px\mr u_{i}-\frac{1}{2} \int v_1|v|^2\left(\mr{\Pi}- \tilde{\Pi}\right) d v- Q_2,
\end{aligned}\right.
\end{align}
where 
\begin{align*}
 &Q_1:=\tilde{Q}_1-\mathfrak{g}_1,\ Q_2:=\tilde{Q}_2-\mathfrak{g}_2,\\
&\mathfrak{g}_1:=(\mathring{\mathcal{G}}_2(U)-\mathcal{G}_2(\mathring{U}))\mathbb{I}_1,\ \mathfrak{g}_2:=\mathring{\mathcal{G}}_3(U)-\mathcal{G}_3(\mathring{U}).
\end{align*}
with $\mathcal{G}_2,\mathcal{G}_3$ given by \eqref{calG} and $\mathbb{I}_1=(1,0,0)$. 
By \cref{prop-decom}, one has
\begin{align}\label{gi}
\begin{aligned}
&\int_{\R}\abs{\mathring{\mathcal{G}}_2(U)-\mathcal{G}_2(\mathring{U})}dx_1\leq C{\norm{\Dn U}_{L^2_x}^2}\leq C\norm{\nabla(U-{\tilde{U}})}_{L^2_x}^2,\\
&\int_{\R}\abs{\mathring{\mathcal{G}}_2(U)-\mathcal{G}_2(\mathring{U})}dx_1\leq C{\norm{\Dn U}_{H^1_x}^2}\leq C\norm{\nabla(U-{\tilde{U}})}_{H^1_x}^2,
\end{aligned}
\end{align}
where $C$ is a positive constant independent of all the small parameters in this paper, such as $\eta_0,\delta,\varepsilon$, and so on. 
Moreover, denoting $(\tilde{\Psi}_1, \tilde{W})$ by \eqref{tiPsi}, then we have the following linearized system for $(\Phi, \tilde{\Psi}_1, \tilde{W})$:
\begin{align}\label{sys-Eu0}
	\left\{\begin{aligned}
		&\pt\Phi+  \tilde{\rho} \px\tilde{\Psi}_{1}+ \tilde{u}_{1} \px\Phi+  \px\tilde{\rho} \tilde{\Psi}_{1}+\px\tilde{u}_{1}\Phi=0 ,\\[2mm]
		  &\tilde{\rho} \pt\tilde{\Psi}_{1}+  \tilde{\rho}  \tilde{u}_{1} \px\tilde{\Psi}_{1}-\frac{1}{3} \tilde{\rho} \px\tilde{u}_{1} \tilde{\Psi}_{1}+\frac{2}{3}\left(   \px\tilde{\rho} \tilde{W}+  \tilde{\rho} \px\tilde{W}+ \tilde{\theta} \px\Phi-\frac{ \tilde{\theta}   \px\tilde{\rho}}{   \tilde{\rho}} \Phi\right) \\[2mm]
		&\qquad\qquad\qquad=\frac{4}{3}\mu(\Tt) \px^2\tilde{\Psi}_{1}-\int v_1^2\left(\mr{\Pi}-\tilde{\Pi}\right) d v+J_{1}+N_{1}-Q_{1}, \\[2mm]
		  &\tilde{\rho} \pt\tilde{W}+  \tilde{\rho}  \tilde{u}_{1} \px\tilde{W}-  \tilde{\rho}  \px\tilde{u}_{1} \tilde{W}+\frac{2}{3} \left(  \tilde{\rho} \tilde{\theta} \px\tilde{\Psi}_{1}-   \tilde{\rho} \px\tilde{\theta} \tilde{\Psi}_{1}\right)=\kappa(\Tt) \px^2\tilde{W} \\[2mm]
		&-\frac{1}{2} \int v_1|v|^2\left(\mr{\Pi}- \tilde{\Pi}\right) d v+\tilde{u}_1\int v_1^2\left(\mr{\Pi}-\tilde{\Pi}\right) d v+J_{4}+N_{4}-\left(Q_{2}-\tilde{u}_{1} Q_{1}\right),
	\end{aligned}\right.\end{align}
where
\begin{align*}
	J_{1}&=\left[\px\left(\int v_1^2\tilde{\Pi} d v-Q_{1}\right)- \frac{4}{3}\frac{\mu(\Tt)}{\tilde{\rho}}  \px\tilde{\rho}  \px\tilde{u}_{1}\right] \frac{\Phi}{ \tilde{\rho}} 
	+\frac{4}{3}\frac{\mu(\Tt)}{\tilde{\rho}}  \px\tilde{u}_{1} \px\Phi
	+\px\left(\frac{4}{3}\frac{\mu(\Tt)}{\tilde{\rho}}  \px\tilde{\rho} \tilde{\Psi}_{1}\right) \\
	&\quad+\frac{4}{3}\left(\left|\px\mr{\ut}_1\right|^2-\frac{\px\mr{\rhot}\px\Tt}{\mr{\rhot}}\right)\mu'(\Tt)\tilde{\Psi}_1+\frac{4}{3}\frac{\mu'(\Tt)}{\tilde{\rho}}\px\mr{\rhot}\px\mr{\ut}_1\tilde{W}+\frac{4}{3}\frac{\mu'(\Tt)}{\tilde{\rho}}\px\mr{\ut}_1\px\tilde{W} ,\\
	J_{4}&=\left(\frac{1}{2} \int v_1|v|^2\px\tilde{\Pi}d v-\tilde{u}_1\int v_1^2\px\tilde{\Pi} d v-\px Q_{2}+ \tilde{u}_{1} \px Q_{1}-\frac{\kappa(\Tt)}{ \tilde{\rho}}  \px\tilde{\rho} \px\tilde{\theta}\right) \frac{\Phi}{ \tilde{\rho}}\\
	&\quad+\frac{\kappa(\Tt)}{ \tilde{\rho}} \px\tilde{\theta} \px\Phi +\left[\px\left(\int v_1^2\tilde{\Pi} d v-Q_{1}\right)- \frac{4}{3}\frac{\mu(\Tt)}{\tilde{\rho}}  \px\tilde{\rho}  \px\tilde{u}_{1}\right]  \tilde{\Psi}_{1}\\
	&\quad+\frac{8}{3}\mu(\Tt)  \px\tilde{u}_{1} \px\tilde{\Psi}_{1} 
	+\px\left[\left(\kappa(\Tt)-\frac{4}{3}\mu(\Tt)\right)  \px\tilde{u}_{1} \tilde{\Psi}_{1}\right]+\px\left(\frac{\kappa(\Tt)}{ \tilde{\rho}}  \px\tilde{\rho} \tilde{W}\right),\\
			N_1:&=O(1)\left|\left(\px\Phi, \px\Psi_1,\px W, \px^2\Phi, \px^2\Psi_{1},\px^2 W\right)\right|^2, \\
	N_4:&=O(1)\left|\left(\px\Phi, \px\Psi_1,\px W, \px^2\Phi, \px^2\Psi, \px^2W\right)\right|^2.
\end{align*}
Next we derive the equation for the microscopic perturbation component $\tilde{\mathbf{G}}(t, x, v)$ and full perturbation $\tilde{f}$ given by \eqref{perturbation}. By \cref{g,Gsi}, we have
\begin{align}\label{gt}
\begin{aligned}
	\tilde{\mathbf{G}}_t-\LL_{\mathbf{M}} \tilde{\mathbf{G}}= & -\mathbf{P}_1\left(v\cdot \nabla\tilde{\mathbf{G}}\right)+Q(\tilde{\mathbf{G}}, \tilde{\mathbf{G}})+2 Q\left(\tilde{\mathbf{G}}, \mathbf{G}^{s_1}+\mathbf{G}^{s_3}\right)+2 Q\left(\mathbf{G}^{s_1}, \mathbf{G}^{s_3}\right) \\
	& -\left[\mathbf{P}_1\left(v\cdot\nabla \mathbf{M}\right)-\mathbf{P}_1^{s_1}\left(v\cdot\nabla \mathbf{M}^{s_1}\right)-\mathbf{P}_1^{s_3}\left(v\cdot\nabla \mathbf{M}^{s_3}\right)\right]+R_1+R_3,
\end{aligned}
\end{align}
where $R_i$ is given by 
\begin{align}\label{ri}
R_i=\left(\LL_{\mathbf{M}}-\LL_{\mathbf{M}^{s_i}}\right) \mathbf{G}^{s_i}-\left[\mathbf{P}_1\left(v\cdot\nabla \mathbf{G}^{s_i}\right)-\mathbf{P}_1^{s_i}\left(v\cdot\nabla \mathbf{G}^{s_i}\right)\right], \quad i=1,3 .
\end{align}
By \cref{MG,Fsi}, one has,
\begin{align}\label{ft}
\begin{aligned}
	{\tilde{f}_t+v\cdot\nabla\tilde{f}=} & \LL_{\mathbf{M}} \tilde{\mathbf{G}}+Q(\tilde{\mathbf{G}}, \tilde{\mathbf{G}})+\left(\LL_{\mathbf{M}}-\LL_{\mathbf{M}^{s_1}}\right)\left(\mathbf{G}^{s_1}\right)+\left(\LL_{\mathbf{M}}-\LL_{\mathbf{M}^{s_3}}\right)\left(\mathbf{G}^{s_3}\right) \\
	& +2 Q\left(\tilde{\mathbf{G}}, \mathbf{G}^{s_1}+\mathbf{G}^{s_3}\right)+2 Q\left(\mathbf{G}^{s_1}, \mathbf{G}^{s_3}\right).
\end{aligned}
\end{align}
\subsection{The local existence theorem and the a priori assumptions}
In this part, we will introduce the local existence theorem, which is similar to \cite[Lemma 2.1]{WW} and we omit the proof. Then, based on local existence, we set up the \emph{a priori} assumptions as in \eqref{priori}. 
\begin{Lem}[Local-in-time existence theorem]  \label{let}
For any suitable small constant $\varepsilon_1>0$, there exists a positive constant $T^*(\varepsilon_1)>0$, such that if the initial values $f_0(x, v) \geq 0$ and
\begin{align}\label{I0}
    \mathcal{I}(0):=\left\|f_0(x, v)-\mathbf{M}_{\left[\tilde{\rho}\left(0, x_1\right), \tilde{u}\left(0, x_1\right), \tilde{\theta}\left(0, x_1\right)\right]}\right\|_{H_x^3\left(L
_v^2\left(\frac{1}{\sqrt{M_*}}\right)\right)} \leq \frac{\varepsilon_1}{2 \sqrt{C_0}},
\end{align}
where $C_0:=\frac{1}{\min \{1, \tilde{\sigma}\}} \geq 1$ and the positive constant $\tilde{\sigma}$ is defined in \cref{linear}. Then the Cauchy problem \cref{be,ini} admits a unique solution on $[0,T(\varepsilon_1)] \times\mathbb{D} \times \mathbb{R}^3$ satisfying $f(t, x, v) \geq 0$ and
$$
\sup _{0 \leq t \leq T^{\star}(\varepsilon_1)}\left\|f(t, x, v)-\mathbf{M}_{\left[\rho\left(t, x_1\right), u\left(t, x_1\right), {\theta}\left(t, x_1\right)\right]}\right\|_{H_x^3\left(L_v^2\left(\frac{1}{\sqrt{M_*}}\right)\right)} \leq \varepsilon_1 .
$$
\end{Lem}

Then we can give our \emph{a priori} estimate as in the Main Theorem \ref{mt}. 
\begin{Prop}\label{ape}
	Under the same assumptions as \cref{mt}, the unique solution $(\phi,\psi,\zeta,\tilde{\G},\tilde{f})$ obtained in \cref{let},  satisfy  the energy estimate
	\begin{align}
		\sup_{0\leq t\leq T}\mathcal{E}(t)+\int_0^T\mathcal{D}(t)dt\leq \mathcal{E}(0)+\delta_0^{\frac{1}{2}}.
	\end{align}  
\end{Prop}

Using this energy estimate, we can obtain the global-in-time solution and asymptotic behavior \eqref{asympt}.

\subsection{Proof of the \cref{mt}}
\label{subproof}
First, by \cref{let} and \cref{ape}, one can obtain a global-in-time solution for \cref{be,ini} by a standard continuity argument. For the large-time behavior, we study the following estimates. By \cref{let} and \cref{ape}, 
\begin{align}\label{L1}
	& \notag\int_0^{+\infty}\norm{\nabla\tilde{f}}_{L^2_v(L^2_x(\frac{1}{\M_*}))}^2 d t \leq \int_0^{+\infty} \norm{\nabla\mathbf{M}-\nabla\mathbf{M}^{s_1}-\nabla\mathbf{M}^{s_3}}_{L^2_v(L^2_x(\frac{1}{\M_*}))}^2+\norm{\nabla\tilde{\mathbf{G}}}_{L^2_v(L^2_x(\frac{1}{\M_*}))}^2 d t \\
	&\notag\quad\leq C \int_0^{+\infty}\left\|\nabla(\phi, \psi, \omega)\right\|^2 d t+C \delta_0 \int_0^{+\infty}\|(\phi, \psi, \omega)\|^2 d t+\int_0^{+\infty} \iint \frac{\left|\nabla\tilde{\mathbf{G}}\right|^2}{\mathbf{M}_*} d v d x d t\\
	&\notag\qquad +C \int_0^{+\infty} \int\left[\left|\px\Theta\right|+\abs{\px u_{1}^{s_1}} \left(\left|u_1^{s_3}-u_{1 \#}\right|+|\Theta|\right)+\abs{\px u_{1 }^{s_3}} \left(\left|u_1^{s_1}-u_{1 \#}\right|+|\Theta|\right)+q\right]^2 d x d t \\
	&\notag\quad\leq  C \int_0^{+\infty}\left\|\nabla(\phi, \psi, \omega)\right\|^2 d t+C \delta_0 \int_0^{+\infty}\|(\phi, \psi, \omega)\|^2 d t+\int_0^{+\infty} \iint \frac{\left|{\nabla\tilde{\mathbf{G}}}\right|^2}{\mathbf{M}_*} d v d x d t+C \delta_0^2 \\
	&\quad\leq  C\left(\mathcal{E}(0)^2+\delta_0^{\frac{1}{2}}\right), 
\end{align}
where $q$ is function satisfying \cref{DefQ}. 
Similarly, by considering second-order derivative, one has
\begin{align}\label{BV}
	\int_{0}^\infty \Big|\frac{d}{dt}\norm{\nabla\tilde{f}}_{L^2_v(L^2_x(\frac{1}{\M_*}))}^2\Big|dt	\leq  C\left(\mathcal{E}(0)^2+\delta_0^{\frac{1}{2}}\right) .
\end{align}
Then by \cref{L1,BV}, we arrive at
\begin{align*}
	\norm{\nabla\tilde{f}}_{L^2_v(L^2_x(\frac{1}{\M_*}))}\rightarrow 0 \qquad \text{as} \qquad t \rightarrow \infty.
\end{align*}
Finally, by \cref{sobolev}, we obtain
\begin{align*}
\|\tilde{f}\|_{L^\infty_x(L^2_v(\frac{1}{\M_*}))}\leq C\|\tilde{f}\|_{L^2_x(L^2_v(\frac{1}{\M_*}))}^{\frac{1}{2}}\|\nabla_x\tilde{f}\|_{L^2_x(L^2_v(\frac{1}{\M_*}))}^{\frac{1}{2}}+C\|\nabla_x\tilde{f}\|_{L^2_x(L^2_v(\frac{1}{\M_*}))}^{\frac{1}{2}}\|\nabla_x^2\tilde{f}\|_{L^2_x(L^2_v(\frac{1}{\M_*}))}^{\frac{1}{2}}
\to 0,
\end{align*}
as $t\rightarrow\infty$. 
Then we have proved \cref{mt} provided that \cref{ape} holds. The rest of this paper will present the proof of \cref{ape}.

\section{The a priori estimate}\label{Sec4}
Due to the compressibility \cref{compressibility} of the shock profile, it is hard to control the lower order terms in the basic energy estimate. In this paper, we apply the method of anti-derivative and control these terms by using the integrated system.

\subsection{Lower-order estimates.}

To ensure clarity, we will first outline our proof strategy. Our goal is to prove the following fundamental estimates.

\subsubsection{Strategy of proof.}
In the energy estimation, the compressibility of the shock profile \cref{compressibility} holds great significance. Thus, for convenience, we use the notations:
\begin{align}
	\px u_1^s:=\px u^{s_1}_1+\px u^{s_3}_1<0,\quad { \px \lambda^s:=\px \lambda^{s_1}+\px \lambda^{s_3}}.
\end{align}
\begin{Lem}\label{be1}
Under the same assumptions as \cref{ape}, one has
\begin{multline*}
\sup_{0\le \tau\le t}\big\|\big(\px\Phi, \Phi, \Psi_1, W, \tilde{\Psi}_1, \tilde{W}\big)(\tau, \cdot)\big\|^2+\int_0^t\Big\|\sqrt{\abs{\px u_1^s}+\left|\px\Theta\right|}(\Phi, \tilde{\Psi}_1, \tilde{W})\Big\|^2d \tau\\
+\sum_{|\alpha^{\prime}|=1}\int_0^t\big\|\partial^{\alpha^{\prime}}(\Phi, \Psi_1, W, \tilde{\Psi}_1, \tilde{W})\big\|^2d \tau  
+\iint \frac{|\tilde{\mathbf{G}}|^2}{\mathbf{M}_*}(t, x, v) d v d x+\int_0^t \iint \frac{\nu(|v|)}{\mathbf{M}_*}|\tilde{\mathbf{G}}|^2 d v d x d \tau \\
\leq C\big(\mathcal{E}(0)^2+\delta_0^{\frac{1}{2}}\big) +C\int_0^t \sum_{\left|\alpha^{\prime}\right|=1}\left\|\partial^{\alpha^{\prime}}\left(\phi, \psi, \omega\right)\right\|_{L_x^2}^2+ \sum_{ |\beta'|=1}\iint \frac{\nu(v)|\partial^{\beta^{\prime}} \tilde{\mathbf{G}}|^2}{\mathbf{M}_*} d v d x d\tau.
\end{multline*}
\end{Lem}
\cref{be1} is the combination of the following Lemmas and estimate of the non-fluid part in \cref{secnonfluid}. First, by \cref{sys-Eu0}, one has the estimate of anti-derivative terms:
\begin{Lem}\label{be2}Under the same assumptions as \cref{ape}, one has,
	\begin{multline}
 \label{Eq42}
		\sup_{0\le \tau\le t}\|(\Phi, \tilde{\Psi}_1, \tilde{W})(\tau, \cdot)\|_{L_x^2}^2+\int_0^t\Big[\Big\|\sqrt{\abs{\px u_1^s}}(\tilde{\Psi}_1, \tilde{W})\Big\|_{L_x^2}^2+\Big\|\big(\px\tilde{\Psi}_1, \px\tilde{W}\big)\Big\|_{L_x^2}^2\Big] d \tau \\
  \leq C\delta_0\int_{0}^{t}\int_\R |\px{u}^{s}_{1}|\Phi^2dx_1d\tau+C(\delta+\chi+\varepsilon)\int_0^t\|\px(\Phi,\tilde{\Psi}_1,\tilde{W})\|^2_{H^1_x}d\tau +\int_0^t\mathcal{K}_1d\tau
 +C(\mathcal{E}(0)+\delta_0^{\frac{1}{2}}), 
	\end{multline}
where
\begin{align}\label{k1}\notag
		\mathcal{K}_1= & -\frac{{\tilde{\Psi}}_1}{\tilde{\theta}} \int v_1^2\left(\mr{\Pi}-\mr{\Pi}^{s_1}-\mr{\Pi}^{s_3}\right) d v\\
		& -\frac{{\tilde{W}}}{\tilde{\theta}^2}\left[\int v_1 \frac{|v|^2}{2}\left(\mr{\Pi}-\mr{\Pi}^{s_1}-\mr{\Pi}^{s_3}\right)d v-\tilde{u}_1 \int v_1^2\left(\mr{\Pi}-\mr{\Pi}^{s_1}-\mr{\Pi}^{s_3}\right)d v\right].
\end{align}
Here, $\Pi$ and $\Pi^{s_i}$ are given by \cref{the} and \cref{the1}, repsectively. 
\end{Lem}
Then we shall control $\mathcal{K}_1$ and obtain: 
\begin{Lem}\label{be3}
	Under the same assumptions as \cref{ape}, one has,
	\begin{multline*}
 \int_0^t\mathcal{K}_1d\tau\le C_\sigma\bigg[ \iint \frac{|\tilde{\mathbf{G}}|^2}{\mathbf{M}_*}(t, x, v) d v d x_1+ \int_0^t \iint \frac{\nu(|v|)}{\mathbf{M}_*}|\tilde{\mathbf{G}}|^2 d v d x_1 d \tau\bigg] \\+C\big(\mathcal{E}(0)^2+\delta_0^{\frac{1}{2}}\big)+
C(\delta_0+\chi+\sigma)\Big(\sup_{0\le \tau\le t}\|(\tilde{\Psi}_1,\tilde{W}_1)(\tau)\|_{L^2_x}^2+\int_0^t\mathcal{D}(\tau)d\tau\Big).
	\end{multline*}
 \end{Lem}
Next, we need to calculate the dissipation $\int_0^t\left\|\sqrt{\abs{\px u_1^s}} \Phi\right\|^2 d \tau$. To fully utilize the compressibility of the shock profile given in \cref{compressibility}, it is necessary to conduct the energy estimate in the diagonalized system.
\begin{Lem}\label{be4}
	Under the same assumptions as \cref{ape}, one has,
	\begin{multline*}
\sup_{0\le \tau\le t}\norm{(\Phi,\Psi_1,W,\tilde{\Psi}_1,\tilde{W})(\tau)}_{L^2_x}^2+\int_0^t\Big\|\sqrt{\abs{\px u_1^s}}(\Phi, \tilde{\Psi}_1, \tilde{W})\Big\|_{L^2_x}^2\\
		\leq C(\sigma+\delta_{0}^{\frac{1}{2}}+\chi)\int_0^t\mathcal{D}(\tau)d\tau+C\int_0^t\left\|\sqrt{\left|\px\Theta\right|}(\Phi, \tilde{\Psi}_1, \tilde{W})\right\|^2d\tau+C(\mathcal{E}(0)+\delta_0^{\frac{1}{2}})\\
 +C\iint \frac{|\tilde{\mathbf{G}}|^2}{\mathbf{M}_*}(t, x, v) d v d x+C_{\sigma}\int_0^t \iint \frac{\nu(|v|)}{\mathbf{M}_*}|\tilde{\mathbf{G}}|^2 d v d x d \tau.
	\end{multline*}
\end{Lem}
Finally, it remains to obtain the dissipation $\int_0^t\|\sqrt{|\px\Theta|}(\Phi, \tilde{\Psi}_1, \tilde{W})\|^2d\tau$: 
\begin{Lem}\label{be5}
		Under the same assumptions as \cref{ape}, one has,
	\begin{multline*}
  \sup_{0\le \tau\le t}\norm{(\px\Phi,\Phi,\Psi_1,W,\tilde{\Psi}_1,\tilde{W})(\tau)}_{L^2_x}^2+\int_0^t\Big\|\sqrt{\abs{\px u_1^s}+\abs{\px\Theta}}(\Phi, \tilde{\Psi}_1, \tilde{W})\Big\|^2_{L^2_x}\,d\tau
  \\+\sum_{|\alpha|=1}\int^t_0\|\partial^\alpha(\Phi,\tilde{\Psi}_1,\tilde{W})\|_{L^2_x}^2\,d\tau
			\leq C\iint \frac{|\tilde{\mathbf{G}}|^2}{\mathbf{M}_*}(t, x, v) d v d x\\+C_{\sigma}\int_0^t \iint \frac{\nu(|v|)}{\mathbf{M}_*}|\tilde{\mathbf{G}}|^2 d v d x d \tau+ C(\sigma+\delta_{0}^{\frac{1}{2}}+\chi)\int_0^t\mathcal{D}(\tau)d\tau+C(\mathcal{E}(0)+\delta_0^{\frac{1}{2}}).
	\end{multline*}
\end{Lem}

Next, we perform the proofs of the above lemmas. 

\subsubsection{Estimates for anti-derivatives}
\begin{proof}[Proof of \cref{be2}]

 Multiplying $\cref{sys-Eu0}_{1}$ by $\frac{2\Phi}{ 3 \tilde{\rho}},\cref{sys-Eu0}_{2}$ by $\frac{\tilde{\Psi}_{1}}{\tilde{\theta}}$, $\cref{sys-Eu0}_{3}$  by $\frac{\tilde{W}}{\tilde{\theta}^{2}}$, respectively, and then integrating the summation of the above-resulting equation over $\R\times[0,t]$, we can get
\begin{align}\label{eq-1}\notag
		&\int_{\R}\Big(\frac{\Phi^{2}}{3   \tilde{\rho}}+\frac{  \tilde{\rho} \tilde{\Psi}_{1}^{2}}{2 \tilde{\theta}}+\frac{  \tilde{\rho} \tilde{W}^{2}}{2 \tilde{\theta}^{2}}\Big)(t)dx_1+\int_0^t\int_{\R}\frac{4\mu(\Tt)}{3\tilde{\theta}} \left|\px\tilde{\Psi}_{1}\right|^{2}+\frac{\kappa(\Tt)}{\tilde{\theta}^{2}} \left|\px\tilde{W}\right|^{2}dx_1d\tau \\
		&\notag \quad -\int_0^t\int_{\R} \tilde{\rho} \px\tilde{u}_{1}\Big(\frac{\tilde{\Psi}_{1}^{2}}{3 \tilde{\theta}}+\frac{\tilde{W}^{2}}{\tilde{\theta}^{2}}\Big)-\Big(\frac{\tilde{\Psi}_{1}^{2}}{2 \tilde{\theta}^{2}}+\frac{\tilde{W}^{2}}{\tilde{\theta}^{3}}\Big)\left(  \tilde{\rho} \pt\tilde{\theta}+  \tilde{\rho} \tilde{u}_{1} \px\tilde{\theta}\right)dx_1d\tau\\
	   =&\notag\int_0^t\int_{\R}\px^2\Big[\frac{2\mu(\Tt)}{3\mr{\tilde{\T}}}\Big]\tilde{\Psi}_1^{2}+\px^2\Big[\frac{\kappa(\Tt)}{2\tilde{\theta}^{2}}\Big]\tilde{W}^2dx_1d\tau +\int_0^t\int_{\R}\frac{\tilde{\Psi}_{1}}{\tilde{\theta}}\left(J_{1}+N_{1}-Q_{1}\right)dx_1d\tau\\
		&\quad+\int_0^t\int_{\R}\big(\frac{\tilde{W}}{\tilde{\theta}^{2}}\left(J_{4}+N_{4}-Q_{2}+\tilde{u}_{1} Q_{1}\right)+\mathcal{K}_1\big)dx_1d\tau
  +\int_{\R}\Big(\frac{\Phi^{2}}{3   \tilde{\rho}}+\frac{  \tilde{\rho} \tilde{\Psi}_{1}^{2}}{2 \tilde{\theta}}+\frac{  \tilde{\rho} \tilde{W}^{2}}{2 \tilde{\theta}^{2}}\Big)(0)dx_1.
\end{align}
where $\mathcal{K}_1$ is defined in \cref{k1}.
 Firstly, to deal with the third term on the left-hand side of \eqref{eq-1}, we use the good sign \cref{compressibility} and equation \eqref{tilderho} to obtain
\begin{align*}
 \tilde{\rho} \tilde{\theta}_{t}+ \tilde{\rho}  \tilde{u}_{1} \px\tilde{\theta}=&-\frac{2}{3}  \tilde{\rho} \tilde{\theta}  \px\tilde{u}_{1}+\frac{4}{3}\mu(\Tt)\left( \px\tilde{u}_{1}\right)^{2}+\px\left(\kappa(\Tt)\px\tilde{\theta}\right)-\int v_1\frac{\abs{v}^2}{2}(\px\Pi^{s_1}+\px\Pi^{s_3})dv\\
&+\tilde{u}_1\int{\abs{v_1}^2}(\px\Pi^{s_1}+\px\Pi^{s_3})dv+\px Q_{2}- \tilde{u}_{1} \px Q_{1}.
\end{align*}
Note that 
\begin{equation*}
 \px\tilde{u}_{1}=\px u_1^s+O(1)\px\Theta=\px u_1^s+q,
\end{equation*}
and 
\begin{multline*}
\frac{4}{3}\mu(\Tt)\left( \px\tilde{u}_{1}\right)^{2}+\px\left(\kappa(\Tt)\px\tilde{\theta}\right)-\int v_1\frac{\abs{v}^2}{2}(\px\Pi^{s_1}+\px\Pi^{s_3})dv\\
+\tilde{u}_1\int{\abs{v_1}^2}(\px\Pi^{s_1}+\px\Pi^{s_3})dv+\px Q_{2}- \tilde{u}_{1} \px Q_{1}=O(1)\delta \abs{\px u^{s}_1}+q. 
\end{multline*}
we can write 
\begin{align*}
	&\int_{0}^{t}\int_\R  \px\tilde{u}_{1}\Big( \tilde{\Psi}_{1}^{2}+
 \tilde{W}^{2}\Big)dx_1d\tau
	=\int_{0}^{t}\int_\R \px{u}^{s}_{1}\Big( \tilde{\Psi}_{1}^{2}
 + \tilde{W}^{2}\Big)dx_1d\tau+\int_{0}^{t}\int_\R q\Big( \tilde{\Psi}_{1}^{2}
 + \tilde{W}^{2}\Big)dx_1d\tau. 
\end{align*}
Similarly, to deal with the first term on the right-hand side of \eqref{eq-1}, we have 
\begin{align*}
    \int_0^t\int_{\R}\px^2\Big[\frac{2\mu(\Tt)}{3\mr{\tilde{\T}}}\Big]\tilde{\Psi}_1^{2}+\px^2\Big[\frac{\kappa(\Tt)}{2\tilde{\theta}^{2}}\Big]\tilde{W}^2dx_1d\tau\leq \int_0^t\int_{\R}(\delta_0\abs{\px u_1^s}+\abs{q})\big(\tilde{\Psi}_1^{2}+\tilde{W}^2\big)dx_1d\tau, 
\end{align*}
and by \eqref{DefQ}, 
\begin{align}\label{es-er0}
	\int_{0}^{t}\int_\R \abs{q}\Big( \tilde{\Psi}_{1}^{2}
 + \tilde{W}^{2}\Big)dx_1d\tau\leq(\varepsilon_0+\delta_0)\sup_{0\leq\tau\leq t}\|(\tilde{\Psi}_1,\tilde{W})\|^2.
\end{align}
By direct calculations, we can obtain the estimate of $J_i$: 
\begin{multline*}
	J_i=O(1)\delta |\px u^{s}|(|\Phi|,|\tilde{\Psi}_1|,|\tilde{W}|)+O(1)|\px u^{s}|(|\px\Phi|,|\px\tilde{\Psi}_1|,|\px\tilde{W}|)\\
	\qquad\qquad+(\varepsilon e^{-\alpha t}+\delta^{\frac{3}{2}}e^{-\delta t})(\p_1+\p^2_1)(|\Phi|,|\tilde{\Psi}_1|,|\tilde{W}|),\qquad i=1,4.
\end{multline*}
Thus, for $i=1,4$,
\begin{multline*}
\int_0^t\int_\R J_i(\Phi,\tilde{\Psi}_1,\tilde{W})dx_1d\tau
\leq C_\epsilon\delta\int_0^{t}\int_\R |\px u^{s}_{1}|\left(|\Phi|^2+|\tilde{\Psi}_1|^2+|\tilde{W}|^2\right)dxd\tau\\
\qquad\qquad\qquad+C_\epsilon(\varepsilon+\delta)\sup_{0\leq\tau\leq t}\|(\tilde{\Phi},\tilde{\Psi}_1,\tilde{W})\|^2+\epsilon\int_0^t\|(\px\Phi,\px\tilde{\Psi}_1,\px\tilde{W})\|_{H^1_x}^2d\tau,
\end{multline*}
where $\epsilon$ is an arbitrary small parameter. 
For the nonlinear terms, we use the a priori assumption
\begin{align}\label{es-er1}
	\int_0^t\int_\R  N_{i}(\tilde{\Psi}_1,\tilde{W})dxd\tau\leq\chi\int_0^t\|\px(\Phi,\tilde{\Psi}_1,\tilde{W})\|^2_{H^1_x}d\tau, 
\end{align}
where we used \eqref{priori} to control $(\tilde{\Psi}_1,\tilde{W})$. 
Finally, substituting the above estimates in \cref{eq-1}, we have
\begin{multline*}
	\|(\Phi,\tilde{\Psi}_1, \tilde{W})(t)\|^{2}-c_1\int_{0}^{t}\int_\R \px{u}^{s}_{1}\Big( \tilde{\Psi}_{1}^{2}
 + \tilde{W}^{2}\Big)dx_1d\tau+c_2\int_0^t\Big(\|\px\tilde{\Psi}_1\|^2+\|\px\tilde{W}\|^2\Big)d\tau\\
	\leq C\delta\int_{0}^{t}\int_\R |\px{u}^{s}_{1}|\Phi^2dx_1d\tau+C(\delta+\chi+\varepsilon)\int_0^t\|\px(\Phi,\tilde{\Psi}_1,\tilde{W})\|^2_{H^1_x}d\tau +\int_0^t\mathcal{K}_1d\tau
 +C(\mathcal{E}(0)+\delta_0^{\frac{1}{2}}).
\end{multline*}
This completes the proof of \cref{be2}. 
\end{proof}
\subsubsection{Wave interaction}
We then estimate the interaction of fluid ($s_i$) and non-fluid parts ($\tilde{G}$), i.e. $\mathcal{K}_1$. 

\begin{proof}[Proof of \cref{be3}]
We shall first study the following 
\begin{align*}
	&\Pi-\Pi^{s_1}-\Pi^{s_3}=\LL_{\mathbf{M}}^{-1}\left[\pt\tilde{\mathbf{G}}+\mathbf{P}_1\left(v\cdot\nabla_x \tilde{\mathbf{G}}\right)-Q(\tilde{\mathbf{G}}, \tilde{\mathbf{G}})\right] \\
	&\quad -\LL_{\mathbf{M}}^{-1}\left[2 Q\left(\tilde{\mathbf{G}}, \mathbf{G}^{s_1}+\mathbf{G}^{s_3}\right)\right]-\LL_{\mathbf{M}}^{-1}\left[2 Q\left(\mathbf{G}^{s_1}, \mathbf{G}^{s_3}\right)\right] \\
	&\qquad +\left(\LL_{\mathbf{M}}^{-1}-\LL_{\mathbf{M}^{s_1}}^{-1}\right)\left[\pt\mathbf{G}^{s_1}-Q\left(\mathbf{G}^{s_1}, \mathbf{G}^{s_1}\right)\right]+\left(\LL_{\mathbf{M}}^{-1} \mathbf{P}_1-\LL_{\mathbf{M}^{s_1}}^{-1} \mathbf{P}_1^{s_1}\right)\left(v\cdot\nabla_x \mathbf{G}^{s_1}\right) \\
	&\qquad +\left(\LL_{\mathbf{M}}^{-1}-\LL_{\mathbf{M}^{s_3}}^{-1}\right)\left[\pt\mathbf{G}^{s_3}-Q\left(\mathbf{G}^{s_3}, \mathbf{G}^{s_3}\right)\right]+\left(\LL_{\mathbf{M}}^{-1} \mathbf{P}_1-\LL_{\mathbf{M}^{s_3}}^{-1} \mathbf{P}_1^{s_3}\right)\left(v\cdot\nabla_x \mathbf{G}^{s_3}\right) \\
	&\quad=  : \LL_{\mathbf{M}}^{-1}\left[\pt\tilde{\mathbf{G}} +\mathbf{P}_1\left(v\cdot\nabla_x \tilde{\mathbf{G}}\right)-Q(\tilde{\mathbf{G}}, \tilde{\mathbf{G}})\right]-\LL_{\mathbf{M}}^{-1}\left[2 Q\left(\tilde{\mathbf{G}}, \mathbf{G}^{s_1}+\mathbf{G}^{s_3}\right)\right] \\
	& \qquad -\LL_{\mathbf{M}}^{-1}\left[2 Q\left(\mathbf{G}^{s_1}, \mathbf{G}^{s_3}\right)\right]+\mathcal{V}_1+\mathcal{V}_3,
\end{align*}
where
\begin{align}
\mathcal{V}_i=\left(\LL_{\mathbf{M}}^{-1}-\LL_{\mathbf{M}^{s_i}}^{-1}\right)\left[\pt\mathbf{G}^{s_i}-Q\left(\mathbf{G}^{s_i}, \mathbf{G}^{s_i}\right)\right]+\left(\LL_{\mathbf{M}}^{-1} \mathbf{P}_1-\LL_{\mathbf{M}^{s_i}}^{-1} \mathbf{P}_1^{s_i}\right)\left(v\cdot\nabla_x \mathbf{G}^{s_i}\right), \quad i=1,3.
\end{align}
Next, we only estimate the first term in $\mathcal{K}_1$ (given in \cref{k1}) since the other terms are similar:
\begin{align*}
		\mathcal{K}_{11}= & -\int_0^t \int \frac{\tilde{\Psi}_1}{\tilde{\theta}} \int v_1^2 \LL_{\mathbf{M}}^{-1}(\tilde{\mathbf{G}})_t d v d x d \tau-\int_0^t \int \frac{\tilde{\Psi}_1}{\tilde{\theta}} \int v_1^2 \LL_{\mathbf{M}}^{-1}\left[\mathbf{P}_1\left(v\cdot\nabla_x \tilde{\mathbf{G}}\right)\right] d v d x d \tau \\
		& +\int_0^t \int \frac{\tilde{\Psi}_1}{\tilde{\theta}} \int v_1^2 \LL_{\mathbf{M}}^{-1}[Q(\tilde{\mathbf{G}}, \tilde{\mathbf{G}})] d v d x d \tau  +\int_0^t \int \frac{\tilde{\Psi}_1}{\tilde{\theta}} \int v_1^2 \LL_{\mathbf{M}}^{-1}\left[2 Q\left(\tilde{\mathbf{G}}, \mathbf{G}^{s_1}+\mathbf{G}^{s_3}\right)\right] d v d x d \tau \\
		& +\int_0^t \int \frac{\tilde{\Psi}_1}{\tilde{\theta}} \int v_1^2 \LL_{\mathbf{M}}^{-1}\left[2 Q\left(\mathbf{G}^{s_1}, \mathbf{G}^{s_3}\right)\right] d v d x d \tau \\
		& -\int_0^t \int \frac{\tilde{\Psi}_1}{\tilde{\theta}} \int v_1^2\left(\LL_{\mathbf{M}}^{-1}-\LL_{\mathbf{M}^{s_1}}^{-1}\right)\left[\mathbf{G}_t^{s_1}-Q\left(\mathbf{G}^{s_1}, \mathbf{G}^{s_1}\right)\right] d v d x d \tau \\
		& -\int_0^t \int \frac{\tilde{\Psi}_1}{\tilde{\theta}} \int v_1^2\left(\LL_{\mathbf{M}}^{-1}-\LL_{\mathbf{M}^{s_3}}^{-1}\right)\left[\mathbf{G}_t^{s_3}-Q\left(\mathbf{G}^{s_3}, \mathbf{G}^{s_3}\right)\right] d v d x d \tau \\
		& -\int_0^t \int \frac{\tilde{\Psi}_1}{\tilde{\theta}} \int v_1^2\left(\LL_{\mathbf{M}}^{-1} \mathbf{P}_1-\LL_{\mathbf{M}^{s_1}}^{-1} \mathbf{P}_1^{s_1}\right)\left(v\cdot\nabla_x \mathbf{G}^{s_1}\right) d v d x d \tau \\
		& -\int_0^t \int \frac{\tilde{\Psi}_1}{\tilde{\theta}} \int v_1^2\left(\LL_{\mathbf{M}}^{-1} \mathbf{P}_1-\LL_{\mathbf{M}^{s_3}}^{-1} \mathbf{P}_1^{s_3}\right)\left(v\cdot\nabla_x \mathbf{G}^{s_3}\right) d v d x d \tau=: \sum_{i=1}^9 \mathcal{K}_{11}^i.
\end{align*}
Firstly, note that the linearized operator $\LL_{\mathbf{M}}^{-1}$ satisfies (see for instance from \cite{XYY}), for any $g \in \mathfrak{N}^{\perp}$,
\begin{align}
    \label{LMminus1}
\p^{\alpha}\left(\LL_{\mathbf{M}}^{-1} g\right)=\LL_{\mathbf{M}}^{-1}\left(\p ^{\alpha}g\right)-2 \LL_{\mathbf{M}}^{-1}\left\{Q\left(\LL_{\mathbf{M}}^{-1} g, \p^{\alpha}\M\right)\right\}, \qquad |\alpha|=1.
\end{align}
Then we have
\begin{align}\label{K11a}
\notag
	\mathcal{K}_{11}^1= & -\int_0^t \int \frac{\tilde{\Psi}_1}{\tilde{\theta}} \int v_1^2\left(\LL_{\mathbf{M}}^{-1} \tilde{\mathbf{G}}\right)_t d v d x d \tau  -2 \int_0^t \int \frac{\tilde{\Psi}_1}{\tilde{\theta}} \int v_1^2 \LL_{\mathbf{M}}^{-1}\left[Q\left(\LL_{\mathbf{M}}^{-1} \tilde{\mathbf{G}}, \M_t\right)\right] d v d x d \tau \\
	= \notag& \iint\left(\frac{\tilde{\Psi}_1}{\tilde{\theta}} v_1^2 \LL_{\mathbf{M}}^{-1} \tilde{\mathbf{G}}\right)(0, x, v) d v d x-\iint\left(\frac{\tilde{\Psi}_1}{\tilde{\theta}} v_1^2 \LL_{\mathbf{M}}^{-1} \tilde{\mathbf{G}}\right)(t, x, v) d v d x \\
	& +\int_0^t \int\left(\frac{\tilde{\Psi}_1}{\tilde{\theta}}\right)_t \int v_1^2 \LL_{\mathbf{M}}^{-1} \tilde{\mathbf{G}} d v d x d \tau -2 \int_0^t \int \frac{\tilde{\Psi}_1}{\tilde{\theta}} \int v_1^2 \LL_{\mathbf{M}}^{-1}\left[Q\left(\LL_{\mathbf{M}}^{-1} \tilde{\mathbf{G}}, \mathbf{M}_t\right)\right] d v d x d \tau .
\end{align}
By \cref{fg} and \cref{linear}, one has 
\begin{align*}
    &\Big|\int v_1^2 \LL_{\mathrm{M}}^{-1} \tilde{\mathbf{G}} d v\Big|^2 \leq C \int \frac{\nu^{-1}(|v|)}{\mathbf{M}_*}|\tilde{\mathbf{G}}|^2 d v,
\end{align*}
and 
\begin{multline*}
	\int v_1^2 \LL_{\mathbf{M}}^{-1}\left[Q\left(\LL_{\mathbf{M}}^{-1} \tilde{\mathbf{G}}, \mathbf{M}_t\right)\right] d v  \leq C\left(\int \frac{\nu(|v|)}{\mathbf{M}_*}\left|\LL_{\mathbf{M}}^{-1}\left[Q\left(\LL_{\mathbf{M}}^{-1} \tilde{\mathbf{G}}, \mathbf{M}_t\right)\right]\right|^2 d v\right)^{\frac{1}{2}} \\
	\leq C\left(\int \frac{\nu(|v|)}{\mathbf{M}_*}\left|\LL_{\mathbf{M}}^{-1} \tilde{\mathbf{G}}\right|^2 d v\right)^{\frac{1}{2}} \cdot\left(\int \frac{\nu(|v|)}{\mathbf{M}_*}\left|\mathbf{M}_t\right|^2 d v\right)^{\frac{1}{2}} \leq C\left|\left(\rho_t, u_t, \theta_t\right)\right|\left(\int \frac{\nu^{-1}(|v|)}{\mathbf{M}_*}|\tilde{\mathbf{G}}|^2 d v\right)^{\frac{1}{2}}.
\end{multline*}
Then \cref{K11a} can be estimated as 
\begin{align}\label{k11b}
\notag\mathcal{K}_{11}^1 &\leq  \sigma\Big(\left\|\tilde{\Psi}_1(t, \cdot)\right\|^2+\int_0^t\left\|\tilde{\Psi}_{1 \tau}\right\|^2 d \tau\Big)+C_\sigma \iint \frac{|\tilde{\mathbf{G}}|^2}{\mathbf{M}_{*}}(t, x, v) d v d x \\
	\notag&\quad +C_\sigma \int_0^t \iint \frac{\nu(|v|)}{\mathbf{M}_*}|\tilde{\mathbf{G}}|^2 d v d x d \tau+C \delta_0 \int_0^t\Big\|\sqrt{\left|\px u_{1}^{s}\right|} \tilde{\Psi}_1\Big\|^2 d \tau \\
	\notag&\quad +C \chi \int_0^t\left\|\left(\phi_\tau, \psi_\tau, \omega_\tau\right)\right\|^2 d \tau+C \int_0^t \int q \tilde{\Psi}_1^2 d x d \tau+C \mathcal{E}(0)^2\\
	\notag&\leq C (\mathcal{E}(0)^2+\delta_{0}^{\frac{1}{2}})+C(\delta_0+\chi+\sigma)\bigg(\mathcal{E}(t)+\int_0^t\mathcal{D}(\tau)d\tau\bigg)\\
	&\qquad +C_{\sigma}\Big(\iint \frac{|\tilde{\mathbf{G}}|^2}{\mathbf{M}_{*}}(t, x, v) d v d x\int_0^t \iint \frac{\nu(|v|)}{\mathbf{M}_*}|\tilde{\mathbf{G}}|^2 d v d x d \tau\Big),
\end{align}
where $\sigma>0$ is a small constant to be determined, $C_\sigma>0$ is a constant depending on $\sigma$ and we have used the following estimate from \cref{sw} and \eqref{derismall}: 
\begin{align}
     \abs{\rhot_t,\ut_t,\Tt_t}^2\leq C\delta|\px\lambda^s|+|q|.
\end{align}
For the term $\mathcal{K}_{11}^2$, note that (see for instance \cite[(6.3)]{LYYZ})
$$
\mathbf{P}_1\left(v_1 \px\tilde{\mathbf{G}}\right)=\px\left[\mathbf{P}_1\left(v_1 \tilde{\mathbf{G}}\right)\right]+\sum_{j=0}^4\langle v_1 \tilde{\mathbf{G}}, \chi_j\rangle\mathbf{P}_1\left(\px\chi_{j }\right),
$$
then by \cref{linear} and \eqref{LMminus1}, one has
\begin{align}\label{k12}
\notag
	\mathcal{K}_{11}^2= & \int_0^t \int\px\left(\frac{\tilde{\Psi}_1}{\tilde{\theta}}\right) \int v_1^2 \LL_{\mathrm{M}}^{-1}\left[\mathbf{P}_1\left(v_1 \tilde{\mathbf{G}}\right)\right] d v d x_1 d \tau \\
	\notag& -\int_0^t \int \frac{\tilde{\Psi}_1}{\tilde{\theta}} \int v_1^2 \LL_{\mathrm{M}}^{-1}\Big(\sum_{j=0}^4\langle v_1 \tilde{\mathbf{G}}, \chi_j\rangle\mathbf{P}_1\left(\px\chi_{j}\right)\Big) d v d x_1 d \tau \\
	\notag& -2 \int_0^t \int \frac{\tilde{\Psi}_1}{\tilde{\theta}} \int v_1^2 \LL_{\mathrm{M}}^{-1}\Big\{Q\big(\LL_{\mathrm{M}}^{-1} \mathbf{P}_1(v_1 \tilde{\mathbf{G}}), \px\mathbf{M}\big)\Big\} d v d x d \tau \\
	\notag\leq & \sigma \int_0^t\big\|\px\tilde{\Psi}_{1}\big\|^2 d \tau+C \delta_0 \int_0^t\big\|\sqrt{|\px u_{1}^{s}|} \tilde{\Psi}_1\big\|^2 d \tau \\
	& +C \chi \int_0^t\left\|\left(\px\phi, \px\psi, \px\omega\right)\right\|^2 d \tau+C_\sigma \int_0^t \iint \frac{\nu(|v|)}{\mathbf{M}_*}|\tilde{\mathbf{G}}|^2 d v d x d \tau .
\end{align}
Similarly, we have 
\begin{align}\label{k13}
	\mathcal{K}_{11}^3 & =\int_0^t \int \frac{\tilde{\Psi}_1}{\tilde{\theta}}\left(\int \frac{\nu(|v|)}{\mathbf{M}_*}\left|\LL_{\mathbf{M}}^{-1}[Q(\tilde{\mathbf{G}}, \tilde{\mathbf{G}})]\right|^2 d v\right)^{\frac{1}{2}} d x d \tau  \leq C \chi \int_0^t \iint \frac{\nu(|v|)}{\mathbf{M}_*}|\tilde{\mathbf{G}}|^2 d v d x d \tau, 
\end{align}
and
\begin{align}\label{k14}\notag
	\mathcal{K}_{11}^4 &  \leq C \int_0^t \int|\tilde{\Psi}_1|\Big(\int \frac{\nu^{-1}(|v|)\big|Q\big(\tilde{\mathbf{G}}, \mathbf{G}^{s_1}+\mathbf{G}^{s_3}\big)\big|^2}{\mathbf{M}_*} d v\Big)^{\frac{1}{2}} d x_1 d \tau \\
	\notag& \leq C \int_0^t \int|\tilde{\Psi}_1|\Big(\int \frac{\nu(|v|) \mid \tilde{\mathbf{G}}^2}{\mathbf{M}_*} d v\Big)^{\frac{1}{2}}\Big(\int \frac{\nu(|v|)(|\mathbf{G}^{s_1}|^2+|\mathbf{G}^{s_3}|^2)}{\mathbf{M}_*} d v\Big)^{\frac{1}{2}} d x_1 d \tau \\
	& \leq C \delta_0 \int_0^t\Big\|\sqrt{|\px u_{1}^{s}|} \tilde{\Psi}_1\Big\|^2 d \tau+C \int_0^t \iint \frac{\nu(|v|)|\tilde{\mathbf{G}}|^2}{\mathbf{M}_*} d v d s d \tau,
\end{align}
where in the last inequality, we have used \cref{sw}. Similarly, for $\mathcal{K}_{11}^5$, one has
\begin{align}\label{k15}\notag
	\mathcal{K}_{11}^5 & =\int_0^t \int \frac{\bar{\Psi}_1}{\bar{\theta}} \int v_1^2 \LL_{\mathbf{M}}^{-1}\left|2 Q\left(\mathbf{G}^{s_1}, \mathbf{G}^{s_1}\right)\right| d v d x_1 d \tau \\
	\notag& \leq C \int_0^t \int|\tilde{\Psi}_1|\Big(\int \frac{\nu^{-1}(|v|)\left|Q\left(\mathbf{G}^{s_1}, \mathbf{G}^{s_3}\right)\right|^2}{\mathbf{M}_*} d v\Big)^{\frac{1}{2}} d x_1 d \tau \\
	& \leq C \int_0^l \int|\tilde{\Psi}_1|\Big(\int \frac{\nu(|v|)\left|\mathbf{G}^{s_1}\right|^2}{\mathbf{M}_*} d v\Big)^{\frac{1}{2}}\Big(\int \frac{\nu(|v|)\left|\mathbf{G}^{s_3}\right|^2}{\mathbf{M}_*} d v\Big)^{\frac{1}{2}} d x_1 d \tau \leq C \delta_0^2. 
\end{align}
For $\mathcal{K}_{11}^6$, we write  
$$
\begin{aligned}
	\mathcal{K}_{11}^6
	& =-\int_0^t \int \frac{\tilde{\Psi}_1}{\hat{\theta}} \int v_1^2\left(\LL_{\mathbf{M}}^{-1}-\LL_{\mathbf{M}^{s_1}}^{-1}\right)\left(\mathbf{G}_t^{s_1}\right) d v d x_1 d \tau \\
	&\quad +\int_0^t \int \frac{\tilde{\Psi}_1}{\bar{\theta}} \int v_1^2\left(\LL_{\mathbf{M}}^{-1}-\LL_{\mathbf{M}^{s_1}}^{-1}\right)\left[Q\left(\mathbf{G}^{s_1}, \mathbf{G}^{s_1}\right)\right] d v d x_1 d \tau :=\mathcal{K}_{11}^{61}+\mathcal{K}_{11}^{62} \text {. }
\end{aligned}
$$
For the first term, noticing $\LL_{\mathbf{M}}^{-1}\LL_{\mathbf{M}}=\LL_{\mathbf{M}^{s_1}}\LL_{\mathbf{M}^{s_1}}^{-1}=I$ is the identity operator, one has
\begin{align}\notag\label{k16a}
	\mathcal{K}_{11}^{61} & =-\int_0^t \int \frac{\bar{\Psi}_1}{\tilde{\theta}} \int v_1^2\left(\LL_{\mathbf{M}}^{-1}-\LL_{\mathbf{M}^{s_1}}^{-1}\right)\left(\mathbf{G}_t^{s_1}\right) d v d x d \tau \\
	&\notag =-\int_0^t \int \frac{\tilde{\Psi}_1}{\tilde{\theta}} \int v_1^2 \LL_{\mathbf{M}}^{-1}\left[2 Q\left(\mathbf{M}^{s_1}-\mathbf{M}, \LL_{\mathbf{M}^{s_1}}^{-1}\left(\mathbf{G}_t^{s_1}\right)\right)\right] d v d x d \tau \\
	&\notag \leq C \int_0^t \int|\tilde{\Psi}_1|\left(\int \frac{\nu(|v|)\left|\mathbf{M}^{s_1}-\mathbf{M}\right|^2}{\mathbf{M}_*} d v\right)^{\frac{1}{2}}\left(\int \frac{\nu(|v|)\left|\LL_{\mathbf{M}^{s_1}}^{-1}\left(\mathbf{G}_t^{s_1}\right)\right|^2}{\mathbf{M}_*} d v\right)^{\frac{1}{2}} d x d \tau \\
	&\notag \leq C \int_0^t \int|\tilde{\Psi}_1||\left(\rho-\rho^{s_1}+u-u^{s_1}, \theta-\theta^{s_1}\right)|\left(\int \frac{\nu^{-1}(|v|)\left|\mathbf{G}_t^{s_1}\right|^2}{\mathbf{M}_*} d v\right)^{\frac{1}{2}} d x d \tau \\
	& \leq C \delta_0 \int_0^t\left\|\sqrt{\left|\px u_{1}^{s}\right|} \tilde{\Psi}_1\right\|^2 d \tau+C \delta_0 \int_0^t\|(\phi, \psi, \omega)\|^2 d \tau+C \int_0^t \int q \tilde{\Psi}_1^2 d x_1 d \tau+C \delta_0 .
\end{align}
Similar estimates hold for $\mathcal{K}_{11}^{62}$, $\mathcal{K}_{11}^i(i=7,8,9)$ since they share similar structure. By collecting all the above estimates \cref{k11b,k12,k13,k14,k15,k16a}, we have
\begin{multline}\label{k11}
	\mathcal{K}_{11} \leq  C_\sigma\bigg[ \iint \frac{|\tilde{\mathbf{G}}|^2}{\mathbf{M}_*}(t, x, v) d v d x_1+ \int_0^t \iint \frac{\nu(|v|)}{\mathbf{M}_*}|\tilde{\mathbf{G}}|^2 d v d x_1 d \tau\bigg] \\+C\big(\mathcal{E}(0)^2+\delta_0^{\frac{1}{2}}\big)+
C(\delta_0+\chi+\sigma)\Big(\sup_{0\le s\le t}\|(\tilde{\Psi}_1,\tilde{W}_1)(s)\|_{L^2_x}^2+\int_0^t\mathcal{D}(\tau)d\tau\Big).
\end{multline}
The other terms in $\mathcal{K}_1$ given by \cref{k1} share similar estimates. This completes the proof of \cref{be3}. 

\end{proof}

\subsubsection{Extra dissipation terms}
Next, we prove the \cref{be4}.

\begin{proof}[Proof of \cref{be4}]
Here, we have to estimate the third right-hand term of \cref{Eq42}. For this, we diagonalize the perturbed system \cref{sys-Eu0}. Let $V=(\Phi, \tilde{\Psi}_{1}, \tilde{W})^{t} ,$ then
\begin{align}\label{eq-dia}
	\pt V+A_{1}\px V+A_{2} V=A_{3}\px^2 V+A_{4},
\end{align}
where
\begin{align*}
	A_{1}=\left(\begin{array}{ccc}
	 \tilde{u}_{1} &  \tilde{\rho} & 0 \\
	\frac{2 \tilde{\theta}}{ 3\tilde{\rho}} &  \tilde{u}_{1} & \frac{2}{3} \\
	0 & \frac{2}{3} \tilde{\theta} &  \tilde{u}_{1}
	\end{array}\right),\quad A_{2}=\left(\begin{array}{ccc}
	 \px\tilde{u}_{1} &  \px\tilde{\rho} & 0 \\
	-\frac{2}{3}\frac{ \tilde{\theta}\px{\rhot}}{ \tilde{\rho}^{2}} & -\frac{\px\tilde{u}_{1}}{3} & \frac{2 \px\tilde{\rho}}{ 3\bar{\rho}} \\
	0 & -\frac{2}{3} \px\tilde{\theta} & - \px\tilde{u}_{1}
	\end{array}\right),\quad	A_{3}=\left(\begin{array}{ccc}
	0 & 0 & 0 \\
	0 & \frac{4\mu(\Tt)}{ 3 \tilde{\rho}} & 0 \\
	0 & 0 & \frac{\kappa(\Tt)}{ \tilde{\rho}}
\end{array}\right),
\end{align*}
and
\begin{align*}
	A_{4}=\left(\begin{array}{c}
	0  \\
	-\int v_1^2\left(\mr{\Pi}-\tilde{\Pi}\right) d v+J_{1}+K_1+N_{1}-Q_{1} \\[2mm]
		-\frac{1}{2} \int v_1|v|^2\left(\mr{\Pi}- \tilde{\Pi}\right) d v+\tilde{u}_1\int v_1^2\left(\mr{\Pi}-\tilde{\Pi}\right) d v+J_{4}+K_4+N_{4}-\left(Q_{2}-\tilde{u}_{1} Q_{1}\right)
\end{array}\right),
\end{align*}
Three eigenvalues of the matrix $A_{1}$ are
\begin{align*}
	\tilde{\lambda}_{1}= \tilde{u}_{1}-\sqrt{\frac{10}{9} \tilde{\theta}}, \quad \tilde{\lambda}_{2}= \tilde{u}_{1}, \quad \tilde{\lambda}_{3}= \tilde{u}_{1}+\sqrt{\frac{10}{9} \tilde{\theta}}
\end{align*}
with corresponding left and right eigenvectors given by
\begin{align*}
	l_{1}=\left(\tilde{\theta},-\sqrt{\frac{5 \tilde{\theta}}{2}} ,  \tilde{\rho}\right), \quad l_{2}=\left(\tilde{\theta}, 0,-\frac{3}{2}\tilde{\rho}\right), \quad l_{3}=\left(\tilde{\theta}, \sqrt{\frac{5 \tilde{\theta}}{2}},  \tilde{\rho}\right)
\end{align*}
and
$$
r_{1}=\frac{3}{10\tilde{\rho} \tilde{\theta}}\left( \tilde{\rho},-\sqrt{\frac{10 \tilde{\theta}}{9}} ,  \frac{2}{3}\tilde{\theta}\right)^{t},\ \  r_{2}=\frac{2}{5  \tilde{\rho} \tilde{\theta}}( \tilde{\rho}, 0,-\tilde{\theta})^{t}, \ \ r_{3}=\frac{3}{10\tilde{\rho} \tilde{\theta}}\left( \tilde{\rho}, \sqrt{\frac{10 \tilde{\theta}}{9}}, \frac{2}{3}\tilde{\theta}\right)^{t},
$$
respectively. Denoting
$$
L=\left(l_{1}, l_{2}, l_{3}\right), \quad R=\left(r_{1}, r_{2}, r_{3}\right),\qquad Z:= L V,
$$
then we have
$$
L R=I d ., \quad L A_{1} R=\Lambda:=\operatorname{diag}\big(\tilde{\lambda}_{1}, \tilde{\lambda}_{2}, \tilde{\lambda}_{3}\big), \quad V=R Z, 
$$
where $Id.$ is the $3\times 3$ identity matrix. Similarly, we denote $R^{s_i}$ by replacing $(\tilde{\rho},\tilde{\theta})$ by $(\rho^{s_i},\theta^{s_i})$ with $i=1,3$. 
Multiplying the system \cref{eq-dia} by $L$ on the left, one obtains the diagonalized system for $Z$
\begin{align}\label{eq-dial}
	\pt Z+\Lambda\px Z-L A_{3} R\px^2 Z=-L\left(\pt R+A_{1}\px R-A_{3} \px^2R\right) Z-L A_{2} R Z+2 L A_{3} \px R \px Z+L A_{4}.
\end{align}
Moreover, we need some weight functions:
\begin{align}
	\alpha(t, x_1)=\frac{\rho^{s_{1}}(t, x_1)}{\rho_{\#}}, \quad \beta(t, x_1)=\frac{\rho^{s_{3}}(t, x_1)}{\rho_{\#}}.
\end{align}
From the positivity of the shock profile to the Boltzmann equation, we have
$$
\px\lambda_{i}^{s_{i}}<0 \quad \text { and }\quad  \px\rho^{s_{i}}<0 \quad(i=1,3).
$$
Thus it holds that
$$
\alpha, \beta<1 \text { and }|\alpha-1|,|\beta-1| \leq \frac{\delta}{\rho_{\#}} \ll 1 \text { if } \delta \ll 1.
$$
Taking the inner product of \eqref{eq-dial} with $$\bar{Z}:=\left(Z_{1}, \alpha^{N} Z_{2}, \alpha^{N} Z_{3}\right)^{t}$$ with $N=\delta^{-\frac{1}{2}}$, and then integrating the resulting equation over $x_1\in\R$, one has,
\begin{align*}
	&\frac{1}{2}{\left[\int_{\R}{Z_{1}^{2}+\alpha^{N}\left(Z_{2}^{2}+Z_{3}^{2}\right)}dx_1\right]_{t}-\int_{\R}\px\tilde{\lambda}_{1} \frac{Z_{1}^{2}}{2}+\alpha^{N} \sum_{i=2}^{3} \frac{\px\tilde{\lambda}_{i} Z_{i}^{2}}{2}} dx_1\\
	&\quad-\int_{\R}\Big(N \alpha^{N-1} \sum_{i=2}^{3} \frac{\big(\pt\alpha+\tilde{\lambda}_{i} \px \alpha \big) Z_{i}^{2}}{2}+\bar{Z} \cdot L A_{3} R \px^2 Z \Big)dx_1 +\int_{\R}\bar{Z} \cdot L A_{2} R Z dx_1\\
	&=-\int_{\R}\big(\bar{Z} \cdot L\left(\pt R+A_{1} \px R-A_{3} \px^2R \right) Z+\bar{Z} \cdot L A_{2} R Z \big)dx_1+\int_{\R}\big(2 \bar{Z} \cdot L A_{3} \px R   \px Z +\bar{Z} \cdot L A_{4}\big)dx_1.
\end{align*}
Note that for $i=2,3$,
$$
\pt\alpha+\tilde{\lambda}_{i} \px \alpha =-s_{1} \px \alpha +\tilde{\lambda}_{i} \px \alpha =\left(\lambda_{i}^{s_{1}}-s_{1}\right) \px \alpha +\left(\tilde{\lambda}_{i}-\lambda_{i}^{s_{1}}\right)\px \alpha 
$$
Then we further have
\begin{align}\label{Z111}\notag
&\frac{1}{2}{\left[\int_{\R}{Z_{1}^{2}+\alpha^{N}\left(Z_{2}^{2}+Z_{3}^{2}\right)}dx_1\right]_{t}-\int_{\R}\px\lambda_{1}^{s_{1}} \frac{Z_{1}^{2}}{2}+\alpha^{N} \sum_{i=2}^{3} \frac{\px\lambda_{i}^{s_{3}} Z_{i}^{2}}{2}} dx_1\\
&\notag\quad-\int_{\R}N \alpha^{N-1} \px \alpha  \sum_{i=2}^{3} \frac{\left(\lambda_{i}^{s_{1}}-s_{1}\right) Z_{i}^{2}}{2}dx_1-\int_{\R}\bar{Z}\cdot L A_{3} R \px^2 Z dx_1 \\
&\notag=-\int_{\R}\big(\bar{Z} \cdot L\left(R_{t}+A_{1} \px R  -A_{3} \px^2R \right) Z+\bar{Z} \cdot L A_{2} R Z\big)dx_1+\int_{\R}\big(2 \bar{Z} \cdot L A_{3} \px R   \px Z +\bar{Z} \cdot L A_{4}\big)dx_1 \\
&\quad+\int_{\R}\left(\px\tilde{\lambda}_{1}-\px\lambda_{1}^{s_{1}}\right) \frac{Z_{1}^{2}}{2}+\alpha^{N} \sum_{i=2}^{3} \frac{\big(\px\tilde{\lambda}_{i}-\px\lambda_{i}^{s_{3}}\big) Z_{i}^{2}}{2}+N \alpha^{N-1}\px \alpha  \sum_{i=2}^{3} \frac{\left(\tilde{\lambda}_{i}-\lambda_{i}^{s_{1}}\right) Z_{i}^{2}}{2}dx_1.
\end{align}
We will calculate the identity \eqref{Z111} terms by terms. For the dissipation terms, we have 
\begin{align*}
&-\bar{Z} \cdot L A_{3} R \px^2 Z =-\px\left(\bar{Z} \cdot L A_{3} R \px Z \right)+\px\bar{Z} \cdot L A_{3} R \px Z +\bar{Z} \cdot\px\left(L A_{3} R\right) \px Z  \\
&\quad=-\px\left(\bar{Z} \cdot L A_{3} R \px Z \right)+\px Z  \cdot L A_{3} R \px Z +\px(\bar{Z}-Z) \cdot L A_{3} R \px Z +\bar{Z} \cdot\px\left(L A_{3} R\right)\px Z. 
\end{align*}
Direct calculation yields that $LA_3 R$ is non-negative, thus,
$$
\px Z  \cdot L A_{3} R \px Z  \geq 0.
$$
The difference between $\bar{Z}$ and $Z$ is
\begin{align}
\begin{aligned}
\px(\bar{Z}-Z) &=\px\left(0,\left(\alpha^{N}-1\right) Z_{2},\left(\alpha^{N}-1\right) Z_{3}\right)^{t} \\
&=\left(\alpha^{N}-1\right)\left(0, \px Z_{2}, \px Z_{3}\right)^{t}+N \alpha^{N-1} \px\alpha\left(0, Z_{2}, Z_{3}\right)^{t}.
\end{aligned}
\end{align}
By the Lax entropy condition \cref{lax} to $1$-shock, we have
$$
\lambda_{2}^{s_{1}}-s_{1}>\lambda_{2}^{s_{1}}-\lambda_{1-}=u_{1}^{s_{1}}-\Big(u_{1-}-\frac{\sqrt{10 \theta_{-}}}{3}\Big) \geq \frac{\sqrt{10 \theta_{-}}}{3}-C \delta>\frac{\sqrt{10 \theta_{-}}}{6}
$$
and
$$
\lambda_{3}^{s_{1}}-s_{1}>\lambda_{2}^{s_{1}}-\lambda_{1-}>\frac{\sqrt{\frac{10}{9} \theta_{-}}}{2}.
$$
Therefore, by choosing $N=\frac{1}{\sqrt{\delta_{0}}}$ with $\delta_{0} \ll 1$, it holds that
\begin{align*}
&\px(\bar{Z}-Z) \cdot L A_{3} R \px Z \\
&\leq\left|\left(\alpha^{N}-1\right)\left(0, \px Z_{2}, \px Z_{3}\right)^{t} \cdot L A_{3} R \px Z \right|+\left|N \alpha^{N-1} \px\alpha\left(0, Z_{2}, Z_{3}\right)^{t} \cdot L A_{3} R \px Z \right| \\
&\leq C \delta\left|\px Z \right|^{2}+N \alpha^{N-1}\left|\px\alpha\right| \sum_{i=2}^{3}\left|Z_{i}\right|\left|\px Z \right| \\
&\leq \frac{N \alpha^{N-1}\left|\px\alpha\right|}{4} \sum_{i=2}^{3} \frac{\left(\lambda_{i}^{s_{1}}-s_{1}\right) Z_{i}^{2}}{2}+C \sqrt{\delta_{0}}\left|\px Z \right|^{2}.
\end{align*} 
Then one has
\begin{align}\label{Z111a}\notag
&\left|Z \cdot\px\left(L A_{3} R\right) \px Z \right| \leq C\left(\left|\px\lambda_{1}^{s_{1}}\right|+\left|\px\lambda_{1}^{s_{3}}\right|+\left|\px \Theta\right|+q\right)|Z|\left|\px Z \right| \\
&\notag\quad \leq C \sqrt{\delta_{0}}\left|\px Z \right|^{2}+\frac{C}{\sqrt{\delta_{0}}}\left(\left|\px\lambda_{1}^{s_{1}}\right|^{2}+\left|\px\lambda_{1}^{s_{3}}\right|^{2}+\left|\px \Theta\right|^{2}+q^{2}\right)|Z|^{2} \\
&\quad \leq C \sqrt{\delta_{0}}\left|\px Z \right|^{2}+C \sqrt{\delta_{0}}\left(\left|\px\lambda_{1}^{s_{1}}\right|+\left|\px\lambda_{1}^{s_{3}}\right|\right)|Z|^{2}+q|Z|^{2},
\end{align}
where $q$ is some function satisfying \eqref{DefQ}. 

\medskip 
Further, for terms involving derivatives of the ansatz, we have
\begin{align*}
L\left(R_{t}+A_{1} \px R  \right) &=L\left(R_{t}^{s_{1}}+R_{t}^{s_{3}}\right)+\Lambda L\left(\px R  ^{s_{1}}+\px R  ^{s_{3}}\right)+O(1)\left|\px \Theta\right|+q \\
&=\left(\Lambda-s_{1} I\right) L \px R  ^{s_{1}}+\left(\Lambda-s_{3} I\right) L \px R  ^{s_{3}}+O(1)\left|\px \Theta\right|+q.
\end{align*}
Using the smallness of $\tilde\lambda_1-s_1$ and $\tilde\lambda_3-s_3$ from Lax entropy condition \cref{lax}, we have 
\begin{multline}\label{Z111b}
-Z \cdot L\left(R_{t}+A_{1} \px R  \right) Z \leq C(\delta+\sigma)\left(\left|\px\lambda_{1}^{s_{1}}\right| Z_{1}^{2}+\left|\px \lambda_{3}^{s_{3}}\right| Z_{3}^{2}\right) \\\quad+C_{\sigma}\left[\left|\px\lambda_{1}^{s_{1}}\right|\left(Z_{2}^{2}+Z_{3}^{2}\right)+\left|\px \lambda_{3}^{s_{3}}\right|\left(Z_{1}^{2}+Z_{2}^{2}\right)\right]+C\left|\px \Theta\right||Z|^{2}+q|Z|^{2}.
\end{multline}
Similarly, noticing from  $\left|\px\lambda_{1}^{s_{1}}\right|+\left|\px\lambda_{1}^{s_{3}}\right|\le (\delta_0)^2$, we have 
\begin{align}\label{Z111c}\notag
2 \bar{Z} \cdot L A_{3} \px R   \px Z  & \leq C\left(\left|\px\lambda_{1}^{s_{1}}\right|+\left|\px\lambda_{1}^{s_{3}}\right|+\left|\px \Theta\right|\right)|Z|\left|\px Z \right| \\
& \leq C \sqrt{\delta_{0}}\left|\px Z \right|^{2}+C \sqrt{\delta_{0}}\left(\left|\px\lambda_{1}^{s_{1}}\right|+\left|\px\lambda_{1}^{s_{3}}\right|\right)|Z|^{2}+q|Z|^{2}.
\end{align}
The second-order derivative is
\begin{align}\label{Z111d}
Z \cdot L A_{3} \px^2R  Z \leq C \delta\left(\abs{\px\lambda^s}\right)|Z|^{2}+q|Z|^{2}.
\end{align}
For the quadratic terms of $Z$, we compute,
\begin{align*}
L A_{2} R=\left(\begin{array}{lll}
a_{11} & a_{12} & a_{13} \\
a_{21} & a_{22} & a_{23} \\
a_{31} & a_{32} & a_{33}
\end{array}\right),
\end{align*}
where $a_{i j}$ $(i, j=1,2,3)$ are linear functions of $ \px\tilde{\rho},  \px\tilde{u}_{1}, \px\tilde{\theta} .$ Then we further have
\begin{align}\label{Z111e}\notag
Z \cdot L A_{2} R Z=&\left(a_{11} Z_{1}^{2}+a_{12} Z_{1} Z_{2}+a_{13} Z_{1} Z_{3}\right)+\alpha^{N}\left(a_{21} Z_{1} Z_{2}+a_{22} Z_{2}^{2}+a_{23} Z_{2} Z_{3}\right) \\
&+\alpha^{N}\left(a_{31} Z_{1} Z_{3}+a_{32} Z_{2} Z_{3}+a_{33} Z_{3}^{2}\right):=a_{11} Z_{1}^{2}+\alpha^{N} a_{33} Z_{3}^{2}+Y,
\end{align}
where $Y$ satisfies 
\begin{align*}
&|Y|=\left|a_{12} Z_{1} Z_{2}+a_{13} Z_{1} Z_{3}+\alpha^{N}\left(a_{21} Z_{1} Z_{2}+a_{22} Z_{2}^{2}+a_{23} Z_{2} Z_{3}+a_{31} Z_{1} Z_{3}+a_{32} Z_{2} Z_{3}\right)\right| \\
&\leq C \sigma\left(\left|\px\lambda_{1}^{s_{1}}\right| Z_{1}^{2}+\left|\px \lambda_{3}^{s_{3}}\right| Z_{3}^{2}\right)+C_{\sigma}\left[\left|\px\lambda_{1}^{s_{1}}\right|\left(Z_{2}^{2}+Z_{3}^{2}\right)+\left|\px \lambda_{3}^{s_{3}}\right|\left(Z_{1}^{2}+Z_{2}^{2}\right)\right]+C\left(\left|\px \Theta\right|+q\right)|Z|^{2}.
\end{align*}
The estimates of $a_{11} Z_{1}^{2}$ and $\alpha^{N} a_{33} Z_{3}^{2}$ are  more subtle. By a direct calculation,
$$
\begin{aligned}
a_{11}= & \frac{1}{15 \tilde{\rho} \tilde{\theta}}\left(-\tilde{\rho} \tilde{\theta} \px\tilde{u}_{1}+\tilde{\rho} \sqrt{10 \tilde{\theta}} \px\tilde{\theta}-\tilde{\theta} \sqrt{10 \tilde{\theta}} \px\rho\right) \\
= & \frac{1}{15 \rho^{s_1 }\theta^{s_1}}\left(-\rho^{s_1} \theta^{s_1} \px u_{1}^{s_1}+\rho^{s_1} \sqrt{10 \theta^{s_1}} \px\theta^{s_1}-\theta^{s_1} \sqrt{10 \theta^{s_1}} \px\rho^{s_1}\right)-\frac{2 \sqrt{10 \theta_{\#}}}{15 \rho_{\#}} \px\Theta \\
& +\frac{1}{15 \rho^{s_3 }\theta^{s_3}}\left(-\rho^{s_3} \theta^{s_3} \px u_{1}^{s_3}+\rho^{s_3} \sqrt{10 \theta^{s_3}} \px\theta^{s_3}-\theta^{s_3} \sqrt{10 \theta^{s_3}} \px\rho^{s_3}\right)+q .
\end{aligned}
$$
By \cref{sw-fluid}, one has
\begin{align}\left\{
    \begin{aligned}
        &\rho^{s_i} \px u_{1}^{s_i}=\left(s_i-u_1^{s_i}\right) \px\rho^{s_i},\\
&-s_i \rho^{s_i} \px u_{1}^{s_i}+\rho^{s_i} u_1^{s_i} \px u_{1}^{s_i}+\frac{2}{3} \rho^{s_i} \px\theta^{s_i}+\frac{2}{3} \theta^{s_i} \px\rho^{s_i}=-\int v_1^2 \px\mathbf{G}^{s_i} d v,\\
&-\rho^{s_i} \px\theta^{s_i}=\frac{3}{2} \rho^{s_i}\left(u_1^{s_i}-s_i\right) \px u_{1}^{s_i}+\theta^{s_i} \px\rho^{s_i}+\frac{3}{2} \int v_1^2 \px\mathbf{G}^{s_i} d v.
    \end{aligned}\right.
\end{align}
Consequently, 
\begin{multline*}
    a_{11}= \sum_{i=1,3}\Big[\frac{\px\rho^{s_i}(\lambda_1^{s_i}-s_i)}{15 \rho^{s_i} \theta^{s_i}}\Big(\theta^{s_i}+\frac{3\sqrt{10 \theta^{s_i}}}{2} \left(\lambda_3^{s_i}-s_i\right)\Big) -\frac{\sqrt{10 \theta^{s_i}}}{10 \rho^{s_i}\theta^{s_i}} \int v_1^2 \px\mathbf{G}^{s_i} d v\Big]\\-\frac{2 \sqrt{10 \theta_{\#}}}{15 \rho_{\#}^{\#}} \px\Theta+q.
\end{multline*}
Therefore, it follows that
\begin{multline}\label{Z111f}
\int_{\R}a_{11} Z_1^2dx_1\\
\leq  C \delta_0\int_{\R}\left|\px \lambda_{1}^{s_1}\right| Z_1^2dx_1+C\int_{\R}\Big[\big(\left|\px \lambda_{3}^{s_3}\right|+\left|\px\Theta\right|+q\big) Z_1^2+\left(\left|\px \lambda_{1}^{s_1}\right|+\left|\px \lambda_{3}^{s_3}\right|\right)\left|Z_1\right|\left|\px Z_{1}\right|\Big] dx_1\\
\leq  C \int_{\R}\big(\delta_0\left|\px \lambda_{1}^{s_1}\right| Z_1^2+C\left|\px \lambda_{3}^{s_3}\right| Z_1^2\big)dx_1+C \delta_0\norm{\px Z}_{L^2_x}^2+C\int_{\R}\big(\left|\px\Theta\right| Z_1^2+q Z_1^2\big)dx_1,
\end{multline}
where we have used the facts that $\lambda_1^{s_1}-s_1=O(\delta)$ and $\px\rho^{s_i}\sim \px\theta^{s_i}\sim \int v_1^2 \mathbf{G}^{s_i} d v \sim \left|\px\lambda_{i}^{s_i}\right|=O(\delta^2)$, $(i=1,3)$. Similarly, one can derive that
\begin{multline}\label{Z111g}
\int_{\R}\alpha^N a_{33} Z_3^2 dx_1 \leq C\int_{\R}\left|\px \lambda_{1}^{s_1}\right| Z_3^2dx_1+C \delta_0\int_{\R}\left|\px \lambda_{3}^{s_3}\right| Z_3^2 dx_1+C \delta_0\norm{\px Z}_{L^2_x}^2\\+C\int_{\R}\big(\left|\px\Theta\right| Z_3^2+q Z_3^2\big) dx_1.
\end{multline}
Finally, substituting \cref{Z111f,Z111g} into \cref{Z111e}, we have,
\begin{align}\label{Z111h}\notag
\int_{\R}Z \cdot L A_{2} R Z\ dx_1 \leq & C\left(\sigma+\delta_{0}\right)\int_{\R}\left(\left|\px\lambda_{1}^{s_{1}}\right| Z_{1}^{2}+\left|\px \lambda_{3}^{s_{3}}\right| Z_{3}^{2}\right)dx_1\\
&\notag+C_{\sigma}\int_{\R}\left[\left|\px\lambda_{1}^{s_{1}}\right|\left(Z_{2}^{2}+Z_{3}^{2}\right)+\left|\px \lambda_{3}^{s_{3}}\right|\left(Z_{1}^{2}+Z_{2}^{2}\right)\right]dx_1 \\
&+C \delta_{0}\norm{\px Z}_{L_x^2}^{2}+C\int_{\R}\left|\px \Theta\right| Z^{2}+q|Z|^{2}dx_1.
\end{align}
Estimates on error term $A_4$ are similar to \cref{es-er0}-\cref{es-er1} and \cref{k11}, which is 
\begin{align}\label{za4}\notag
&\int_0^t\int_{\R}Z \cdot L A_{4}dx_1d\tau \leq \epsilon\left[\|Z(t, \cdot)\|^{2}+\int_{0}^{t}\left\| \px Z \right\|^{2} d \tau\right]+C \delta_{0} \int_{0}^{t}\left\|\sqrt{\abs{\px\lambda^s}} Z\right\|^{2} d \tau \\
&\notag\quad+\int_{0}^{t} \int q|Z|^{2} d x_1 d \tau+C_\sigma\bigg[ \iint \frac{|\tilde{\mathbf{G}}|^2}{\mathbf{M}_*}(t, x, v) d v d x_1+ \int_0^t \iint \frac{\nu(|v|)}{\mathbf{M}_*}|\tilde{\mathbf{G}}|^2 d v d x_1 d \tau\bigg]\\
&+C\left(\delta_{0}+\chi\right) \int_{0}^{t}\|( \phi,  \psi, {w})\|^{2} d \tau +C \chi  \int_{0}^{t}\left\|\partial_1( \phi,  \psi, {w})\right\|^{2} d \tau+C\big(\mathcal{E}(0)^{2}+\delta_{0}^{\frac{1}{2}}\big).
\end{align}
Lastly, we need to estimate the difference between $\px\tilde{\lambda}_{i}$ and $\px\lambda^{s_1}_{i}$, and we calculate one of these cases:
\begin{align}\label{Z111i}\notag
\left(\px\tilde{\lambda}_{1}-\px\lambda_{1}^{s_{1}}\right) \frac{Z_{1}^{2}}{2} &=\left[ \px\tilde{u}_{1}-\px u_{1}^{s_{1}}-\sqrt{\frac{10}{9}}\px\left(\sqrt{\tilde{\theta}}-\sqrt{\theta^{s_{1}}}\right)\right] \frac{Z_{1}^{2}}{2} \\
& \leq C\left(\left|\px \lambda_{3}^{s_{3}}\right|+\left|\px \Theta\right|+q\right) Z_{1}^{2}.
\end{align}
Now we have done the estimates, integrating  \eqref{Z111} over $[0, t]$, combining the estimates \cref{Z111a,Z111b,Z111c,Z111d,Z111h,za4,Z111i}, choosing $\delta_{0} \ll 1$, and using Gr\"onwall's inequality, we arrive at
\begin{align}\label{es-zb}\notag
&\|Z(t, \cdot)\|^{2}+\int_{0}^{t} \int\left[\left|\px\lambda_{1}^{s_{1}}\right| Z_{1}^{2}+\sum_{i=2}^{3}\left|\px\lambda_{i}^{s_{3}}\right| Z_{i}^{2}+N\left|\px\alpha\right| \sum_{i=2}^{3} Z_{i}^{2}\right] d x_1 d \tau \\
&\notag\leq C\left(\sigma+\sqrt{\delta_{0}}\right) \int_{0}^{t}\left\|\px Z \right\|_{L^2_x}^{2} d \tau +C \int_{0}^{t} \int\left|\px \lambda_{3}^{s_{3}}\right|\left(Z_{1}^{2}+Z_{2}^{2}\right) d x_1 d \tau+C\left(\delta_{0}+\chi\right) \int_{0}^{t}\|( \phi,  \psi, {w})\|_{L^2_x}^{2} d \tau\\
&\notag\qquad+C \chi  \int_{0}^{t}\left\|\px( \phi,  \psi, {w})\right\|_{L^2_x}^{2} d \tau +C \int_{0}^{t}\left\|\sqrt{\left|\px \Theta\right|} Z\right\|_{L^2_x}^{2} d \tau+C\big(\mathcal{E}(0)^{2}+\delta_{0}^{\frac{1}{2}}\big)\\
&\qquad +C_\sigma\bigg[ \iint \frac{|\tilde{\mathbf{G}}|^2}{\mathbf{M}_*}(t, x, v) d v d x_1+ \int_0^t \iint \frac{\nu(|v|)}{\mathbf{M}_*}|\tilde{\mathbf{G}}|^2 d v d x_1 d \tau\bigg].
\end{align}
Similarly, if we take the inner product of \eqref{eq-dial} with $\tilde{Z}:=\left(\beta^{-N} Z_{1}, \beta^{-N} Z_{2}, Z_{3}\right)^{t}$, then we have
\begin{align*}
&\frac{1}{2}{\left[\int_{\R}{\beta^{-N}\left(Z_{1}^{2}+Z_{2}^{2}\right)+Z_{3}^{2}}dx_1\right]_{t}-\int_{\R}\px \lambda_{3}^{s_{3}} \frac{Z_{3}^{2}}{2}+\beta^{-N} \sum_{i=1}^{2} \frac{\px\lambda_{i}^{s_{1}} Z_{i}^{2}}{2}}dx_1\\
&\notag\quad+\int_{\R}N \beta^{-N-1} \px \beta \sum_{i=1}^{2} \frac{\left(\lambda_{i}^{s_{3}}-s_{3}\right) Z_{i}^{2}}{2}-\tilde{Z} \cdot L A_{3} R \px^2 Z dx_1  \\
&\notag=-\int_{\R}\tilde{Z} \cdot L\left(R_{t}+A_{1} \px R  -A_{3} \px^2R \right) Z+\tilde{Z} \cdot L A_{2} R Zdx_1+\int_{\R}2 \tilde{Z} \cdot L A_{3} \px R   \px Z +\tilde{Z} \cdot L A_{4} dx_1\\
&+\int_{\R}\left(\px \tilde{\lambda}_{3}-\px \lambda_{3}^{s_{3}}\right) \frac{Z_{3}^{2}}{2}+\beta^{-N} \sum_{i=1}^{2} \frac{\left(\px\tilde{\lambda}_{i}-\px\lambda_{i}^{s_{1}}\right) Z_{i}^{2}}{2}-N \beta^{-N-1} \px \beta \sum_{i=1}^{2} \frac{\left(\tilde{\lambda}_{i}-\lambda_{i}^{s_{3}}\right) Z_{i}^{2}}{2}dx_1.
\end{align*}
Then similar to \cref{es-zb}, one can get
\begin{align}\label{es-zt}\notag
&\|Z(t, \cdot)\|^{2}+\int_{0}^{t} \int\left[\left|\px \lambda_{3}^{s_{3}}\right| Z_{3}^{2}+\sum_{i=1}^{2}\left|\px\lambda_{i}^{s_{1}}\right| Z_{i}^{2}+N\left|\px \beta\right| \sum_{i=2}^{3} Z_{i}^{2}\right] d x_1 d \tau \\
\leq&\notag C\left(\sigma+\sqrt{\delta_{0}}\right) \int^{t}\left\|\px Z \right\|^{2} d \tau+C \int_{0}^{t} \int\left|\px\lambda_{1}^{s_{1}}\right|\left(Z_{2}^{2}+Z_{3}^{2}\right) d x_1 d \tau \\
&\notag+C\left(\delta_{0}+\chi\right) \int_{0}^{t}\|( \phi,  \psi, {w})\|^{2} d \tau+C\chi\int_{0}^{t}\|( \phi,  \psi, {w})\|^{2} d \tau +C \int_{0}^{t} \| \sqrt{\left|\px \Theta\right|} Z \|^{2} d \tau\\
&+C\big(\mathcal{E}(0)^{2}+\delta_{0}^{\frac{1}{2}}\big)+C_\sigma\bigg[ \iint \frac{|\tilde{\mathbf{G}}|^2}{\mathbf{M}_*}(t, x, v) d v d x_1+ \int_0^t \iint \frac{\nu(|v|)}{\mathbf{M}_*}|\tilde{\mathbf{G}}|^2 d v d x_1 d \tau\bigg].
\end{align}
Combining \cref{es-zb} and \cref{es-zt} and choosing $\delta_{0}$ sufficiently small such that $N$ large enough, we have
\begin{align}\label{450}
\begin{aligned}
\|Z(t, &\cdot)\|^{2}+C\int_{0}^{t}\int_{\R}\abs{\px\lambda^s}|Z|^2dx_1 d \tau 
\leq\left(\sigma+\sqrt{\delta_{0}}+\chi\right) \int_{0}^{t}\mathcal{D}(\tau) d \tau+C \int_{0}^{t}\left\|\sqrt{\left|\px \Theta\right|} Z\right\|^{2} d \tau \\
&+C\big(\mathcal{E}(0)^{2}+\delta_{0}^{\frac{1}{2}}\big)+C_\sigma\bigg[ \iint \frac{|\tilde{\mathbf{G}}|^2}{\mathbf{M}_*}(t, x, v) d v d x_1+ \int_0^t \iint \frac{\nu(|v|)}{\mathbf{M}_*}|\tilde{\mathbf{G}}|^2 d v d x_1 d \tau\bigg].
\end{aligned}
\end{align}
Note that $\|Z(t)\|^{2}_{L^2_{x_1}}$ is quivalent to $\|(\Phi,\tilde{\Psi}_1,\tilde{W})\|_{L^2_{x_1}}$. This completes the proof of \cref{be4}. 
\end{proof}

\begin{proof}[Proof of \cref{be5}]
Note the decay rate of $\|\px\Theta\|_{L^\infty_{x_1}}$ is $(1+t)^{-1}$, which is critical. Thus, a sharp estimate of the heat kernel is essential. Firstly, continuing the proof of \cref{be4}, one has
\begin{align}\label{451}
\left|\px \Theta \| Z\right|^{2} \leq C\left|\alpha_{2}\right|(1+t)^{-1} e^{-\frac{\left(x-u_{1 \#} t\right)^{2}}{8 a(1+t)}}\left(Z_{1}^{2}+Z_{2}^{2}+Z_{3}^{2}\right) \leq C\left|\alpha_{2}\right|\left[h\left(Z_{1}^{2}+Z_{3}^{2}\right)+h^{2} Z_{2}^2\right],
\end{align}
where
\begin{align}
	 h=\frac{1}{\sqrt{16 \pi a(1+t)}} \exp\Big(-\frac{\left(x-u_{1 \#} t\right)^{2}}{16 a(1+t)}\Big),
\end{align}
satisfying
$$
h_{t}+u_{1 \#} \px h=a \px^2 h, \quad a=\frac{3\kappa}{5\rho_{\#}}.
$$
Denoting
 $\eta_{1}=\exp \left(\int_{-\infty}^{x} h(y, t) d y\right), $
 it holds that
\begin{align}\label{etabound1}
1 \leq \eta_{1} \leq e,
\end{align}
and
$$
\eta_{1 t}=\eta_{1} \int_{-\infty}^{x_1} h_{t}(y, t) d y=\eta_{1}\left(a \px h-u_{1 \#} h\right), \quad \px\eta_{1}=\eta_{1} h.
$$
Set the matrices $\left(c_{i j}\right)_{n \times n}$ and $\left(b_{i j}\right)_{n \times n}$ as
$$
L A_{3} R:=\left(c_{i j}\right)_{n \times n}, \quad L\left(R_{t}+A_{1} \px R  -A_{3} \px^2R +A_{2} R\right):=\left(b_{i j}\right)_{n \times n},
$$
then we have estimates for $|\alpha|=1$
\begin{align}\label{escb}
    |c_{ij}|\le C,\qquad
  |\p^{\alpha}c_{ij}|+  |b_{ij}|\le C\delta|\px\lambda^s|+|\p^{\alpha}\Theta|,
\end{align}
where $C>0$ is a universal constant.
Multiplying \cref{eq-dia} by $L$ from the left, we have the equation for $Z_{1}$,
\begin{align}\label{eq-z1}
Z_{1 t}+\tilde{\lambda}_{1} \px Z_{1}=\sum_{j=1}^{3} c_{1 j} \px^2Z_{j}-\sum_{j=1}^{3} b_{1 j} Z_{j}+\left(2 L A_{3} \px R   \px Z +L A_{4}\right)_{1},
\end{align}
where $(\cdot)_{i}\ (i=1,2,3)$ denotes the $i$ th component of the vector $(\cdot)$.
Multiplying \cref{eq-z1} by $\eta_{1} Z_{1}$ and integrating over $x_1\in\R$, one has 
\begin{align*}
&\left(\int_{\R}\eta_{1} \frac{Z_{1}^{2}}{2}dx_1\right)_{t}-\int_{\R}\left(\eta_{1 t}+\tilde{\lambda}_{1} \px\eta_{1}\right) \frac{Z_{1}^{2}}{2}+\px \tilde{\lambda}_{1} \eta_{1} \frac{Z_{1}^{2}}{2}\ dx_1 \\
&=\sum_{j=1}^{3} \int_{\R} c_{1 j}\px^2Z_{j} \eta_{1} Z_{1}dx_1+\int_{\R}\left[-\sum_{j=1}^{3} b_{1 j} Z_{j}+\left(2 L A_{3} \px R   \px Z +L A_{4}\right)_{1}\right] \eta_{1} Z_{1} dx_1 \\
&=-\int_{\R}\sum_{j=1}^{3} \px Z_{j}\px\left(c_{1 j} \eta_{1} Z_{1}\right)dx_1+\int_{\R}\left[-\sum_{j=1}^{3} b_{1 j} Z_{j}+\left(2 L A_{3} \px R   \px Z +L A_{4}\right)_{1}\right] \eta_{1} Z_{1}dx_1.
\end{align*}
Note that
$$
\begin{aligned}
\eta_{1 t}+\tilde{\lambda}_{1} \px\eta_{1} &=\left(\tilde{\lambda}_{1}-u_{1 \#}\right) \eta_{1} h+\eta_{1} a \px h \\
&=-\frac{\sqrt{10\tilde{\theta}}}{3} \eta_{1} h+\eta_{1} a \px h+\left( \tilde{u}_{1}-u_{1 \#}\right) \eta_{1} h,
\end{aligned}
$$
then one further has
\begin{multline}\label{455}
\left(\int_{\R}\eta_{1} \frac{Z_{1}^{2}}{2}dx_1\right)_{t}+\int_{\R}\frac{\sqrt{10\tilde{\theta}}}{6} \eta_{1} h Z_{1}^{2}dx_1=-\sum_{j=1}^{3} \int_{\R}\px Z_{j}\px\left(c_{1 j} \eta_{1} Z_{1}\right)dx_1-\sum_{j=1}^{3}\int_{\R} b_{1 j} \eta_{1} Z_{j} Z_{1}dx_1 \\
+\int_{\R}\left(2 L A_{3} \px R   \px Z +L A_{4}\right)_{1} \eta_{1} Z_{1}dx_1+\int_{\R}a \eta_{1} \px h \frac{Z_{1}^{2}}{2}+\px\tilde{\lambda}_{1} \eta_{1} \frac{Z_{1}^{2}}{2}+\left( \tilde{u}_{1}-u_{1 \#}\right) \eta_{1} h \frac{Z_{1}^{2}}{2}dx_1.
\end{multline}
Integrating the above equation over $[0, t]$, for each right-hand terms, we have from \eqref{escb} that
\begin{align}\label{442a}
&\int_{0}^{t} \int-\sum_{j=1}^{3} \px Z_{j}\px\left(c_{1 j} \eta_{1} Z_{1}\right) d x_1 d \tau\nonumber \\
=&\nonumber\int_{0}^{t} \int-\sum_{j=1}^{3}\left(\eta_{1} h Z_{1} \px Z_{j} c_{1 j}+\eta_{1} Z_{1} \px Z_{j} \px c_{1 j}+c_{1 j} \eta_{1} \px Z_{j} \px Z_{1}\right) d x_1 d \tau \\
\leq &\int_{0}^{t} \int\frac{\sqrt{10 \tilde{\theta}}}{12} \eta_{1} h Z_{1}^{2} d x_1 d \tau+C \delta_{0} \int_{0}^{t} \| \sqrt{\abs{\px\lambda^s}} Z_{1} \|^{2} d \tau 
+\int_{0}^{t} \int q Z_{1}^{2} d x_1 d \tau+C \int_{0}^{t}\left\|\px Z \right\|^{2} d \tau,
\end{align}
and
\begin{multline}\label{442b}
\int_{0}^{t} \int-\sum_{j=1}^{3} b_{1 j} \eta_{1} Z_{j} Z_{1} d x_1 d \tau \\
\leq C \int_{0}^{t}\left[\left\|\sqrt{\abs{\px\lambda^s}} Z\right\|^{2}+\left\|\px Z \right\|^{2}\right] d \tau +C \int_{0}^{t} \int\left|\px \Theta \| Z\right|^{2} d x_1 d \tau+\int_{0}^{t} \int q|Z|^{2} d x_1 d \tau.
\end{multline}
For the third right-hand term of \eqref{455}, similar to \cref{za4}, one has,
\begin{multline}\label{442c}
\int_{0}^{t} \int\left|\left(2 L A_{3} \px R   \px Z +L A_{4}\right)_{1}\right| \eta_{1} Z_{1} d x_1 d \tau 
\leq \sigma\left\|Z_{1}(t, \cdot)\right\|^{2}+C \delta_{0} \int_{0}^{t} \| \sqrt{\abs{\px\lambda^s} }Z \|^{2} d \tau\\ +C_\sigma\bigg[ \iint \frac{|\tilde{\mathbf{G}}|^2}{\mathbf{M}_*}(t, x, v) d v d x_1+ \int_0^t \iint \frac{\nu(|v|)}{\mathbf{M}_*}|\tilde{\mathbf{G}}|^2 d v d x_1 d \tau\bigg]\\
+C \int_{0}^{t}\left\|\px Z \right\|^{2} d \tau+C\left(\delta_{0}+\chi\right) \int_{0}^{t}\|( \phi,  \psi, {w})\|^{2} d \tau+C \chi  \int_{0}^{t}\left\|\px( \phi,  \psi, {w})\right\|^{2} d \tau \\
+ \int_{0}^{t}\int q|Z|^{2} d x_1 d \tau+C\big(\mathcal{E}(0)^{2}+\delta_{0}^{\frac{1}{2}}\big).
\end{multline}
Lastly, for the last term in \eqref{455}, we estimate
\begin{align}\label{459}
\notag&\int_{0}^{t} \int\left[a \px h+\px\tilde{\lambda}_{1}+\left( \tilde{u}_{1}-u_{1 \#}\right) h\right] \eta_{1} \frac{Z_{1}^{2}}{2} d x_1 d \tau \\
\notag&\quad\leq C\left(\delta+\left|\alpha_{2}\right|\right) \int_{0}^{t} \int \eta_{1} h Z_{1}^{2} d x_1 d \tau+C \int_{0}^{t} \int\left(\abs{\px\lambda^s}+\left|\px h\right|+q\right) Z_{1}^{2} d x_1 d \tau \\
&\quad\leq C \delta_{0} \int_{0}^{t} \int \eta_{1} h Z_{1}^{2} d x_1 d \tau+C \int_{0}^{t} \int\left(\abs{\px\lambda^s}+\left|\px h\right|+q\right) Z_{1}^{2} d x_1 d \tau
\end{align}
To deal with the last second term in \eqref{459}, notice that
\begin{align}\label{442e}\notag
&\int_{0}^{t} \int\left|\px h\right| Z_{1}^{2} d x_1 d \tau=\int_{0}^{t} \int \frac{1}{\sqrt{16 \pi a(1+\tau)}} \frac{2\left|x-u_{1 \#} \tau\right|}{16 a(1+\tau)} e^{-\frac{\left(x-u_{1 \#}\tau\right)^{2}}{16 a(1+\tau)}} Z_{1}^{2} d x_1 d \tau\\
\notag&\quad\leq \sigma \int_{0}^{t} \int \frac{1}{\sqrt{16 \pi a(1+\tau)}} e^{-\frac{\left(x-u_{1 \#} \tau\right)^{2}}{16 a(1+\tau)}} Z_{1}^{2} d x_1 d \tau+C_{\sigma} \int_{0}^{t} \int\frac{\left(x-u_{1 \#} \tau\right)^{2}}{16 a(1+\tau)^{1+3/2}} e^{-\frac{\left(x-u_{1\#} \tau\right)^{2}}{16 a(1+\tau)}} Z_{1}^{2} d x_1 d \tau\\
&\quad\leq \sigma \int_{0}^{t} \int h Z_{1}^{2} d x_1 d \tau+C_{\sigma} \int_{0}^{t}(1+\tau)^{-3 / 2}\left\|Z_{1}\right\|^{2}_{L^2_{x_1}} d \tau.
\end{align}
Integrating \eqref{455} over $[0,t]$, combining the above estimates \cref{etabound1,442a,442b,442c,459,442e}, choosing $\delta_{0}$ suitably small, and using Gr\"onwall's inequality, we finally arrive at
\begin{align}\label{46Z1}
&\notag\int Z_{1}^{2} d x_1+\int_{0}^{t} \int h Z_{1}^{2} d x_1 d \tau\\
 &\notag\quad\leq C\big(\mathcal{E}(0)^{2}+\delta_{0}^{\frac{1}{2}}\big)+C \int_{0}^{t}\left[\left\|\sqrt{\abs{\px\lambda^s}} Z\right\|^{2}+\left\|\px Z \right\|^{2}\right] d \tau\\
&\notag\qquad+C\left(\delta_{0}+\chi\right) \int_{0}^{t}\|( \phi,  \psi, {w})\|^{2} d \tau +C\chi\int_{0}^{t}\|( \phi,  \psi, {w})\|^{2} d \tau+C \int_{0}^{t} \int\left|\px \Theta\right|\left(Z_{2}^{2}+Z_{3}^{2}\right) d x_1 d \tau\\
&\qquad+C_\sigma\bigg[ \iint \frac{|\tilde{\mathbf{G}}|^2}{\mathbf{M}_*}(t, x, v) d v d x_1+ \int_0^t \iint \frac{\nu(|v|)}{\mathbf{M}_*}|\tilde{\mathbf{G}}|^2 d v d x_1 d \tau\bigg]. 
\end{align}
Note from \eqref{DefQ} and \eqref{priori} that $\int_{0}^{t}\int q|Z|^{2} d x_1 d \tau\le C\delta_0^{\frac{1}{2}}$. 
Similarly, for the term $Z_3$, we have 
\begin{align}\label{46Z3}\notag
&\int Z_{3}^{2} d x+\int_{0}^{t} \int h Z_{3}^{2} d x_1 d \tau \leq C\big(\mathcal{E}(0)^{2}+\delta_0^{\frac{1}{2}}\big)+C \int_{0}^{t}\left[\| \sqrt{\abs{\px\lambda^s}} Z\|^{2}+\| \px Z  \|^{2}\right] d \tau\\
&\notag+C\left(\delta_{0}+\chi\right) \int_{0}^{t}\|( \phi,  \psi, {w})\|^{2} d \tau 
+C\chi\int_{0}^{t}\|( \phi,  \psi, {w})\|^{2} d \tau+C \int_{0}^{t} \int\left|\px \Theta\right|\left(Z_{1}^{2}+Z_{2}^{2}\right) d x_1 d \tau\\
&+C_\sigma\bigg[ \iint \frac{|\tilde{\mathbf{G}}|^2}{\mathbf{M}_*}(t, x, v) d v d x_1+ \int_0^t \iint \frac{\nu(|v|)}{\mathbf{M}_*}|\tilde{\mathbf{G}}|^2 d v d x_1 d \tau\bigg].
\end{align}
The only term to estimate in \eqref{46Z1} and \eqref{46Z3} is 
\begin{align}\label{Z2222}
\int_{0}^{t} \int\left|\px \Theta\right| Z_{2}^{2} d x_1 d \tau.
\end{align}
Setting
$$
n(x, t)=\int_{-\infty}^{x} h(y, t) d y,
$$
we have the following equation with direct calculation,
$$
n_{t}=a \px h-u_{1 \#} h, \quad\|n\|_{L^{\infty}}=1.
$$
Here, we introduce a useful lemma to estimate terms involving the heat kernel,
\begin{Lem}[\cite{HLM}]\label{HLM}
For $0<\tau \leq \infty$, suppose that $Z(t, x)$ satisfies
$$
Z \in L^{\infty}\left(0, t ; L^{2}(\mathbb{R})\right), \quad \px Z  \in L^{2}\left(0, t ; L^{2}(\mathbb{R})\right), \quad Z_{t} \in L^{2}\left(0, t ; H^{-1}(\mathbb{R})\right)
$$
Then the following estimate holds for any $\tau \in(0, t]$,
\begin{multline}
\int_{0}^{t} \int h^{2} Z^{2} d x_1 d \tau \leq  \int  Z^{2}(x, 0) d x+4 a\int_{0}^{t}\left\|\px Z \right\|^{2} d \tau \\+\frac{2}{a}\left(\int_{0}^{t}\int Z_{\tau}Z n^{2}dxd \tau-u_{1 \#} \int_{0}^{t} \int Z^{2} n h d x_1 d \tau\right).
\end{multline}
\end{Lem}
To deal with the term \eqref{Z2222}, we use \cref{eq-dia} to write the equations for $Z_2$:
\begin{align}\label{eq-z2}
Z_{2 t}+\tilde{\lambda}_{2} \px Z_{2}=\sum_{j=1}^{3} c_{2 j} \px^2Z_{j}-\sum_{j=1}^{3} b_{2 j} Z_{j}+\left(2 L A_{3} \px R   \px Z +L A_{4}\right)_{2}.
\end{align}
Direct computation yields,
\begin{multline}\label{466}
\int_{0}^{t}\int Z_{2 \tau} Z_{2} n^{2}dx d \tau-u_{1 \#} \int_{0}^{t} \int Z_{2}^{2} n h d x_1 d \tau \\
=\int_{0}^{t} \int\Big\{\Big(\sum_{j=1}^{3} c_{2 j} \px^2Z_{j}-\sum_{j=1}^{3} b_{2 j} Z_{j}+\left(2 L A_{3} \px R   \px Z +L A_{4}\right)_{2}\Big) Z_{2} n^{2}\\
-\Big(\tilde{\lambda}_{2} \px Z_{2} Z_{2} n^{2}+u_{1 \#} Z_{2}^{2} n h\Big)\Big\} d x_1 d \tau.
\end{multline}
The estimates are similar to the above estimates of $Z_1$ and $Z_3$. Firstly, using \eqref{escb}, we have 
\begin{align*}
&\int_{0}^{t} \int\Big(\sum_{j=1}^{3} c_{2 j} \px^2Z_{j}-\sum_{j=1}^{3} b_{2 j} Z_{j}\Big) Z_{2} n^{2} d x_1 d \tau \\
&\quad=-\sum_{j=1}^{3} \int_{0}^{t} \int\left(2 h Z_{2} c_{2 j} n \px Z_{j}+c_{2 j x} n^{2} Z_{2} \px Z_{j}+c_{2 j} n^{2} \px Z_{j} \px Z_{2}+b_{2 j} Z_{j} Z_{2} n^{2}\right) d x_1 d \tau \\
&\quad\leq \frac{a}{4} \int_{0}^{t} \int h^{2} Z_{2}^{2} d x_1 d \tau+C \int_{0}^{t}\left[\| \sqrt{\abs{\px\lambda^s}} Z\left\|^{2}+\right\| \px Z  \|^{2}\right] d \tau\\
&\qquad+C \int_{0}^{t} \int\left|\px \Theta\right| Z^{2} d x_1 d \tau+\int_{0}^{t} \int q|Z|^{2} d x_1 d \tau.
\end{align*}
Also, a similar estimate as \cref{za4} yields
\begin{align*}
&\int_{0}^{t} \int\left|\left(2 L A_{3} \px R   \px Z +L A_{4}\right)_{2}\right| Z_{2} n^{2} d x_1 d \tau \\
&\quad\leq \epsilon\left\|Z_{2}(t, \cdot)\right\|^{2}+C \delta_{0} \int_{0}^{t} \| \sqrt{\abs{\px\lambda^s}} Z \|^{2} d \tau +C \int_{0}^{t}\left\|\px Z \right\|^{2} d \tau+C\left(\delta_{0}+\chi\right) \int_{0}^{t}\|( \phi,  \psi, {w})\|^{2} d \tau\\
& \qquad+C\chi\int_{0}^{t}\|( \phi,  \psi, {w})\|^{2} d \tau+\int_{0}^{t} \int q|Z|^{2} d x_1 d \tau+C\big(\mathcal{E}(0)^{2}+\delta_{0}^{\frac{1}{2}}\big)\\
&\qquad+C_\sigma\bigg[ \iint \frac{|\tilde{\mathbf{G}}|^2}{\mathbf{M}_*}(t, x, v) d v d x_1+ \int_0^t \iint \frac{\nu(|v|)}{\mathbf{M}_*}|\tilde{\mathbf{G}}|^2 d v d x_1 d \tau\bigg].
\end{align*}
Lastly, we find the cancellation for the last term in \eqref{466} as follows. 
\begin{align*}
&-\int_{0}^{t} \int\left(\tilde{\lambda}_{2} \px Z_{2} Z_{2} n^{2}+u_{1 \#} Z_{2}^{2} n h\right) d x_1 d \tau=\int_{0}^{t} \int\left[\px\tilde{\lambda}_{2} n^{2} \frac{Z_{2}^{2}}{2}+\left(\tilde{\lambda}_{2}-u_{1 \#}\right) n h Z_{2}^{2}\right] d x_1 d \tau \\
&=\int_{0}^{t} \int\left[\left(\px u_{1}^{s_{1}}+\px u_{1}^{s_{3}}+q\right) n^{2} \frac{Z_{2}^{2}}{2}+\left(u_{1}^{s_{1}}-u_{1 \#}+u_{1}^{s_{3}}-u_{1 \#}-\frac{a}{\rho_{\#}} \px \Theta+q\right) n h Z_{2}^{2}\right] d x_1 d \tau \\
&\leq C \int_{0}^{t} \| \sqrt{\abs{\px\lambda^s}} Z_{2} \|^{2} d \tau+\int_{0}^{t} \int q Z_{2}^{2} d x_1 d \tau .
\end{align*}
Therefore, combining the above estimates and using energy estimate in \cref{HLM}, we have
\begin{align}\label{pte}
\begin{aligned}
&\int_{0}^{t} \int h^{2} Z_{2}^{2} d x_1 d \tau \leq C(\epsilon+\delta)\left\|Z_{2}\right\|^{2}+C \int_{0}^{t}\| \sqrt{\abs{\px\lambda^s}} Z\|^{2}+\| \px Z  \|^{2} d \tau\\
&+C\left(\delta_{0}+\chi\right) \int_{0}^{t}\|( \phi,  \psi, {w})\|^{2} d \tau+C \chi  \int_{0}^{t}\left\|\partial_{1}( \phi,  \psi, {w})\right\|^{2} d \tau \\
&+C \int_{0}^{t} \int\left|\px \Theta\right|\left(Z_{1}^{2}+Z_{3}^{2}\right) d x_1 d \tau+\int_{0}^{t} \int q|Z|^{2} d x_1 d \tau+C\big(\mathcal{E}(0)^{2}+\delta_{0}^{\frac{1}{2}}\big) \\
&+C \int_{0}^{t}(1+\tau)^{-\frac{3}{2}}\left\|Z_{2}\right\|^{2} d \tau+C_\sigma\bigg[ \iint \frac{|\tilde{\mathbf{G}}|^2}{\mathbf{M}_*}(t, x, v) d v d x_1+ \int_0^t \iint \frac{\nu(|v|)}{\mathbf{M}_*}|\tilde{\mathbf{G}}|^2 d v d x_1 d \tau\bigg],
\end{aligned}
\end{align}
where we used \eqref{451} and the smallness of $\alpha_2$ to eliminate the term $|\px\Theta|Z_2^2$. 
Combining the estimates of $(Z_1,Z_2,Z_3)$ from \cref{46Z1,46Z3,pte}, using Gr\"onwall's inequality and choosing $\delta_{0}$ suitably small, we can obtain
\begin{align}\label{470}\notag
&\left\|\left(Z_{1}, Z_{3}\right)\right\|^{2}+\int_{0}^{t} \int\left[h\left(Z_{1}^{2}+Z_{3}^{2}\right)+h^{2} Z_{2}^{2}\right] d x_1 d \tau \\
&\notag\quad\leq C(\epsilon+\delta)\left\|Z_{2}\right\|^{2}+C \int_{0}^{t}\left[\| \sqrt{\abs{\px\lambda^s}} Z\|^{2}+\| \px Z  \|^{2}\right] d \tau +C\left(\delta_{0}+\chi\right) \int_{0}^{t}\|( \phi,  \psi, {w})\|_{H^1_x}^{2} d \tau\\
&\notag\qquad+C\big(\mathcal{E}(0)^{2}+\delta_{0}^{\frac{1}{2}}\big) +C \int_{0}^{t}(1+\tau)^{-3 / 2}\left\|Z_{2}\right\|^{2} d \tau\\
&\qquad+C_\sigma\bigg[ \iint \frac{|\tilde{\mathbf{G}}|^2}{\mathbf{M}_*}(t, x, v) d v d x_1+ \int_0^t \iint \frac{\nu(|v|)}{\mathbf{M}_*}|\tilde{\mathbf{G}}|^2 d v d x_1 d \tau\bigg], 
\end{align}
where we also used \eqref{451} and the smallness of $\alpha_2$ to eliminate the term $|\px\Theta|(Z_1^2+Z_3^2)$. 
Now we have completed the estimate for $\|Z\|^2$. Combining \eqref{451}, \eqref{450} and \eqref{470}, and choosing suitably small $\delta_0,\sigma,\chi,\alpha_2>0$, we deduce 
\begin{align}\label{471}\notag
&\|Z(t, \cdot)\|^{2}+\int_{0}^{t}\left\|\sqrt{\abs{\px\lambda^s}+\left|\px \Theta\right|} Z\right\|^{2} d \tau \\
&\notag\quad\leq C\big(\mathcal{E}(0)^{2}+\delta_{0}^{\frac{1}{2}}\big)+C(\delta^{\frac{1}{2}}_0+\chi+\sigma)\int_{0}^t\mathcal{D}(\tau)d\tau\\
&\qquad+C_\sigma\bigg[ \iint \frac{|\tilde{\mathbf{G}}|^2}{\mathbf{M}_*}(t, x, v) d v d x_1+\int_0^t \iint \frac{\nu(|v|)}{\mathbf{M}_*}|\tilde{\mathbf{G}}|^2 d v d x_1 d \tau\bigg], 
\end{align}
where $|\alpha_2|\le \chi$. 
Noting that $Z$ is a linear combination of $(\Phi, \tilde{\Psi}_1, \tilde{W})$ with non-dengenerate coefficients, combining estimates \eqref{471} with \cref{be2} and \cref{be3}, we obtain 
\begin{align}\label{Z}\notag
&\|(\Phi, \tilde{\Psi}_1, \tilde{W})(t, \cdot)\|^{2}+\int_{0}^{t}\left[\left\|\sqrt{\abs{\px u_{1}^s}+\left|\px \Theta\right|}(\Phi, \tilde{\Psi}_1, \tilde{W})\right\|^{2}+\left\|\left(\px\tilde{\Psi}_1, \px\tilde{W}\right)\right\|^{2}\right] d \tau \\
&\quad\notag\leq C\left(\sigma+\chi+\sqrt{\delta_{0}}\right) \int_{0}^{t}\big(\| \px\Phi\|^{2}+\mathcal{D}(\tau)\big)d \tau+C\big(\mathcal{E}(0)^{2}+\delta_{0}^{\frac{1}{2}}\big)\\
&\qquad+C_\sigma\bigg[ \iint \frac{|\tilde{\mathbf{G}}|^2}{\mathbf{M}_*}(t, x, v) d v d x_1+ \int_0^t \iint \frac{\nu(|v|)}{\mathbf{M}_*}|\tilde{\mathbf{G}}|^2 d v d x_1 d \tau\bigg].
\end{align}
Here we write $\| \px\Phi\|^{2}$ to emphasize that it's about to be estimated below. 

\smallskip 
	\noindent\textbf{Estimate on $\left\|\px\Phi(t, \cdot)\right\|^{2}$.}
We rewrite \cref{sys-Eu0}$_2$ as
\begin{multline}\label{phi}
 \frac{4}{3} \frac{\mu(\tilde{\theta})}{\tilde{\rho}} \px\Phi_{t}+\tilde{\rho} \tilde{\Psi}_{1t}+\tilde{\rho} \tilde{u}_1\px\tilde{\Psi}_{1}+\frac{2}{3} \tilde{\theta}\px\Phi=-\frac{4 \mu(\tilde{\theta})}{3 \tilde{\rho}}\left(2 \px\tilde{\rho} \px\tilde{\Psi}_{1}+\px^2\tilde{\rho} \tilde{\Psi}_1\right)-\frac{4 \mu(\tilde{\theta})}{3 \tilde{\rho}}\px^2\left(\tilde{u}_1 \Phi\right) \\
 +\frac{1}{3} \tilde{\rho} \px\tilde{u}_{1} \tilde{\Psi}_1-\frac{2}{3}\px \tilde{\rho} \tilde{W}-\frac{2}{3} \tilde{\rho}\px \tilde{W}+\frac{2 \tilde{\theta} \px\tilde{\rho}}{3 \tilde{\rho}} \Phi-\int v_1^2\left(\Pi-\Pi^{s_1}-\Pi^{s_3}\right) d v+J_1+N_1-Q_1.
\end{multline}
Multiplying \cref{phi} by $\px\Phi$ and then integrating the resulting equation over $\R$, one has
\begin{multline}\label{phiss}
 \frac{d}{dt}\int_{\R}\Big(\frac{2 \mu(\tilde{\theta})}{3 \tilde{\rho}} \abs{\px\Phi}^2+\tilde{\rho} \px\Phi \tilde{\Psi}_1\Big)dx_1+\int_{\R}\frac{2 \tilde{\theta}}{3} \abs{\px\Phi}^2dx_1=\int_{\R}\abs{\px\left(\tilde{\rho} \tilde{\Psi}_1\right)}^2+\px\left(\tilde{\rho} \tilde{\Psi}_1\right) \px\tilde{u}_{1} \Phi dx_1 \\
 +\int_{\R}\left[-\frac{4 \mu(\tilde{\theta})}{3 \tilde{\rho}}\left(2 \px\tilde{\rho} \px\tilde{\Psi}_{1}+\px^2\tilde{\rho} \tilde{\Psi}_1\right)-\frac{4 \mu(\tilde{\theta})}{3 \tilde{\rho}}\px^2\left(\tilde{u}_1 \Phi\right)-\frac{2}{3} \tilde{\rho} \px\tilde{u}_{1} \tilde{\Psi}_1-\frac{2}{3} \px\tilde{\rho} \tilde{W}-\frac{2}{3} \tilde{\rho} \px\tilde{W}\right] \px\Phi dx_1\\
 +\int_{\R}\left[\frac{2 \tilde{\theta} \px\tilde{\rho}}{3 \tilde{\rho}} \Phi-\int v_1^2\left(\Pi-\Pi^{s_1}-\Pi^{s_3}\right) d v+J_1+N_1-Q_1\right] \px\Phi+\left(\frac{2 \mu(\tilde{\theta})}{3 \tilde{\rho}}\right)_t \abs{\px\Phi}^2 dx_1,
\end{multline}
where we have used the fact that by \cref{sys-Eu0}$_1$, 
$$
\px\Phi\big(\tilde{\rho} \tilde{\Psi}_{1 t}+\tilde{\rho} \tilde{u}_1 \px\tilde{\Psi}_{1}\big)=\big(\tilde{\rho} \tilde{\Psi}_1 \px\Phi\big)_t-\px\big(\tilde{\rho} \tilde{\Psi}_1 \Phi_t\big)-\abs{\px\big(\tilde{\rho} \tilde{\Psi}_1\big)}^2+\tilde{\rho} \px\tilde{u}_{1} \tilde{\Psi}_1 \px\Phi-\px\big(\tilde{\rho} \tilde{\Psi}_1\big) \px\tilde{u}_{1} \Phi.
$$
Similar to \cref{k11}, one has
$$
\begin{aligned}
& \int_0^t \int\left|\int v_1^2\left(\Pi_1-\Pi^{s_1}-\Pi^{s_3}\right) d v\right|^2 d x d \tau \leq C \delta_0^2+C \delta_0 \int_0^t\|(\phi, \psi, \omega)\|^2 d \tau \\
& \quad+C\left(\delta_0+\chi\right) \int_0^t \iint \frac{\nu(|v|)}{\mathbf{M}_*}|\tilde{\mathbf{G}}|^2 d v d x d \tau+C \sum_{\left|\alpha^{\prime}\right|=1} \int_0^t \iint \frac{\nu(|v|)}{\mathbf{M}_*}\left|\partial^{\alpha^{\prime}} \tilde{\mathbf{G}}\right|^2 d v d x d \tau.
\end{aligned}
$$
In addition, by integrating \cref{phiss} over $[0,t]$, it is direct to obtain 
\begin{align}
\label{Phix}\notag
& \left\|\px\Phi(t, \cdot)\right\|^2+\int_0^t\left\|\px\Phi\right\|^2 d \tau \leq C\left\|\tilde{\Psi}_1(t, \cdot)\right\|^2+C \int_0^t\left\|\left(\px\tilde{\Psi}_1, \px\tilde{W}\right)\right\|^2 d \tau \\
&\notag \quad+C \delta_0 \int_0^t\left\|\sqrt{\left|\px u_{1}^{s}\right|+\left|\px\Theta\right|}\left(\Phi, \tilde{\Psi}_1, \tilde{W}\right)\right\|^2 d \tau+C \chi \int_0^t\left\|\left(\px\Phi, \px\psi_{1}\right)\right\|^2 d \tau \\
&\notag \quad+C \delta_0 \int_0^t\|(\psi, \omega)\|^2 d \tau+C\left(\delta_0+\chi\right) \int_0^t \iint \frac{\nu(|v|)}{\mathbf{M}_*}|\tilde{\mathbf{G}}|^2 d v d x d \tau+C\left(\mathcal{E}(0)^2+\delta_0^{\frac{1}{2}}\right) \\
& \quad+C \sum_{\left|\alpha^{\prime}\right|=1} \int_0^t \iint \frac{\nu(|v|)}{\mathbf{M}_*}\left|\partial^{\alpha^{\prime}} \tilde{\mathbf{G}}\right|^2 d v d x d \tau+\int_0^t \int q|(\Phi, \tilde{\Psi}_1, \tilde{W})|^2 d x d \tau.
\end{align}
Finally, directly by \cref{sys-Eu0}, one has
\begin{multline}
\label{phit}
\int_0^t\left\|\left(\Phi_\tau, \tilde{\Psi}_{1\tau}, \tilde{W}_\tau\right)\right\|_{L^2_x}^2 \mathrm{~d} \tau \leqq C \int_0^t\left\|\left(\px\Phi, \px\tilde{\Psi}_1, \px\tilde{W}\right)\right\|_{H^1_x}^2 \mathrm{~d} \tau \\+C \delta_0 \int_0^t\left\|\sqrt{\left|\px u_{1}^{s}\right|+\left|\px\Theta\right|}(\Phi, \tilde{\Psi}_1, \tilde{W})\right\|_{L^2_x}^2\mathrm{d} \tau 
 +C\left(\delta_0+\chi\right) \int_0^t \iint \frac{\nu(|v|)}{\mathbf{M}_*}|\tilde{\mathbf{G}}|^2 \mathrm{d} v \mathrm{d} x \mathrm{d} \tau\\ +C \sum_{\left|\alpha^{\prime}\right|=1} \int_0^t \iint \frac{\nu(|v|)}{\mathbf{M}_*}\left|\partial^{\alpha^{\prime}} \tilde{\mathbf{G}}\right|^2 \mathrm{d} v \mathrm{d} x \mathrm{d} \tau+C \delta_0^2,
\end{multline}
where the last term follows from the estimate of $q$. 
This completes the proof of \cref{be5}. 
\end{proof}

\subsection{Estimate on non-fluid part}\label{secnonfluid}
Next, we will estimate $\G$ and the higher-order derivatives.

\medskip 
For any $0\leq|\alpha|\leq 2$, applying $\p^{\alpha}$ to \cref{gt}, one has
\begin{align}\label{gtal}
\begin{aligned}
& \p^{\alpha}\tilde{\mathbf{G}}_{t}-\p^{\alpha}\left(\LL_{\mathbf{M}} \tilde{\mathbf{G}}\right)=\p^{\alpha}\left\{-\mathbf{P}_1\left(v\cdot\nabla_x \tilde{\mathbf{G}}\right)+Q(\tilde{\mathbf{G}}, \tilde{\mathbf{G}})+2 Q\left(\tilde{\mathbf{G}}, \mathbf{G}^{s_1}+\mathbf{G}^{s_3}\right)\right. \\
& \quad+2 Q\left(\mathbf{G}^{s_1}, \mathbf{G}^{s_3}\right)-\left[\mathbf{P}_1\left(v\cdot\mathbf{M}\right)-\mathbf{P}_1^{s_1}\left(v\cdot \mathbf{M}^{s_1}\right)-\mathbf{P}_1^{s_3}\left(v\cdot \mathbf{M}^{s_3}\right)\right]+R_1+R_3\bigg\}.
\end{aligned}
\end{align}
Multiplying \cref{gtal} by $\frac{\p^{\alpha}\tilde{\mathrm{G}}}{\mathrm{M}_*}$ and then integrating the resulting equation over ${\R}^3\times\mathbb{D}\times[0,t]$, one has
\begin{align}\label{nonfluid}\notag
        &\iint \frac{|\p^{\alpha}\tilde{\mathbf{G}}(t)|^2}{\mathbf{M}_*} d v d x+\int_0^t \iint \frac{\nu(|v|)}{\mathbf{M}_*}\left|\p^{\alpha}\tilde{\mathbf{G}}\right|^2 d v d x d \tau\\
        &\quad\notag \leq C \int_0^t \iint \sum_{|\beta|=|\alpha|+1}\frac{\nu(|v|)}{\mathbf{M}_*}\left|\p^{\beta}\tilde{\mathbf{G}}\right|^2 d v d x d \tau+C \int_0^t\sum_{|\beta|=|\alpha|+1}\left\|\left(\p^{\beta}\phi, \p^{\beta}\psi, \p^{\beta}\omega\right)\right\|^2 d \tau \\
& \qquad+C\left(\chi+\delta_0\right) \int_0^t \mathcal{D}(\tau) d \tau+C\big(\mathcal{E}(0)^2+\delta_0^{\frac{1}{2}}\big) .
\end{align}
where we have used \cref{fg}, \cref{linear} and wave interaction estimates as \cref{k11}. 


\medskip 
Finally, we obtain \cref{be1} by combining \cref{Z,Phix,phit,nonfluid}. 

\subsection{Entropy-entropy flux pair}\label{EEP}
In this part, we study the energy estimate in the original system. We first rewrite \cref{NSG} as
\begin{align}\left\{
\begin{aligned}
\rho_{t}&+u \cdot \nabla \rho+\rho \operatorname{div} u=0, \\[-1mm]
\rho u_{t}&+\rho u \cdot \nabla u+\frac{2 }{3 }\theta \nabla \rho+\frac{2}{3}\rho \nabla \theta=\mu(\theta)\big(\Delta u+\frac{1}{3} \nabla \operatorname{div} u\big) \\[-1mm]
&+\mu^{\prime}(\theta) \nabla \theta \cdot\big(\nabla u+(\nabla u)^{t}-\frac{2}{3}\mathbb{I}_1 \operatorname{div} u\big)-\int v \otimes v \cdot \nabla_{x} L_{\M}^{-1}\Pi d v, \\[-1mm]
\rho\theta_{t}&+\rho u \cdot \nabla \theta+\frac{2}{3} \rho\theta \operatorname{div} u=\kappa(\theta) \Delta \theta+{\mu(\theta)}\Big[\frac{\left(\nabla u+(\nabla u)^{t}\right)^{2}}{2}-\frac{2}{3}(\operatorname{div} u)^{2}\Big]  \\[-1mm]
&+{\kappa^{\prime}(\theta)}|\nabla \theta|^{2}-\int \frac{1}{2}|v|^{2} v \cdot \nabla_{x} L_{\M}^{-1}\Pi d v+u \cdot \int v \otimes v \cdot \nabla_{x} L_{\M}^{-1}\Pi d v .
\end{aligned}\right.
\end{align}
By this and \cref{tilderho}, recalling the definition of $(\phi,\varphi,\zeta)$ in \cref{perturbation,tipsi}, we have
\begin{align}\label{sys-perturb}
	\begin{cases}
		\phi_{t}+u \cdot \nabla \phi+\rho \operatorname{div} \varphi+\varphi \cdot \nabla \tilde{\rho}+\phi \operatorname{div} \tilde{u}=0, \\
		\rho\varphi_{t}+\rho u \cdot \nabla \varphi+\nabla(p-\tilde{p})+
		\rho\varphi \cdot \nabla \tilde{u}-\frac{\phi}{\tilde{\rho}}\nabla\tilde{p}={\mu(\thetat)}\left(\Delta \varphi+\frac{1}{3} \nabla \operatorname{div} \varphi\right)+\tilde{R},\\
		\rho\zeta_{t}+\rho u \cdot \nabla \zeta+\frac{2}{3} \rho\theta \operatorname{div} \varphi+\rho\varphi \cdot \nabla \tilde{\theta}+\frac{2}{3}\rho \zeta \operatorname{div} \tilde{u} =\kappa(\thetat)\Delta \zeta+\tilde{R}_4.
	\end{cases}
\end{align}
where 
\begin{align}\nonumber
	\tilde{R}&:=-\frac{\rho}{\rhot}\px Q_1\mathbb{I}_1+\left({\mu(\theta)}-\mu(\thetat)\right)\left(\Delta u+\frac{1}{3} \nabla \operatorname{div} u\right)+\mu(\thetat)\left(1-\frac{\rho}{\rhot}\right)\left(\Delta \tilde{u}+\frac{1}{3} \nabla \operatorname{div} \tilde{u}\right)\label{R}\\
	&\quad+\mu^{\prime}(\thetat) \nabla \thetat \cdot\left(\nabla \varphi+(\nabla \varphi)^{t}-\frac{2}{3}\mathbb{I} \operatorname{div} \varphi\right)-\mu'(\thetat)\frac{\nabla\thetat\phi}{\rhot}\left(\nabla \tilde{u}+(\nabla \tilde{u})^{t}-\frac{2}{3}\mathbb{I} \operatorname{div} \tilde{u}\right)\nonumber\\
	&\quad+\left({\mu^{\prime}(\theta)} \nabla \theta-\mu'(\thetat)\nabla\thetat\right) \cdot\left(\nabla {u}+(\nabla {u})^{t}-\frac{2}{3}\mathbb{I} \operatorname{div} {u}\right)\nonumber\\
	&\quad+ \int v \otimes v \cdot \nabla_{x} \left(L_{\M}^{-1}\Pi-\tilde{Y}\right) d v-\frac{\phi}{\rhot} \int v \otimes v \cdot \nabla_{x} \tilde{Y} d v,\\
	\tilde{R}_4&\nonumber:=-\frac{\rho}{\rhot}\px Q_2+\frac{\rho}{\rhot}\px Q_1\tilde{u}_1+\left(\kappa(\theta)-\kappa(\thetat)\right) \Delta \theta+{\kappa^{\prime}(\theta)}|\nabla \theta|^{2}-\kappa^{\prime}(\thetat)|\nabla \thetat|^{2}+\kappa^{\prime}(\thetat)\left(1-\frac{\rho}{\rhot}\right)|\nabla \thetat|^{2}\label{R_4} \\
	&\quad+\mu(\theta)\left[\frac{\left(\nabla u+(\nabla u)^{t}\right)^{2}}{2}-\frac{2}{3}(\operatorname{div} u)^{2}\right]-\mu(\thetat)\left[\frac{\left(\nabla \tilde{u}+(\nabla \tilde{u})^{t}\right)^{2}}{2}-\frac{2}{3}(\operatorname{div} \tilde{u})^{2}\right] \nonumber\\
	&\quad+\mu(\thetat)\left(1-\frac{\rho}{\rhot}\right)\left[\frac{\left(\nabla \tilde{u}+(\nabla \tilde{u})^{t}\right)^{2}}{2}-\frac{2}{3}(\operatorname{div} \tilde{u})^{2}\right]\nonumber\\
	&\quad+\left(-\int \frac{1}{2}|v|^{2} v \cdot (\nabla_{x} L_{\M}^{-1}\Pi-\nabla_{x} \tilde{Y} )d v+\tilde{u} \cdot \int v \otimes v \cdot(\nabla_{x} L_{\M}^{-1}\Pi-\nabla_{x} \tilde{Y})  d v\right)\nonumber\\
	&\quad+\varphi\cdot \int v \otimes v \cdot\nabla_{x} L_{\M}^{-1}\Pi  d v+\frac{\phi}{\rhot}\left(\int \frac{1}{2}|v|^{2} v \cdot \nabla_{x} \tilde{Y} d v-\tilde{u} \cdot \int v \otimes v \cdot \nabla_{x} \tilde{Y}d v\right),\\
\tilde{Y}:&= L_{{\M^{s_1}}}^{-1}{\Pi^{s_1}}+ L_{{\M^{s_3}}}^{-1}{\Pi^{s_3}}. 
\end{align}
We define the following functional as \cite{HLM},
\begin{equation}
	\begin{aligned}
	&\Gamma(s)=s-lns-1,\qquad\mathbb{E}=\frac{2}{3}\tilde{\T}\Gamma(\frac{ \tilde{\rho}}{\rho})+\frac{1}{2}|\varphi|^2+\tilde{\T}\Gamma(\frac{{\T}}{\tilde{\T}}).
	\end{aligned}
	\end{equation}
 Then direct calculations yield
 $$\rho\mathbb{E}=O(1)|\phi|^2+|\varphi|^2+|\zeta|^2. $$
	Calculating  $\cref{sys-perturb}_1\times (\frac{2}{3}\tilde{\T}(1-\frac{ \tilde{\rho}}{\rho}))+\cref{sys-perturb}_2\times\varphi+\cref{sys-perturb}_3\times\frac{\zeta}{\T}$, and then integrating the resulting equation over $\mathbb{D}\times[0,T]$, we get 
	\begin{equation}\label{*1}
	\begin{aligned}\sup_{0\leq \tau\leq t}\int_{\mathbb{D}}\rho\mathbb{E}dx+\int_0^t\int_\mathbb{D}^{}(|\nabla\varphi|^2+|\nabla\zeta|^2)dxd\tau
	\leq\int_0^t\int_\mathbb{D}^{} (\mathcal{H}-\mathcal{Q}) dxd\tau.
	\end{aligned}
	\end{equation}
	where $\mathcal{Q},\mathcal{H}$ equal to 	
	\begin{align*}
	\mathcal{Q}:=&\big[\Gamma(\frac{\tilde{\T}}{\T})-\frac{2}{3} \Gamma(\frac{ \tilde{\rho}}{\rho})\big]\rho u\cdot\nabla\tilde{\T}+\frac{2}{3}\tilde{\T}(\Gamma(\frac{ \tilde{\rho}}{\rho})+\Gamma(\frac{\rho}{ \tilde{\rho}})) u\cdot\nabla \tilde{\rho}\\
	&+\left[-\rho\varphi\cdot\nabla\tilde{u}+\frac{2}{3} \nabla\tilde{\rho}(\frac{\phi}{ \tilde{\rho}}\tilde{\T}-\zeta)-\frac{\phi}{ \tilde{\rho}}\dv \tilde{S}\right]\cdot\varphi\\
	&+\frac{2}{3}\tilde{\T}\frac{\phi}{\rho}(\varphi\cdot\nabla \tilde{\rho}+\phi \dv\tilde{u})+\zeta\nabla(\kappa(\Tt))\cdot\nabla\zeta +\varphi\cdot\nabla(\mu(\Tt))\cdot\nabla\varphi+\frac{1}{3}\nabla\mu(\Tt)\cdot\nabla\varphi \dv \varphi\\
	&+\left\{\big[\Gamma(\frac{\tilde{\T}}{\T})-\frac{2}{3} \Gamma(\frac{ \tilde{\rho}}{\rho})\big]\rho\tilde{\T}_{t}+\frac{2}{3}\tilde{\T}(\Gamma(\frac{ \tilde{\rho}}{\rho})+\Gamma(\frac{\rho}{ \tilde{\rho}})) \tilde{\rho}_{t}\right\},\\
	\mathcal{H}:=&\tilde{R}\cdot\varphi_1+\tilde{R}_4\frac{\zeta}{\T}.
	\end{align*}
	Direct estimate on $\mathcal{Q}$ shows that 
	\begin{align}
		|\mathcal{Q}|\leq O(1)(\delta+\chi)\left[(\phi^2+|\varphi|^2+\zeta^2)+|\nabla\phi|^2+|\nabla\varphi|^2+|\nabla\zeta|^2)\right].
	\end{align}
	To deal with the zero-order terms, we compute 
	\begin{align*}
&\int_0^t\int_\mathbb{D}^{}\phi^2+|\varphi|^2+\zeta^2dxdx'd\tau\\
		&\quad\leq\int_0^t\int_\mathbb{D}^{}| \phim|^2+| \mr{\varphi}|^2+|\zetam|^2+|\acute{\rho}|^2+|\acute{u}|^2+|\acute{\T}|^2dxd\tau\\		&\quad\leq\mathcal{E}(0)+\delta_0^{\frac{1}{2}}+\int_0^t\|\nabla(\phi,\varphi,\zeta)\|^2d\tau,
	\end{align*}
	where we have used \cref{prop-decom}. To estimate $\mathcal{H}$, one should notice the following fact
 \begin{align}\label{alpha}
     \abs{\p^{\alpha}(\rho,u,\theta)}\leq \abs{\p^{\alpha}(\phi,\varphi,\zeta)}+\abs{\p^{\alpha}(\rhot,\ut,\Tt)}\leq C(\delta_0+\chi),\qquad 1\leq|\alpha|\leq 2.
 \end{align}
 Then one has the following estimates for $\tilde{R}$
 \begin{align}\label{riaaa}
 \begin{aligned}
          |\tilde{R}_i|&\leq \abs{q}+\abs{\nabla_x(\phi,\varphi,\zeta)}^2+\abs{\nabla^2_x(\phi,\varphi,\zeta)}\abs{(\phi,\varphi,\zeta)}+(\delta_0+\chi)\abs{\phi,\varphi,\zeta}\\
          &\qquad\qquad+\Big|\int v \otimes v \cdot \nabla_{x} \left(L_{\M}^{-1}\Pi-\tilde{Y}\right) d v\Big|,\qquad  i=1,2,3,\\
          |\tilde{R}_4|&\leq \abs{q}+\abs{\nabla_x(\phi,\varphi,\zeta)}^2+\abs{\nabla^2_x(\phi,\varphi,\zeta)}\abs{(\phi,\varphi,\zeta)}+(\delta_0+\chi)\abs{\phi,\varphi,\zeta}\\
          &+\Big|-\int \frac{1}{2}|v|^{2} v \cdot (\nabla_{x} L_{\M}^{-1}\Pi-\nabla_{x} \tilde{Y} )d v+\tilde{u} \cdot \int v \otimes v \cdot(\nabla_{x} L_{\M}^{-1}\Pi-\nabla_{x} \tilde{Y})  d v\Big|.
      \end{aligned}
 \end{align}
 Then we can use a similar argument as $\mathcal{Q}$ and the wave interactions estimate similar to \cref{be3} to control $\mathcal{H}$. Then we arrive at
 \begin{align}\label{psi0}\notag
&\sup_{\tau\in[0,t]}\int_\mathbb{D}^{}(|\phi|^2,|\varphi|^2,|\zeta|^2)(\tau,x)dx+ \int_0^t\int_\mathbb{D}^{}(|\nabla\varphi|^2+|\nabla\zeta|^2)dxd\tau\\
	&\notag\quad\leq C(\delta_0+\chi+\sigma)\int_{0}^t \mathcal{D}(\tau) d\tau+\big(\mathcal{E}(0)^2+\delta^{\frac{1}{2}}\big)+C\left(\delta_0+\chi\right) \int_0^t \iint \frac{\nu(|v|)}{\mathbf{M}_*}|\tilde{\mathbf{G}}|^2 \mathrm{~d} v \mathrm{~d} x \mathrm{~d} \tau \\
 &\qquad+C \sum_{\left|\alpha^{\prime}\right|=1} \int_0^t \iint \frac{\nu(|v|)}{\mathbf{M}_*}\left|\partial^{\alpha^{\prime}} \tilde{\mathbf{G}}\right|^2 \mathrm{~d} v \mathrm{~d} x \mathrm{~d} \tau.
 \end{align}


\subsection{Higher-order estimate}
In this section, we estimate the higher order derivatives of $(\phi,\varphi,\zeta)$. 
\subsubsection{Dissipation estimate on $\p^{\alpha}(\phi,\varphi,\zeta)$ with $\alpha_0=0$ and $1\le|\alpha|\le 2$.}

We first introduce the systems for higher-order derivatives,
    Applying $\p^{\alpha}$ to \cref{sys-perturb} with $\alpha=(0,\alpha_1,\alpha_2,\alpha_3)$ and $1\leq|\alpha|\leq 2$, one has
    \begin{align}
        \label{sys-perturb-alpha}
	\begin{cases}
		\p^{\alpha}\phi_{t}+u \cdot \nabla\p^{\alpha} \phi+\rho \operatorname{div} \p^{\alpha}\varphi=J^1_{\alpha}, \\
		\rho\p^{\alpha}\varphi_{t}+\rho u \cdot \nabla\p^{\alpha} \varphi+\frac{2}{3}\rho\nabla\p^{\alpha}\zeta+\frac{2}{3}\T\nabla\p^{\alpha}\phi={\mu(\thetat)}\left(\Delta\p^{\alpha}\varphi+\frac{1}{3} \nabla \operatorname{div} \p^{\alpha}\varphi\right)+J_{\alpha}^2,\\
		\rho\p^{\alpha}\zeta_{t}+\rho u \cdot \nabla \p^{\alpha}\zeta+\frac{2}{3} \rho\theta \operatorname{div} \p^{\alpha}\varphi=\kappa(\thetat)\Delta \p^{\alpha}\zeta+J_{\alpha}^{3},
	\end{cases}
    \end{align}
where
\begin{align}
\begin{aligned}
    J^{1}_{\alpha}:=&\p^{\alpha}\big(u \cdot \nabla\p^{\alpha} \phi+\rho \operatorname{div} \p^{\alpha}\varphi\big)-u \cdot \nabla\p^{\alpha} \phi-\rho \operatorname{div} \p^{\alpha}\varphi+\p^{\alpha}\big(\varphi \cdot \nabla \tilde{\rho}+\phi \operatorname{div} \tilde{u}\big),\\
    J^{2}_{\alpha}:=&-\p^{\alpha}\Bigg[\rho u \cdot \nabla \varphi+\frac{2}{3}\rho\nabla\zeta+\frac{2}{3}\T\nabla\phi+{\mu(\thetat)}\left(\Delta\varphi+\frac{1}{3} \nabla \operatorname{div} \varphi\right) \Bigg]\\
    &+\rho u \cdot \nabla\p^{\alpha} \varphi+\frac{2}{3}\rho\nabla\p^{\alpha}\zeta+\frac{2}{3}\T\nabla\p^{\alpha}\phi+{\mu(\thetat)}\left(\Delta\p^{\alpha}\varphi+\frac{1}{3} \nabla \operatorname{div} \p^{\alpha}\varphi\right)\\
	&-\p^{\alpha}(	\rho\varphi \cdot \nabla \tilde{u}-\frac{\phi}{\tilde{\rho}}\nabla\tilde{p}-\tilde{R})-\p^{\alpha}\rho\varphi_t-\frac{2}{3}\p^{\alpha}\big(\nabla\rhot\zeta+\phi\nabla\tilde{\T}\big),\\
 J^{3}_{\alpha}:=&-\p^{\alpha}\bigg(\rho\zeta_{t}+\rho u \cdot \nabla \zeta+\frac{2}{3} \rho\theta \operatorname{div} \varphi\bigg)+\rho\p^{\alpha}\zeta_{t}+\rho u \cdot \nabla \p^{\alpha}\zeta+\frac{2}{3} \rho\theta \operatorname{div} \p^{\alpha}\varphi\\
 &-\p^{\alpha}\big(\rho\varphi \cdot \nabla \tilde{\theta}+\frac{2}{3}\rho \zeta \operatorname{div} \tilde{u}+ \kappa(\thetat)\Delta \zeta\big)-\kappa(\thetat)\Delta \p^{\alpha}\zeta+\p^{\alpha}\tilde{R}_4-\p^{\alpha}\rho\zeta_t.
 \end{aligned}
\end{align}
Multiplying \cref{sys-perturb-alpha}$_1$ by $\frac{2}{3}\frac{\theta}{\rho}\p^{\alpha}\phi$, multiplying \cref{sys-perturb-alpha}$_2$ by $\p^{\alpha}\varphi$, multiplying \cref{sys-perturb-alpha}$_3$ by $\frac{\p^{\alpha}\zeta}{\theta}$, adding them, and then integrating the resulting equation over $\mathbb{D}^{}$, one has
\begin{align}\label{alpha=0}
\begin{aligned}
    &\frac{d}{dt}\int_{\mathbb{D}^{}}\frac{\T}{3\rho}\abs{\p^{\alpha}\phi}^2+\frac{\rho}{2}\abs{\p^{\alpha}\varphi}^2+\frac{\rho}{2\T}\abs{\p^{\alpha}\zeta}^2dx+\int_{\mathbb{D}^{}}\abs{\nabla\p^{\alpha}\varphi}^2+\abs{\nabla\p^{\alpha}\zeta}^2dx\\
    \leq&\int_{\mathbb{D}^{}}\big(\abs{\pt\rho}+\abs{\pt\T}+\abs{\nabla_x\rho}+\abs{\nabla_x u}\big)\big(\abs{\p^{\alpha}\phi}^2+\abs{\p^{\alpha}\varphi}^2+\abs{\p^{\alpha}\zeta}^2\big)dx\\
    &+\int_{\mathbb{D}^{}} \big(\abs{\nabla\mu(\Tt)}+\abs{\nabla\kappa(\Tt)}\big)\big(\abs{\nabla\p^{\alpha}\varphi}\abs{\p^{\alpha}\varphi}+\abs{\nabla\p^{\alpha}\zeta}\abs{\p^{\alpha}\zeta}  \big)dx\\
    &+C\int_{\mathbb{D}^{}}\abs{J^1_{\alpha}}\abs{\p^{\alpha}\phi}+\abs{J^2_{\alpha}}\abs{\p^{\alpha}\varphi}+\abs{J^3_{\alpha}}\abs{\p^{\alpha}\zeta}dx.
    \end{aligned}
\end{align}
We only need to estimate the highest order terms in the third line on the right-hand side of \cref{alpha=0}. 
\begin{align}\label{111a}\notag
    \int_{\mathbb{D}^{}}\bigg(\kappa(\T)-\kappa(\Tt)\bigg)\Delta\p^{\alpha}\T\frac{\p^{\alpha}\zeta}{\theta}dx&\leq C\int_{\mathbb{D}^{}}\abs{\nabla\phi}\abs{\nabla^{|\alpha|+1}\zeta}\abs{\nabla^{|\alpha|}\zeta}+(\delta_0+\chi)\abs{\nabla^{|\alpha|+1}\zeta}^2dx\\
    &\leq C(\delta_0+\chi)\big(\norm{\nabla^{|\alpha|+1}\zeta}_{L_x^2}^2+\norm{\nabla^{|\alpha|}\zeta}_{L_x^2}^2\big).
\end{align}
and
\begin{align}\label{222}\notag
&\sum_{\substack{|\beta|+|\beta'|\leq|\alpha|+2\\\max\{|\beta|,|\beta'|\leq |\alpha|\}}}\int_{\mathbb{D}^{}}\abs{\p^{\beta}(\rho,u,\theta)}\abs{\p^{\beta'}(\rho,u,\theta)}\abs{\p^{\alpha}(\rho,u,\theta)}dx\\
&\quad\notag\leq \sum_{\substack{|\beta|+|\beta'|\leq|\alpha|+2\\\max\{|\beta|,|\beta'|\leq |\alpha|\}}}\big\|\nabla^{|\beta|}(\rho,u,\theta)\big\|_{L^4_x}\big\|\nabla^{|\beta'|}(\rho,u,\theta)\big\|_{L^4_x}\big\|\p^{\alpha}(\rho,u,\theta)\big\|_{L_x^2}\\
&\quad\notag\leq \sum_{\substack{|\beta|+|\beta'|\leq|\alpha|+2\\\max\{|\beta|,|\beta'|\leq |\alpha|\}}}\big\|\p^{\beta}(\rho,u,\theta)\big\|_{H^1_x}\big\|\p^{\beta'}(\rho,u,\theta)\big\|_{H^1_x}\big\|\p^{\alpha}(\rho,u,\theta)\big\|_{L_x^2}\\
&\quad\leq C(\delta_0+\chi)\big(\big\|\nabla^{|\alpha|+1}(\phi,\varphi,\zeta)\big\|_{L_x^2}^2+\big\|\nabla^{|\alpha|}(\phi,\varphi,\zeta)\big\|_{L_x^2}^2\big).
\end{align}
Then by \cref{111a,222} and wave interaction estimate similar to \cref{be3}, we have 
\begin{multline}\label{phial}
\norm{\big(\nabla^{|\alpha|}\phi,\nabla^{|\alpha|}\varphi,\nabla^{|\alpha|}\zeta\big)}_{L_x^2}^2+\int_{0}^t\norm{\big(\nabla^{|\alpha|+1}\varphi,\nabla^{|\alpha|+1}\zeta\big)}_{L_x^2}^2d\tau\\
\leq C\sum_{|\beta|=|\alpha|+1}\int_{0}^t\Big(\norm{(\p^{\beta}\phi,\p^{\beta}\varphi,\p^{\beta}\zeta)}_{L_x^2}^2+\iint\frac{\nu(|v|)}{\M_*}\abs{\p^{\beta}\tilde{\G}}^2dvdx\Big)d\tau+C(\delta+\chi+\sigma)\int_{0}^t\mathcal{D}(\tau)d\tau.
\end{multline}

\subsubsection{Dissipation estimate on $\p^{\alpha}(\phi,\psi,\zeta)$ with $\alpha_0\neq 0$ and $1\le|\alpha|\le 3$.}

    For convenience, we estimate $\p^{\alpha}\pt(\phi,\psi,\zeta)$ with $|\alpha|=0,1,2$. We introduce the following system,
    \begin{align}\label{sys-nonuni}
        \left\{\begin{aligned}
&\p^{\alpha}\phi_t+  \dv \p^{\alpha}\psi=0, \\
&\p^{\alpha}\psi_{1 t}  +\dv\p^{\alpha}\left(\rho u_1\uf-\tilde{\rho} \tilde{u}_1\tilde{\uf}\right)+\px\p^{\alpha}(p-\tilde{p})= -\int v_1v\cdot \nabla\p^{\alpha}\tilde{\mathbf{G}} d v,\\
& \qquad\qquad-\frac{4}{3}\px\p^{\alpha}\left(\mu(\tilde{\theta}) \px\tilde{u}_{1}-\mu\left(\theta^{s_1}\right) \px u_{1}^{s_1}-\mu\left(\theta^{s_3}\right)\px u_{1}^{s_3}\right)-\px\p^{\alpha} Q_{1}, \\
&\p^{\alpha}\psi_{i t}  +\dv\p^{\alpha}\left(\rho  u_i\uf\right)+\p_i\p^{\alpha}(p-\tilde{p})=-\int v_i v\cdot\nabla \p^{\alpha}\tilde{\mathbf{G}} d v, \quad i=2,3, \\
&\p^{\alpha}\omega_t+  \dv\p^{\alpha}\left(\rho \uf \theta-\tilde{\rho} \tilde{\uf} \tilde{\theta}+\rho \uf \frac{|u|^2}{2}-\tilde{\rho} \tilde{\uf} \frac{|\tilde{u}|^2}{2}+p \uf-\tilde{p} \tilde{\uf}\right) \\
 &\qquad\qquad= -\px\p^{\alpha}\left(k(\tilde{\theta}) \px\tilde{\theta}-k\left(\theta^{s_1}\right) \px\theta^{s_1}-k\left(\theta^{s_3}\right) \px\theta^{s_3}\right) \\
&\qquad\qquad\quad\ -\frac{4}{3}\px\p^{\alpha}\left(\mu(\tilde{\theta}) \tilde{u}_1 \px\tilde{u}_{1}-\mu\left(\theta^{s_1}\right) u_1^{s_1} \px u_{1}^{s_1}-\mu\left(\theta^{s_3}\right) u_1^{s_3}\px u_{1 }^{s_3}\right) \\
&\qquad\qquad\quad\ -\frac{1}{2} \int |v|^2v\cdot \nabla\p^{\alpha}\tilde{\mathbf{G}} d v-\px\p^{\alpha} Q_{2} .
\end{aligned}\right.
    \end{align}
    Multiplying \cref{sys-nonuni}$_1$ by $\p^{\alpha}\phi_t$, \cref{sys-nonuni}$_2$ by $\p^{\alpha}\psi_t$, \cref{sys-nonuni}$_3$ by $\p^{\alpha}\psi_{it}$, \cref{sys-nonuni}$_4$ by $\p^{\alpha}\zeta_t$, adding them up and then integrating the resulting equation over $\mathbb{D}^{}\times[0,t]$, one has,
    \begin{multline}\label{phitau}
        \int_0^t\int_{\mathbb{D}^{}} (\abs{\p^{\alpha}\phi_\tau}^2+\abs{\p^{\alpha}\psi_\tau}^2+\abs{\p^{\alpha}\zeta_\tau}^2) dxd\tau
        \leq C\int_0^t\Big(\norm{\nabla^{|\alpha|+1}(\phi,\psi,\zeta)}^2\\+\int_{\mathbb{D}}\abs{q}\abs{\nabla^{|\alpha|+1}(\phi,\psi,\zeta)}dx+\iint_{\mathbb{D}^{}}\frac{\nu(v)}{\M_*}\abs{\nabla^{|\alpha|+1}\tilde{\G}}^2dxdv\Big)d\tau+C(\delta_0+\chi)\int_{0}^{t}\mathcal{D}(\tau) d\tau,
    \end{multline}
    where we have used \cref{u1,Q1}.

\subsubsection{Dissipation estimates on $\p^{\alpha}\phi$ with $\alpha_0= 0$ and $1\le|\alpha|\le2$.}
Observing from the estimate \eqref{phial}, we need further estimates on $\p^{\alpha}\phi$. 
We start with \cref{sys-nonuni}$_2$, we have
$$
\p_i\big(p-\tilde{p}\big)=\frac{2}{3}\big(\rho\p_i\zeta+\T\p_i\phi\big)+\frac{2}{3}\big(\p_i\rhot\zeta+\phi\p_i\tilde{\T}\big).
$$
Thus, one can rewrite \cref{sys-nonuni}$_2$ and \cref{sys-nonuni}$_3$ as
\begin{align}\label{phi-alpha}
    \p^{\alpha}\psi_t+\frac{2}{3}\T\nabla\p^{\alpha}\phi=J_4^{\alpha},
\end{align}
where
\begin{multline*}
J_4^{\alpha}:=-\dv\p^{\alpha}\big(\rho \uf\otimes\uf-\rhot\ut_1^2\mathbb{I}_1\big)-\nabla\p^{\alpha}\big(p-\tilde{p}\big)+\frac{2}{3}\big(\rho\nabla\p^{\alpha}\zeta+\T\nabla\p^{\alpha}\phi\big)-\int v\otimes v\cdot \nabla\p^{\alpha}\tilde{\mathbf{G}} d v,\\
 \qquad\qquad-\frac{4}{3}\px\p^{\alpha}\left(\mu(\tilde{\theta}) \px\tilde{u}_{1}-\mu\left(\theta^{s_1}\right) \px u_{1}^{s_1}-\mu\left(\theta^{s_3}\right)\px u_{1}^{s_3}\right)-\px\p^{\alpha} Q_{1}-\frac{2}{3}\rho\nabla\p^{\alpha}\zeta,
\end{multline*}
Then multiplying \cref{phi-alpha} by $\p^{\alpha}\nabla\phi$ and integrating the resulting equation over $\mathbb{D}^{}\times[0,t]$, one has
\begin{multline*}
    \iint \p^{\alpha}\psi_1(t)\p^{\alpha}\phi(t)+\int_0^t\iint\p^{\alpha}\psi_{1}\px\dv \p^{\alpha}\psi
    dxd\tau+\int_0^t\iint\frac{2}{3}\T\abs{\nabla\p^{\alpha}\phi}^2 dxd\tau\\=\int_0^t\iint J_4^{\alpha}\cdot\p^{\alpha}\nabla\phi dxd\tau.
\end{multline*}
Then we arrive at 
\begin{multline}\label{phi0}
    \iint \p^{\alpha}\psi_1(t)\p^{\alpha}\phi(t) dx+\int_0^t\iint\frac{2}{3}\T\abs{\nabla\p^{\alpha}\phi}^2 dxd\tau\\
    \leq C\int_{0}^t\Big(\norm{\nabla^{|\alpha|+1}(\phi,\psi,\zeta)}_{L^2_x}^2+\iint_{\mathbb{D}^{}}\frac{\nu(v)}{\M_*}\abs{\nabla^{|\alpha|+1}\tilde{\G}}dxdv\Big)d\tau+C(\delta_0+\chi)\int_{0}^{t}\mathcal{D}(\tau) d\tau.
\end{multline}
Combining \cref{be1}, \cref{phial,phitau,phi0}, and the estimates on $\tilde{\G}$ \cref{nonfluid}, we obtain the following estimate,
\begin{align}\label{final-1}
   	&\norm{\big(\Phi, \tilde{\Psi}_1, \tilde{W}\big)(t, \cdot)}_{L^2}^2+\norm{(\phi, \varphi, \zeta)(t, \cdot)}_{H^2}^2+\sum_{0\leq \left|\alpha^{\prime}\right|\leq2} \iint \frac{\left|\partial^{\alpha^{\prime}} \tilde{\mathbf{G}}\right|^2}{\mathbf{M}_*} d v d x\nonumber\\
&+\int_0^t\sum_{|\alpha|=1}\left\|\p^{\alpha}\left({\Phi},\tilde{\Psi}_1, \tilde{W}\right)\right\|_{L_x^2}^2+\left\|\sqrt{\left|\px u^{s_1}\right|+\left|\px\Theta\right|}(\Phi, \tilde{\Psi}_1, \tilde{W})\right\|_{L_x^2}^2+\sum_{1 \leq\left|\alpha^{\prime}\right| \leq 2}\left\|\partial^{\alpha^{\prime}}\left(\phi_{\tau}, \varphi_\tau, \zeta_\tau\right)\right\|_{L_x^2}^2d\tau\nonumber\\
& +\int_{0}^{t}\sum_{1 \leq\left|\alpha^{\prime}\right| \leq 2}\left\|\partial^{\alpha^{\prime}}\left(\phi, \varphi, \zeta\right)\right\|_{L_x^2}^2+ \sum_{|\alpha'|=3}\left\|\partial^{\alpha^{\prime}}\left( \varphi, \zeta\right)\right\|_{L_x^2}^2+\sum_{1\leq |\beta'|\leq 2}\iint \frac{\nu(v)\left|\partial^{\beta^{\prime}} \tilde{\mathbf{G}}\right|^2}{\mathbf{M}_*} d v d x d\tau\nonumber\\
&\leq C(\mathcal{E}(0)^2+\delta^{\frac{1}{2}})+C\sum_{|\alpha|=3}\norm{\p^{\alpha}\phi}_{L^2_x}^2+C\sum_{|\beta'|=3}\int_0^t\iint \frac{\nu(v)\left|\partial^{\beta^{\prime}} \tilde{\mathbf{G}}\right|^2}{\mathbf{M}_*} d v d x d\tau.
\end{align}

\subsection{Highest-order estimate.}
The higher-order estimate \eqref{final-1} implies that we should further estimate the highest-order derivative term
\begin{align}
C\sum_{|\alpha|=3}\norm{\p^{\alpha}\phi}_{L^2_x}^2+C\sum_{|\beta'|=3}\int_0^t\iint \frac{\nu(v)|\partial^{\beta^{\prime}} \tilde{\mathbf{G}}|^2}{\mathbf{M}_*} d v d x d\tau,
\end{align}
for which we return back to the original system.

    
To obtain the estimate of the highest order of derivatives, we will use the equations \eqref{ft} for $\tilde{f}$, which is  
\begin{multline}\label{ftalpha}
(\partial^\alpha \tilde{f})_t+v\cdot\nabla(\partial^\alpha \tilde{f})=\partial^\alpha \LL_{\mathbf{M}} \tilde{\mathbf{G}}+\partial^\alpha Q(\tilde{\mathbf{G}}, \tilde{\mathbf{G}})+\partial^\alpha\Big[\left(\LL_{\mathbf{M}}-\LL_{\mathbf{M}^{s_1}}\right)\left(\mathbf{G}^{s_1}\right) \\
+(\LL_{\mathbf{M}}-\LL_{\mathbf{M}^{s_3}})(\mathbf{G}^{s_3})+2 Q(\tilde{\mathbf{G}}, \mathbf{G}^{s_1}+\mathbf{G}^{s_3})+2 Q(\mathbf{G}^{s_1}, \mathbf{G}^{s_3})\Big].
\end{multline}
Note that 
$$\frac{\partial^\alpha \tilde{f}}{\mathbf{M}_*}=\frac{\partial^\alpha\left(\mathbf{M}-\mathbf{M}^{s_1}-\mathbf{M}^{s_3}\right)}{\mathbf{M}_*}+\frac{\partial^\alpha \tilde{{\G}}}{\mathbf{M}_*},$$ 
then multiplying \cref{ftalpha} by $\frac{\partial^\alpha \tilde{f}}{\mathbf{M}_*}$, one has
\begin{align}
\label{111}
& \Big(\frac{|\partial^\alpha \tilde{f}|^2}{2 \mathbf{M}_*}\Big)_t-\frac{\partial^\alpha \tilde{\mathbf{G}}}{\mathbf{M}_*} \LL_{\mathbf{M}} \partial^\alpha \tilde{\mathbf{G}}=\frac{\partial^\alpha\left(\mathbf{M}-\mathbf{M}^{s_1}-\mathbf{M}^{s_3}\right)}{\mathbf{M}_*} \LL_{\mathbf{M}} \partial^\alpha \tilde{\mathbf{G}}-v\cdot\nabla\Big(\frac{ |\partial^\alpha \tilde{f}|^2}{2 \mathbf{M}_*}\Big) \nonumber\\
& +2 \frac{\partial^\alpha\left(\mathbf{M}-\mathbf{M}^{s_1}-\mathbf{M}^{s_3}\right)+\partial^\alpha \tilde{\mathbf{G}}}{\mathbf{M}_*}\bigg\{\frac{1}{2}\left(\LL_{\mathbf{M}}-\LL_{\mathbf{M}^{s_1}}\right) \partial^\alpha \mathbf{G}^{s_1}+\frac{1}{2}\left(\LL_{\mathbf{M}}-\LL_{\mathbf{M}^{s_3}}\right) \partial^\alpha \mathbf{G}^{s_3}\nonumber\\
& \nonumber+\sum_{|\beta|=1}^{|\alpha|-1}\left[Q\left(\partial^{\alpha-\beta} \mathbf{M}, \partial^\beta \tilde{\mathbf{G}}\right)+Q\left(\partial^{\alpha-\beta} \tilde{\mathbf{G}}, \partial^\beta \tilde{\mathbf{G}}\right)\right]+\frac{1}{2} Q\left(\partial^\alpha \mathbf{M}, \tilde{\mathbf{G}}\right)+Q\left(\partial^\alpha \tilde{\mathbf{G}}, \tilde{\mathbf{G}}\right) \\
& +\sum_{|\beta|=1}^{|\alpha|-1}\left[Q\left(\partial^{\alpha-\beta}\left(\mathbf{M}-\mathbf{M}^{s_1}\right), \partial^\beta \mathbf{G}^{s_1}\right)+Q\left(\partial^{\alpha-\beta}\left(\mathbf{M}-\mathbf{M}^{s_3}\right), \partial^\beta \mathbf{G}^{s_3}\right)\right] \nonumber\\
& +\frac{1}{2} Q\left(\partial^\alpha\left(\mathbf{M}-\mathbf{M}^{s_1}\right), \mathbf{G}^{s_1}\right)+\frac{1}{2} Q\left(\partial^\alpha\left(\mathbf{M}-\mathbf{M}^{s_3}\right), \mathbf{G}^{s_3}\right)+Q\left(\partial^\alpha \tilde{\mathbf{G}}, \mathbf{G}^{s_1}+\mathbf{G}^{s_3}\right) \nonumber\\
& +2 \sum_{|\beta|=1}^{|\alpha|-1}\left[Q\left(\partial^{\alpha-\beta} \tilde{\mathbf{G}}, \partial^\beta\left(\mathbf{G}^{s_1}+\mathbf{G}^{s_3}\right)\right)+Q\left(\partial^{\alpha-\beta} \mathbf{G}^{s_1}, \partial^\beta \mathbf{G}^{s_3}\right)\right] \nonumber\\
& +Q\left(\tilde{\mathbf{G}}, \partial^\alpha\left(\mathbf{G}^{s_1}+\mathbf{G}^{s_3}\right)\right)+Q\left(\partial^\alpha \mathbf{G}^{s_1}, \mathbf{G}^{s_3}\right)+Q\left(\mathbf{G}^{s_1}, \partial^\alpha \mathbf{G}^{s_3}\right)\bigg\}.
\end{align}
We first note that for $|\alpha'|=1$
\begin{align*}
\p^{\alpha'} \M & =\p^{\alpha'}\Big(\frac{\rho}{(\frac{4}{3} \pi  \theta)^{3 / 2}} \mathrm{e}^{-\frac{|v-u|^2}{2 R \theta}}\Big) \\
& =\M\Big(\frac{\p^{\alpha'}\rho}{\rho}-\frac{3 \p^{\alpha'}\theta}{2 \theta}+\frac{(v-u)^2 \p^{\alpha'}\theta}{ \frac{4}{3} \theta^2}+\sum_{i=1}^3 \frac{\p^{\alpha'}u_{i}\left(v_i-u_i\right)}{\frac{4}{3} \theta}\Big).
\end{align*}
Thus,
\begin{align*}
        &\notag\p^{\alpha'} \big(\M-\M^{s_1}-\M^{s_3}\big)   \\
 &\quad=\notag\M\Big(\frac{\p^{\alpha'}\phi}{\rho}-\frac{3 \p^{\alpha'}\zeta}{2 \theta}+\frac{(v-u)^2 \p^{\alpha'}\zeta}{ \frac{4}{3} \theta^2}+\sum_{i=1}^3 \frac{\p^{\alpha'}\varphi_{i}\left(v_i-u_i\right)}{\frac{4}{3} \theta}\Big)\\
 &\qquad+\notag\Big[\M\Big(\frac{\p^{\alpha'}\rho^{s_1}}{\rho}-\frac{3 \p^{\alpha'}\theta^{s_1}}{2 \theta}+\frac{(v-u)^2 \p^{\alpha'}\theta^{s_1}}{ \frac{4}{3} \theta^2}+\sum_{i=1}^3 \frac{\p^{\alpha'}u^{s_1}_{i}\left(v_i-u_i\right)}{\frac{4}{3} \theta}\Big)\\
 &\qquad-\notag\M^{s_1}\Big(\frac{\p^{\alpha'}\rho^{s_1}}{\rho^{s_1}}-\frac{3 \p^{\alpha'}\theta^{s_1}}{2 \theta^{s_1}}+\frac{(v-u^{s_1})^2 \p^{\alpha'}\theta^{s_1}}{ \frac{4}{3} (\theta^{s_1})^2}+\sum_{i=1}^3 \frac{\p^{\alpha'}u^{s_1}_{i}\left(v_i-u^{s_1}_i\right)}{\frac{4}{3} \theta^{s_1}}\Big)\Big]\\
 &\qquad+\notag\Big[\M\Big(\frac{\p^{\alpha'}\rho^{s_3}}{\rho}-\frac{3 \p^{\alpha'}\theta^{s_3}}{2 \theta}+\frac{(v-u)^2 \p^{\alpha'}\theta^{s_3}}{ \frac{4}{3} \theta^2}+\sum_{i=1}^3 \frac{\p^{\alpha'}u^{s_3}_{i}\left(v_i-u_i\right)}{\frac{4}{3} \theta}\Big)\\
 &\qquad\notag- \M^{s_3}\Big(\frac{\p^{\alpha'}\rho^{s_3}}{\rho^{s_3}}-\frac{3 \p^{\alpha'}\theta^{s_3}}{2 \theta^{s_3}}+\frac{(v-u^{s_3})^2 \p^{\alpha'}\theta^{s_3}}{ \frac{4}{3} (\theta^{s_3})^2}+\sum_{i=1}^3 \frac{\p^{\alpha'}u^{s_3}_{i}\left(v_i-u^{s_3}_i\right)}{\frac{4}{3} \theta^{s_3}}\Big)\Big],
\end{align*}
where the last two terms in the bracket are controllable small terms.
Then integrating \cref{111} over $\R^3\times\mathbb{D}\times[0,t]$, one has
\begin{align}
\label{ftau}\notag 
& \sum_{|\beta|=|\alpha|} \iint \frac{|\partial^\beta\tilde{f}|^2}{\mathbf{M}_*}(t, x, v) d v d x+\sum_{|\beta|=|\alpha|} \int_0^t \iint \frac{\nu(|v|)}{\mathbf{M}_*}\left|\partial^\beta\tilde{\mathbf{G}}\right|^2 d v d x d \tau \\
& \qquad\qquad\qquad\qquad\qquad\quad\leq C\left(\eta_0+\delta_0+\chi\right) \int_0^{t}\mathcal{D}(\tau)d\tau+C\left(\mathcal{E}(0)^2+\delta_0\right),
\end{align}
where we have used the wave interaction and \cref{fg} and \cref{linear}; cf. \cite{Li2017}. The small constant $\eta_0$ defined in \cref{linear} can be chosen as $\eta_0=O(1)\left(\delta_0+\chi\right)$.
Note that
$$
\left\|\p^{\alpha}\left(\phi, \varphi,\zeta\right)\right\|_{L^2_x}^2 \leq C \sum_{|\beta|=|\alpha|} \iint \frac{|\partial^\beta\tilde{f}|^2}{\mathbf{M}_*}(t, x, v) d v d x+C \delta_0,
$$
then we combine \cref{final-1,ftau} to finish the proof of \cref{ape}, and hence, the main result can be proved as in \cref{subproof}. 


\medskip
\noindent {\bf Acknowledgements.}
D.-Q. Deng was partially supported by the National Research Foundation of Korea (NRF) grant funded by the Korea government (MSIT) No. RS-2023-00210484 and No. RS-2023-00212304.

\bibliography{Bib_Ref}
\bibliographystyle{amsplain}

\end{document}